\documentclass[10pt]{amsart}
\usepackage{latexsym}   
\usepackage{amssymb}    
\usepackage{amsmath}    
\usepackage{amsthm}     
\usepackage[all]{xy}
\usepackage{hyperref}
\newcommand{\pa}{\partial}\newcommand{\al}{\alpha}

\newcommand{\Ga}{\Gamma}\newcommand{\del}{\delta}

\newcommand{\om}{\omega}
\newcommand{\Om}{\Omega}
\newcommand{\ti}{\tilde}

\renewcommand{\thefootnote}

\newtheorem{theorem}{Theorem}[section]

\theoremstyle{definition}
\newtheorem{definition}[theorem]{Definition}

\theoremstyle{remark}

\numberwithin{equation}{section}

\title[On the isometric deformation of surfaces via the B\"{a}cklund transformation] {
On the isometric deformation of surfaces via the B\"{a}cklund transformation}

\author[  Ion I. Dinc\u{a}]{Ion I. Dinc\u{a}}
\address{Department of Mathematical Methods and Models,
Faculty of Applied Sciences, University Politehnica of Bucharest
313 Spl.Independentei 060042 Bucharest, Romania}
 \email{ion.dinca@mathem.pub.ro}

\begin{document}


\begin{abstract}
In trying to generalize Bianchi's B\"{a}cklund transformation of quadrics to B\"{a}cklund transformations
of isometric deformations of other (classes of) surfaces, we investigate basic features of the
isometric deformation of surfaces via the B\"{a}cklund transformation with isometric correspondence
of leaves of a general nature (independent of the shape of the seed).
\end{abstract}

\maketitle

\tableofcontents \pagenumbering{arabic}

\section{Introduction}

The classical problem of finding the isometric deformations of surfaces (see Eisenhart \cite{E1})
was stated in 1859 by the French Academy of Sciences as

\begin{center}
{\it To find all surfaces applicable to a given one.}
\end{center}

Probably the most successful researcher of this problem is Bianchi, who in 1906 in \cite{B1}
solved the problem for quadrics by introducing the B\"{a}cklund transformation of surfaces
isometric to quadrics and the isometric correspondence provided by the Ivory affine transformation.
By 1909 Bianchi had a fairly complete treatment \cite{B2} and in 1917 he proved in \cite{B3} the
rigidity of the B\"{a}cklund transformation of isometric deformations of quadrics in the case of
auxiliary surface plane: the only B\"{a}cklund transformation with defining surface and auxiliary
surface plane appears as the singular B\"{a}cklund transformation of isometric deformations of
quadrics.

In \cite{D3} we have managed to improve Bianchi's result \cite{B3} to arbitrary auxiliary surface:
the only B\"{a}cklund transformation with defining surface is Bianchi's B\"{a}cklund transformation
of isometric deformations of quadrics.

The question remained to consider isometric deformation of surfaces via the B\"{a}cklund
transformation with isometric correspondence of leaves of a general
nature (independent of the shape of the seed) and without the assumption of defining surface;
it is the purpose of this note to address this question.

We shall consider the complexification
$$(\mathbb{C}^3,\langle\cdot,\cdot\rangle),\ \langle x,y\rangle :=x^Ty,\
|x|^2:=x^Tx,\ x,y\in\mathbb{C}^3$$
of the real $3$-dimensional Euclidean space; in this setting surfaces are $2$-dimensional objects of
$\mathbb{C}^3$ depending on two real or complex parameters.

{\it Isotropic} (null) vectors are those vectors of length $0$; since most vectors are not isotropic
we call a vector simply vector and we shall emphasize isotropic for isotropic vectors. The same
denomination will apply in other settings: for example we call quadric a non-degenerate quadric (a
quadric projectively equivalent to the complex unit sphere).

Consider Lie's viewpoint: one can replace a surface $x\subset\mathbb{C}^3$ with a $2$-dimensional
distribution of {\it contact elements} (pairs of points and planes passing through those points; the
classical geometers call them {\it facets}): the collection of its tangent planes (with the points
of tangency highlighted); thus a contact element is the infinitesimal version of a surface (the
integral element $(x,dx)|_{\mathrm{pt}}$ of the surface). Conversely, a $2$-dimensional distribution
of contact elements is not always the collection of the tangent planes of a surface (with the points
of tangency highlighted), but the condition that a $2$-dimensional distribution of contact elements
is integrable (that is it is the collection of the tangent planes of a leaf (sub-manifold)) does not
distinguish between the cases when this sub-manifold is a surface, curve or point, thus allowing the
collapsing of the leaf.

A $3$-dimensional distribution of contact elements is integrable if it is the collection of the
tangent planes of an $1$-dimensional family of leaves.

Two {\it rollable} (applicable or isometric) surfaces can be {\it rolled} (applied) one onto the
other such that at any instant they meet tangentially and with same differential at the tangency
point.

\begin{definition}
The rolling of two isometric surfaces $x_0,\ x\subset\mathbb{C}^3$ (that is $|dx_0|^2=|dx|^2$) is
the surface, curve or point $(R,t)\subset\mathbf{O}_3(\mathbb{C})\ltimes\mathbb{C}^3$ such that
$(x,dx)=(R,t)(x_0,dx_0):=(Rx_0+t,Rdx_0)$.
\end{definition}

The rolling introduces the flat connection form (it encodes the difference of the second fundamental
forms of $x_0,\ x$ and it being flat encodes the difference of the Gau\ss-Codazzi-Mainardi-
Peterson equations of $x_0,\ x$).

\begin{definition}
Consider an integrable $3$-dimensional distribution of contact elements $\mathcal{F}=(p,m)$ centered
at $p=p(u,v,w)$, with normal fields $m=m(u,v,w)$ and distributed along the surface $x_0=x_0(u,v)$.
If we roll $x_0$ on an isometric surface $x$ (that is $(x,dx)=(R,t)(x_0,dx_0):=(Rx_0+t,Rdx_0)$), then
the rolled distribution of contact elements is $(Rp+t,Rm)$ and is distributed along $x$; if it
remains integrable for any rolling, then the distribution is called {\it integrable rolling
distribution of contact elements}.
\end{definition}

Since by infinitesimal rolling in an arbitrary tangential infinitesimal direction $\del$ an initial
contact element $\mathcal{F}$ which is common tangent plane to two isometric surfaces is replaced
with an infinitesimally close contact element $\mathcal{F}'$ having in common with $\mathcal{F}$ the
direction $\del$, in the actual rolling problem we have contact elements centered on each other (the
symmetric {\it tangency configuration}) and contact elements centered on another one $\mathcal{F}$
reflect in $\mathcal{F}$; note that we assumed a finite law of a general nature as a consequence of
an infinitesimal law via discretization (the converse is clear); see Bianchi \cite{B3}.

Thus for a theory of isometric deformation of surfaces with the assumption above we are led to
consider, via rolling, certain $4$-dimensional distributions of contact elements centered on the
tangent planes of the considered surface $x_0$ and passing through the origin of the tangent planes
(each point of each tangent plane is the center of finitely many contact elements) and their rolling
counterparts on the isometric surface $x$.

After a thorough study of infinitesimal laws (in particular infinitesimal isometric deformations)
and their iterations the classical geometers were led to consider the B\"{a}cklund transformation (a
finite law of a general nature) as a consequence of an infinitesimal law via discretization.

The B\"{a}cklund transformation in the isometric deformation problem naturally appears by splitting
the $4$-dimensional distribution of contact elements above into an $1$-dimensional family of
$3$-dimensional integrable rolling distributions of contact elements, thus introducing a spectral
parameter $z$ (the B\"{a}cklund transformation is denoted $B_z$); each $3$-dimensional distribution
of contact elements is integrable (with {\it leaves} $x^1$) regardless of the shape of the
{\it seed} surface $x^0$.

The question appears if for the $4$-dimensional distribution of contact elements above being
integrable for any surface $x$ isometric to $x_0$ and with a $2$-dimensional family of leaves (in
general surfaces), then it splits into an $1$-dimensional family of $3$-dimensional integrable
rolling distributions of contact elements.

By imposing the initial collapsing ansatz of leaves of the $1$-dimensional family of $3$-dimensional
integrable rolling distributions of contact elements to be a $2$-dimensional family of curves (this
ansatz individuates the {\it defining surface} $x_0$) the $2$-dimensional family of curves naturally
splits into an $1$-dimensional family of {\it auxiliary surfaces} $x_z$ ($1$-dimensional families of
curves) such that the contact elements centered on these auxiliary surfaces form a $3$-dimensional
integrable rolling distribution of contact elements whose integrability is independent of the shape
of the seed surface $x^0$ isometric to the surface $x_0^0\subset x_0$. For $x_0$ quadric the
auxiliary surfaces $x_z$ are ((isotropic) singular) confocal quadrics doubly ruled by collapsed
leaves and with $z$ being the spectral parameter of the confocal family (this includes in a general
sense the (isotropic) planes of the pencil of (isotropic) planes through the axis of revolution for
$x_0$ quadric of revolution (excluding (pseudo-)spheres) (in this case the actual isotropic singular
quadric
of the confocal family is the two isotropic planes of this pencil) or (isotropic) planes of the
pencil of (isotropic) planes through an isotropic line (which is the actual isotropic singular
quadric of its confocal family) for $x_0$ Darboux quadric with tangency of order $3$ with the circle
at infinity $C(\infty)$).

This collapsing ansatz allows us to simplify the denomination {\it B\"{a}cklund transformation of
surfaces isometric to quadrics} to {\it B\"{a}cklund transformation of quadrics}; Bianchi's
B\"{a}cklund transformation of quadrics is just the metric-projective generalization of Lie's point
of view on the B\"{a}cklund transformation of the pseudo-sphere (one replaces 'pseudo-sphere' with
'quadric' and 'circle' with 'conic').

Note however that according to our note \cite{D3} there are no other instances of B\"{a}cklund
transformations with defining surface except Bianchi's B\"{a}cklund transformation of quadrics.

Bianchi considered the most general form of a B\"{a}cklund transformation as the focal surfaces (one
transform of the other) of a {\it Weingarten} congruence (congruence upon whose two focal surfaces
the asymptotic directions correspond; equivalently the second fundamental forms are proportional).
Note that although the correspondence provided by the Weingarten congruence does not give the
applicability (isometric) correspondence, the B\"{a}cklund transformation is the tool best suited to
attack the isometric deformation problem via transformation, since it provides correspondence of the
characteristics of the isometric deformation problem (according to Darboux these are the asymptotic
directions), it is directly linked to the infinitesimal isometric deformation problem (Darboux
proved that infinitesimal isometric deformations generate Weingarten congruences and Guichard proved
the converse: there is an infinitesimal isometric deformation of a focal surface of a Weingarten
congruence in the direction normal to the other focal surface; see Darboux (\cite{D1},\S\ 883-\S\
924)) and it admits a version of the Bianchi Permutability Theorem for its second iteration.

\subsection{Preliminaries and results from \cite{D3}}\noindent

\noindent

Let $(u,v)\in D$ with $D$ domain of $\mathbb{R}^2$ or $\mathbb{C}^2$ and $x:D\mapsto\mathbb{C}^3$ be
a surface.

For $\om_1,\om_2\ \mathbb{C}^3$-valued $1$-forms on $D$ and $a,b\in\mathbb{C}^3$ we have
\begin{eqnarray}\label{eq:fund}
a^T\om_1\wedge b^T\om_2=((a\times b)\times\om_1+b^T\om_1a)^T\wedge\om_2=
(a\times b)^T(\om_1\times\wedge\om_2)+b^T\om_1\wedge a^T\om_2;\nonumber\\
\mathrm{in\ particular}\ a^T\om\wedge b^T\om=\frac{1}{2}(a\times b)^T(\om\times\wedge\om).
\end{eqnarray}

Since both $\times$ and $\wedge$ are skew-symmetric, we have $2\om_1\times\wedge\om_2=
\om_1\times\om_2+\om_2\times\om_1=2\om_2\times\wedge\om_1$.

Consider the scalar product $\langle\cdot,\cdot\rangle$ on $\mathbf{M}_3(\mathbb{C}):\
\langle X,Y\rangle :=\frac{1}{2}\mathrm{tr}(X^TY)$. We have the isometry

$$\al:\mathbb{C}^3\mapsto\mathbf{o}_3(\mathbb{C}),\
\al\left(\begin{bmatrix}x^1\\x^2\\x^3\end{bmatrix}\right)
=\begin{bmatrix}0&-x^3&x^2\\x^3&0&-x^1\\-x^2&x^1&0\end{bmatrix},\
x^Ty=\langle\al(x),\al(y)\rangle=
\frac{1}{2}\mathrm{tr}(\al(x)^T\al(y)),$$
$$\al(x\times y)=[\al(x),\al(y)]=\al(\al(x)y)=yx^T-xy^T,\ \al(Rx)=R\al(x)R^{-1},\
x,y\in\mathbb{C}^3,\ R\in\mathbf{O}_3(\mathbb{C}).$$

Let $x\subset\mathbb{C}^3$ be a surface non-rigidly applicable (isometric) to a surface
$x_0\subset\mathbb{C}^3$:

\begin{eqnarray}\label{eq:x0x}
(x,dx)=(R,t)(x_0,dx_0):=(Rx_0+t,Rdx_0),
\end{eqnarray}
where $(R,t)$ is a sub-manifold in $\mathbf{O}_3(\mathbb{C})\ltimes\mathbb{C}^3$ (in general
surface, but it is a curve if $x_0,\ x$ are ruled and the rulings correspond under isometry or a
point if $x_0,\ x$ differ by a rigid motion). The sub-manifold $R$ gives the rolling of $x_0$ on $x$,
that is if we rigidly roll $x_0$ on $x$ such that points corresponding under the isometry will have
the same differentials, $R$ will dictate the rotation of $x_0$; the translation $t$ will satisfy
$dt=-dRx_0$.

For $(u,v)$ parametrization on $x_0,\ x$ and outside the locus of isotropic (degenerate) induced
metric of $x_0,\ x$ we have $N_0:=\frac{\pa_ux_0\times\pa_vx_0}{|\pa_ux_0\times\pa_vx_0|},\
N:=\frac{\pa_ux\times\pa_vx}{|\pa_ux\times\pa_vx|}$ respectively positively oriented unit normal
fields of $x_0,\ x$ and $R$ is determined by $R=[\pa_ux\ \ \pa_vx\ \ N][\pa_ux_0\ \ \pa_vx_0\ \
\det(R)N_0]^{-1}$; we take $R$ with $\det(R)=1$; thus the rotation of the rolling with the other
face of $x_0$ (or on the other face of $x$) is $R':=R(I_3-2N_0N_0^T)=(I_3-2NN^T)R,\ \det(R')=-1$.

With $\om:=\al^{-1}(R^{-1}dR)=N_0\times R^{-1}dRN_0$ we have

\begin{eqnarray}\label{eq:om}
d\om+\frac{1}{2}\om\times\wedge\om=0,\ \om\times\wedge dx_0=0,\ (\om)^\bot=0.
\end{eqnarray}

Here we shall recall some results from \cite{D3}.

Consider a surface $x_0=x_0(u,v)\subset\mathbb{C}^3$ with unit normal field $N_0=N_0(u,v)$.

Consider a $3$-dimensional distribution of contact elements with the symmetry of the tangency
configuration, that is the contact elements are centered at $x_0+V,\ V=V(u,v,w),\ N_0^TV=0,\
du\wedge dv\wedge dw\neq 0$ and have non-isotropic normal fields $m=V\times N_0+\mathbf{m}N_0,\
\mathbf{m}=\mathbf{m}(u,v,w)\subset\mathbb{C}$.

With $\ti d\cdot:=\pa_u\cdot du+\pa_v\cdot dv+\pa_w\cdot dw=d\cdot+\pa_w\cdot dw$ if the
distribution of contact elements is integrable and the rolled distribution remains integrable if we
roll $x_0$ on an isometric surface $x,\ (x,dx)=(R,t)(x_0,dx_0)$ (that is we replace $x_0,\ V,\ m$
with $x,\ RV,\ Rm$), then along the leaves we have
$$0=(Rm)^T\ti d(RV+x)=m^T(\om\times V+d(V+x_0)+\pa_wVdw),$$
or, assuming $N_0^T(\pa_wV\times V)\neq 0$,

\begin{eqnarray}\label{eq:dw}
dw=\frac{N_0^T[V\times d(V+x_0)]}{N_0^T(\pa_wV\times V)}+
\mathbf{m}\frac{V^T(\om\times N_0+dN_0)}{N_0^T(\pa_wV\times V)}.
\end{eqnarray}
By applying the compatibility condition $\ti d$ to (\ref{eq:dw}) and using the equation itself we get
the integrability condition

\begin{eqnarray}\label{eq:dm}
V\times\pa_wV\neq 0,\ \frac{2[\pa_wV\times d(V+x_0)]^T\wedge(V\times dx_0)}
{(\pa_wV\times V)^T(dx_0\times\wedge dx_0)}+(\mathbf{m}^2+|V|^2)K=0,\nonumber\\
d\mathbf{m}=-\pa_w\mathbf{m}\frac{N_0^T[V\times d(V+x_0)]}{N_0^T(\pa_wV\times V)}+
\mathbf{m}\frac{N_0^T(\pa_wV\times dV)}{N_0^T(\pa_wV\times V)},\ du\wedge dv\wedge dw\neq 0.
\end{eqnarray}
Assume $\mathbf{m}\neq 0$; excluding the case $x_0$ developable we can take $\mathbf{m}^2$ from
the
first equation of (\ref{eq:dm}) and replace it into the second one; applying the compatibility
condition $d$ to this equation and using the equation itself we get a relation, which coupled to
differentiating the first equation of (\ref{eq:dm}) with respect to $w$ gives

\begin{eqnarray}\label{eq:pawV}
\frac{2(\pa_wV\times dV)^T\wedge(\pa_wV\times dx_0)}{(\pa_wV\times V)^T(dx_0\times\wedge dx_0)}
+(\mathbf{m}\pa_w\mathbf{m}+V^T\pa_wV)K=0,\ du\wedge dv\wedge dw\neq 0
\end{eqnarray}
and (\ref{eq:dm}) becomes
\begin{eqnarray}\label{eq:intrdst}
\pa_w\frac{N_0^T[\pa_wV\times d(V+x_0)]}{N_0^T(\pa_wV\times V)}\wedge N_0^T(V\times dx_0)-
\frac{N_0^T(\pa_wV\times dV)}{N_0^T(\pa_wV\times V)}\wedge N_0^T(\pa_wV\times dx_0)=0,\nonumber\\
\frac{1}{2}d(2\frac{[\pa_wV\times d(V+x_0)]^T\wedge(V\times dx_0)}
{K(\pa_wV\times V)^T(dx_0\times\wedge dx_0)}+|V|^2)=
-(2\frac{[\pa_wV\times dV]^T\wedge(\pa_wV\times dx_0)}{K(\pa_wV\times V)^T(dx_0\times\wedge dx_0)}
+V^T\pa_wV)\nonumber\\
\frac{N_0^T[V\times d(V+x_0)]}{N_0^T(\pa_wV\times V)}+
(2\frac{[\pa_wV\times d(V+x_0)]^T\wedge(V\times dx_0)}{K(\pa_wV\times V)^T(dx_0\times\wedge dx_0)}
+|V|^2)\frac{N_0^T(\pa_wV\times dV)}{N_0^T(\pa_wV\times V)},\nonumber\\
du\wedge dv\wedge dw\neq 0;
\end{eqnarray}
applying $\pa_w$ to the second equation of (\ref{eq:intrdst}) and using the first one we get another
second order equation in $V$.

From (\ref{eq:pawV}) we get the Weingarten congruence property of $x$ and any leaf of the
$3$-dimensional integrable rolling distribution of contact elements with the symmetry of the tangency
configuration (and thus B\"{a}cklund transformation according to Bianchi's definition).

We consider now the question whether the leaves of a $3$-dimensional integrable rolling
distribution of contact elements are isometric deformations of surfaces; we exclude the case $m$
isotropic because all leaves are isotropic developables with degenerate $2$-dimensional metric.

Consider a surface $x_0=x_0(u,v)\subset\mathbb{C}^3$ with unit normal field $N_0=N_0(u,v)$.

Consider a $3$-dimensional distribution of contact elements distributed along $x_0$, that is the
contact elements are centered at $x_0+V,\ V=V(u,v,w),\ du\wedge dv\wedge dw\neq 0$ and have normal
fields $m=m(u,v,w)\subset\mathbb{C}^3$.

With $\ti d\cdot:=\pa_u\cdot du+\pa_v\cdot dv+\pa_w\cdot dw=d\cdot+\pa_w\cdot dw$ if the distribution
of contact elements is integrable and the rolled distribution remains integrable if we roll $x_0$ on
an isometric surface $x,\ (x,dx)=(R,t)(x_0,dx_0)$ (that is we replace $x_0,\ V,\ m$ with $x,\ RV,\
Rm$), then along the leaves we have

\begin{eqnarray}\label{eq:int}
0=(Rm)^T\ti d(RV+x)=m^T(\om\times V+d(V+x_0)+\pa_wVdw).
\end{eqnarray}
Since we shall not need the integrability condition of this general integrable rolling distribution
of contact elements, we shall not derive it.

The leaves for a particular isometric deformation $x_0$ are isometric to the same surface or to an
$1$-dimensional family of surfaces.

Since we already have a correspondence between the contact elements of leaves given by rolling, it is
natural to require that the isometric correspondence is independent of the shape of $x$.

However, for $(m\times V)\times N_0\neq 0$ the distribution of contact elements into leaves changes
with the shape of $x$; thus all leaves of all rolled distributions of contact elements must be
isometric to the same surface (and we have a submersion from the $3$-dimensional integrable rolling
distribution of contact elements to the distribution of tangent planes of the fixed surface), or
$(m\times V)\times N_0=0$ and we have only a $2$-dimensional integrable rolling distribution of
contact elements with the leaf having metric independent of the shape of $x$.

Consider the general case $m^T\pa_wV(m\times V)\times N_0\neq 0$ and all leaves of all rolled
distributions of contact elements are isometric to the same surface $y=y(u_1,v_1)=y(u_1(u,v,w),
v_1(u,v,w))$; if for a particular isometric deformation of $x_0$ one knows the isometric
correspondence of all leaves to the surface $y$, then one finds the isometric correspondence to $y$
of all rolled leaves (including the case of leaves with degenerate metric) by composing the rolling
of the particular leaves on $y$ with the inverse of the rolling of $x_0$ on $x$. For $u_1,v_1=$ct
from the three independent variables $u,v,w$ only one remains independent, thus giving the
submersion from the integrable rolling distribution of contact elements to the distribution of
tangent planes of $y$ (equivalently we count each contact element $(y,dy)$ with simple $\infty$
multiplicity).

The function $w=w(u,v,c)$ is given by the integration of (\ref{eq:int}); thus a-priori $w$ depends
also on $\om$ and we have

\begin{eqnarray}\label{eq:iso}
|(I_3-\frac{\pa_wVm^T}{m^T\pa_wV})[\om\times V+d(V+x_0)]|^2=
|dy-\frac{m^T[\om\times V+d(V+x_0)]}{m^T\pa_wV}\pa_wy|^2,\nonumber\\
\forall\om\ \mathrm{satisfying}\ (\ref{eq:om}).
\end{eqnarray}
The leaves are isometric to different regions of $y$ because the constant $c$ in $w=w(u,v,c)$ changes
for $\om$ fixed and for different $\om\ w$ changes; however in order to be determined by
(\ref{eq:int}), $w$ is not allowed to be linked to $\om$ by any other relation, either functional
(as a-priori (\ref{eq:iso}) is) or differential; thus in (\ref{eq:iso}) $\om$ cancels independently
of $w$ and outside $w$ we can replace $\om$ with any other $\hat\om$ satisfying (\ref{eq:om}).

After some manipulation of these premises we get the tangency configuration $V^TN_0=0$ and

\begin{eqnarray}\label{eq:tidy}
\ti dy=R_1(I_3-\frac{N_0m^T}{m^TN_0})\ti d(V+x_0),\ R_1\subset\mathbf{O}_3(\mathbb{C}),\
du\wedge dv\wedge dw\neq 0.
\end{eqnarray}
Note that since along the leaves we have $m^T\ti d(V+x_0)=0,\ R_1$ is the rotation of the rolling of
the leaves on $y$; if we replace $x_0$ with an isometric surface $x=Rx_0+t$, then $R_1$ is replaced
with $R_1R^{-1}$.

Imposing the compatibility condition $R_1^{-1}\ti d$ on (\ref{eq:tidy}) we get

\begin{eqnarray}\label{eq:R1}
0=[R_1^{-1}\ti dR_1(I_3-\frac{N_0m^T}{m^TN_0})-\ti d\frac{N_0m^T}{m^TN_0}]\wedge\ti d(V+x_0);
\end{eqnarray}
with $R^{-1}\ti dR_1=:
\Om_1du+\Om_2dv+\Om_3dw$ this constitutes a linear system of $9$ equations in the $9$ entries of
$\Om_j,\ j=1,2,3$ with the rank of the matrix of the system being $6$, so the rank of the augmented
matrix of the system must be also $6$ and the solution $R_1^{-1}\ti dR_1$ must also satisfy the
compatibility condition

\begin{eqnarray}\label{eq:dR1}
\ti d(R_1^{-1}\ti dR_1)+\frac{1}{2}[R_1^{-1}\ti dR_1,\wedge
R_1^{-1}\ti dR_1]=0;
\end{eqnarray}
these are the necessary and sufficient conditions on the $3$-dimensional
tangential integrable rolling distribution of contact elements in order to obtain isometric
correspondence of leaves of a general nature.

From (\ref{eq:tidy}) $R_1$ is determined modulo the rotation of the rolling of $y$ on an isometric
surface, so from
the three constants introduced in the solution of (\ref{eq:R1}) (which are actually functions of
three variables $u,v,w$) satisfying (\ref{eq:dR1}) we must get dependence on two functions of a
variable related to the isometric deformation problem.

By applying, if necessary, a change of variables $w=w(\ti w, u,v)$ we have
$N_0^T[d(V+x_0)\times\wedge d(V+x_0)]\neq 0$ and the above considered linear system is consistent
for
\begin{eqnarray}\label{eq:cons}
d(V+x_0)^T(I_3-\frac{mN_0^T}{m^TN_0})\odot[d(\frac{N_0m^T}{m^TN_0})\pa_wV-
\pa_w(\frac{N_0m^T}{m^TN_0})d(V+x_0)+\nonumber\\
2N_0^T[\pa_wV\times d(V+x_0)]\frac{d(\frac{N_0m^T}{m^TN_0})\wedge
d(V+x_0)}{N_0^T[d(V+x_0)\times\wedge d(V+x_0)]}]=0,\ du\wedge dv\wedge dw\neq 0;
\end{eqnarray}
if we further assume the symmetric tangency configuration $m^TV=0$, then by using (\ref{eq:dm})
(\ref{eq:cons}) is satisfied.

\subsection{The main program}

Our purpose is to generalize Bianchi's B\"{a}cklund transformation of quadrics to B\"{a}cklund
transformation of isometric deformations of other (classes of) surfaces.

The question remained to consider isometric deformations of surfaces via the B\"{a}cklund
transformation with isometric correspondence of leaves of a general
nature (independent of the shape of the seed) and without the assumption of defining surface.

Similarly to the quadrics case, we should get

$\bullet$ (Existence of B\"{a}cklund transformation) necessary and sufficient conditions on the
metric of the
seed surface in order to get B\"{a}cklund transformation with isometric correspondence of leaves of
a general nature (independent of the shape of the seed). Once these conditions are satisfied, a
certain differential system subjacent to such B\"{a}cklund transformations introduces a constant
(denoted $z$ in $B_z$) and a second constant is introduced by finding the leaves of the integrable
rolling distribution of contact elements.

If (classes of) surfaces other than quadrics and B\"{a}cklund transformation with defining
surface are found to satisfy these requirements, then we continue the program with

$\bullet$ (Inversion of the B\"{a}cklund transformation) the metric of the leaves satisfies the
necessary and sufficient conditions in order to get B\"{a}cklund transformation with isometric
correspondence of leaves of a general nature (independent of the shape of the seed); thus the seed
can be considered leaf and the leaf seed.

$\bullet$ (Bianchi Permutability Theorem) If $x^1=B_{z_1}(x^0),\ x^2=B_{z_2}(x^0)$, then one can
find only by algebraic computations a surface $B_{z_1}(x^2)=x^3=B_{z_2}(x^1)$; thus
$B_{z_2}\circ B_{z_1}=B_{z_1}\circ B_{z_2}$ and
once all B transforms of the seed $x^0$ are found, the B transformation can be iterated using only
algebraic computations.

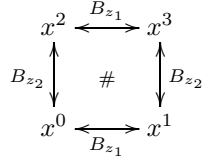
\begin{center}
$\xymatrix{\ar@{}[dr]|{\#}x^2\ar@{<->}[d]_{B_{z_2}}\ar@{<->}[r]^{B_{z_1}}&x^3\ar@{<->}[d]^{B_{z_2}}\\
x^0\ar@{<->}[r]_{B_{z_1}}&x^1}$
\end{center}

$\bullet$ (Existence of $3$-M\"{o}bius moving configurations) If $x^1=B_{z_1}(x^0),\
x^2=B_{z_2}(x^0),\ x^4=B_{z_3}(x^0)$ and by use of the Bianchi Permutability
Theorem one finds $B_{z_3}(x^2)=x^6=B_{z_2}(x^4),\ B_{z_1}(x^4)=x^5=B_{z_3}(x^1),\
B_{z_2}(x^1)=x^3=B_{z_1}(x^2),\ B_{z_3}(x^3)={x'}^7=B_{z_2}(x^5),\
B_{z_1}(x^6)={x''}^7=B_{z_3}(x^3),\\ B_{z_2}(x^5)={x'''}^7=B_{z_1}(x^6)$, then
${x'}^7={x''}^7={x'''}^7=:x^7$; thus once all $B$
transforms of the seed $x^0$ are found, the $B$ transformation can be further iterated using only
algebraic computations.

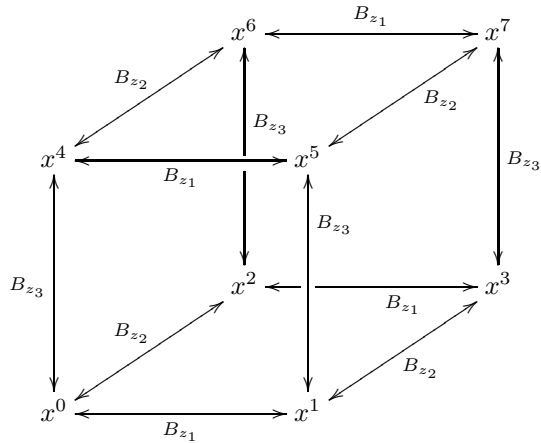
\begin{center}
$\xymatrix@!0{&&&x^6\ar@{<->}[rrrr]^{B_{z_1}}&&&&x^7\\
&&&&&&&\\
x^4\ar@{<->}[uurrr]^{B_{z_2}} \ar@{<->}[rrrr]_{B_{z_1}}&&&&
x^5\ar@{<->}[uurrr]_>>>>>>>>>{B_{z_2}}&&&\\
&&&&&&&\\
&&&x^2\ar@{<->}'[r][rrrr]_{B_{z_1}}
\ar@{<->}'[uu][uuuu]_<<<<<{B_{z_3}}&&&&
x^3\ar@{<->}[uuuu]_{B_{z_3}}\\
&&&&&&&\\
x^0\ar@{<->}[rrrr]_{B_{z_1}}\ar@{<->}[uurrr]^{B_{z_2}}
\ar@{<->}[uuuu]^{B_{z_3}}&&&& x^1\ar@{<->}[uurrr]_{B_{z_2}}
\ar@{<->}[uuuu]_<<<<<<<<<<<<<<<<<<<<<<{B_{z_3}}&&&}$
\end{center}

Out of this program we were able to partially fulfill the requirements of the first item: we found
necessary and sufficient conditions on $V,\mathbf{m}$ (and their derivatives) and the metric of
the seed surface $x_0$ which using
(\ref{eq:dm}) and adjoined to (\ref{eq:intrdst}) constitutes a differential system on $V$ (here
compatibility conditions may impose conditions on the metric of the seed surface).

\section{The B\"{a}cklund transformation with isometric correspondence of leaves of a general
nature (independent of the shape of the seed)}

From (\ref{eq:R1}) with $m=V\times N_0+\mathbf{m}N_0,\ R_1^{-1}\ti dR_1=:\al(\ti\om_udu+\ti\om_vdv
+\ti\om_wdw)=\al(\ti\om+\ti\om_wdw)$ and using (\ref{eq:fund}) and (\ref{eq:dm})
(in particular $-\frac{1}{\mathbf{m}}\frac{N_0^T(V\times dV)}{N_0^T(\pa_wV\times V)}\wedge
N_0^T[\pa_wV\times d(V+x_0)]
-\frac{|V|^2}{\mathbf{m}}\frac{1}{2}N_0^T(dN_0\times\wedge dN_0)=
\frac{1}{2\mathbf{m}}N_0^T[d(V+x_0)\times\wedge d(V+x_0)]
+\frac{\mathbf{m}}{2}N_0^T(dN_0\times\wedge dN_0)$) we get

\begin{eqnarray}\label{eq:tiom}
N_0^T[\ti\om\times\wedge d(V+x_0)]-dN_0^T\wedge dV
+\frac{1}{2\mathbf{m}}N_0^T[d(V+x_0)\times\wedge d(V+x_0)]
+\frac{\mathbf{m}}{2}N_0^T(dN_0\times\wedge dN_0)=0,\nonumber\\
-\ti\om^TN_0\wedge d(V+x_0)\times N_0+\ti\om\times N_0\wedge\frac{N_0^T[V\times d(V+x_0)]}{\mathbf{m}}
-\frac{1}{2}N_0^T(dN_0\times\wedge dN_0)N_0\times V\nonumber\\
+dN_0\wedge\frac{N_0^T[V\times d(V+x_0)]}{\mathbf{m}}=0,\nonumber\\
N_0^T[\ti\om_w\times d(V+x_0)]-N_0^T(\ti\om\times\pa_wV)+dN_0^T\pa_wV+
\frac{N_0^T[\pa_wV\times d(V+x_0)]}{\mathbf{m}}=0,\nonumber\\
-\ti\om_w^TN_0d(V+x_0)\times N_0+\ti\om_w\times N_0\frac{N_0^T[V\times d(V+x_0)]}{\mathbf{m}}
-\ti\om^TN_0N_0\times\pa_wV\nonumber\\
+\frac{\ti\om\times N_0+dN_0}{\mathbf{m}}
N_0^T(\pa_wV\times V)=0.
\end{eqnarray}
With $\mathcal{M}:=-\frac{1}{2\mathbf{m}}
+\frac{dN_0^T\wedge dV-\frac{\mathbf{m}}{2}N_0^T(dN_0\times\wedge dN_0)}
{N_0^T[d(V+x_0)\times\wedge d(V+x_0)]},\ \mathcal{U}:=\pa_u(V+x_0),\ \mathcal{V}:=\pa_v(V+x_0)$
we have from the first equation of (\ref{eq:tiom})
$$\ti\om_u=(c_1+\mathcal{M})(N_0\times\mathcal{U})\times N_0+c_2(N_0\times\mathcal{V})\times N_0
+c_3N_0,$$
$$\ti\om_v=c_4(N_0\times\mathcal{U})\times N_0-(c_1-\mathcal{M})(N_0\times\mathcal{V})\times N_0
+c_5N_0;$$
and multiplying the second equation of (\ref{eq:tiom}) on the left with $\mathcal{U}^T,\
\mathcal{V}^T$ we get
$$c_3=c_2\frac{N_0^T(V\times\mathcal{V})}{\mathbf{m}}+
(c_1-\mathcal{M})\frac{N_0^T(V\times\mathcal{U})}{\mathbf{m}}$$
$$+\mathcal{U}^T\frac{-N_0^T(\pa_uN_0\times\pa_vN_0)N_0\times V
+\pa_uN_0\frac{N_0^T(V\times\mathcal{V})}{\mathbf{m}}-
\pa_vN_0\frac{N_0^T(V\times\mathcal{U})}{\mathbf{m}}}{N_0^T(\mathcal{U}\times\mathcal{V})},$$
$$c_5=-(c_1+\mathcal{M})\frac{N_0^T(V\times\mathcal{V})}{\mathbf{m}}+
c_4\frac{N_0^T(V\times\mathcal{U})}{\mathbf{m}}$$
$$+\mathcal{V}^T\frac{-N_0^T(\pa_uN_0\times\pa_vN_0)N_0\times V
+\pa_uN_0\frac{N_0^T(V\times\mathcal{V})}{\mathbf{m}}-
\pa_vN_0\frac{N_0^T(V\times\mathcal{U})}{\mathbf{m}}}{N_0^T(\mathcal{U}\times\mathcal{V})}.$$
With $\ti\om_w=:c_6(N_0\times\mathcal{U})\times N_0+c_7(N_0\times\mathcal{V})\times N_0+c_8N_0$
from the third equation of (\ref{eq:tiom}) we get
$$c_6=-\frac{c_4N_0^T(\pa_wV\times\mathcal{U})-(c_1-\mathcal{M})N_0^T(\pa_wV\times\mathcal{V})
+\pa_vN_0^T\pa_wV+\frac{N_0^T(\pa_wV\times\mathcal{V})}
{\mathbf{m}}}{N_0^T(\mathcal{U}\times\mathcal{V})},$$
$$c_7=\frac{(c_1+\mathcal{M})N_0^T(\pa_wV\times\mathcal{U})+c_2N_0^T(\pa_wV\times\mathcal{V})
+\pa_uN_0^T\pa_wV+\frac{N_0^T(\pa_wV\times\mathcal{U})}
{\mathbf{m}}}{N_0^T(\mathcal{U}\times\mathcal{V})}.$$
Multiplying the last equation of (\ref{eq:tiom}) on the left with $\mathcal{U}^T,\
\mathcal{V}^T$ we get four scalar equations; from one of them we get
$$c_8=-\frac{c_1-\mathcal{M}}{\mathbf{m}}N_0^T(\pa_wV\times V)
+\frac{c_7}{\mathbf{m}}N_0^T(V\times\mathcal{V})-\frac{c_5N_0^T(\pa_wV\times\mathcal{U})
-\frac{\pa_vN_0^T\mathcal{U}}{\mathbf{m}}N_0^T(\pa_wV\times V)}
{N_0^T(\mathcal{U}\times\mathcal{V})}$$
and the remaining three are identically satisfied.

Thus

\begin{eqnarray}\label{eq:tiomuvw}
\ti\om_u=c_1((N_0\times\mathcal{U})\times N_0+\frac{N_0^T(V\times\mathcal{U})}{\mathbf{m}}N_0)
+c_2((N_0\times\mathcal{V})\times N_0
+\frac{N_0^T(V\times\mathcal{V})}{\mathbf{m}}N_0)\nonumber\\
+\mathcal{M}((N_0\times\mathcal{U})\times N_0-\frac{N_0^T(V\times\mathcal{U})}{\mathbf{m}}N_0)
\nonumber\\
+\mathcal{U}^T\frac{-N_0^T(\pa_uN_0\times\pa_vN_0)N_0\times V
+\pa_uN_0\frac{N_0^T(V\times\mathcal{V})}{\mathbf{m}}-
\pa_vN_0\frac{N_0^T(V\times\mathcal{U})}{\mathbf{m}}}{N_0^T(\mathcal{U}\times\mathcal{V})}N_0,
\nonumber\\
\ti\om_v=c_4((N_0\times\mathcal{U})\times N_0+\frac{N_0^T(V\times\mathcal{U})}{\mathbf{m}}N_0)
-c_1((N_0\times\mathcal{V})\times N_0
+\frac{N_0^T(V\times\mathcal{V})}{\mathbf{m}}N_0)\nonumber\\
+\mathcal{M}((N_0\times\mathcal{V})\times N_0
-\frac{N_0^T(V\times\mathcal{V})}{\mathbf{m}}N_0)\nonumber\\
+\mathcal{V}^T\frac{-N_0^T(\pa_uN_0\times\pa_vN_0)N_0\times V
+\pa_uN_0\frac{N_0^T(V\times\mathcal{V})}{\mathbf{m}}-
\pa_vN_0\frac{N_0^T(V\times\mathcal{U})}{\mathbf{m}}}{N_0^T(\mathcal{U}\times\mathcal{V})}N_0,
\nonumber\\
\ti\om_w=c_1(\frac{N_0^T(\pa_wV\times\mathcal{V})}{N_0^T(\mathcal{U}\times\mathcal{V})}
((N_0\times\mathcal{U})\times N_0+\frac{N_0^T(V\times\mathcal{U})}{\mathbf{m}}N_0)\nonumber\\
+\frac{N_0^T(\pa_wV\times\mathcal{U})}{N_0^T(\mathcal{U}\times\mathcal{V})}
((N_0\times\mathcal{V})\times N_0+\frac{N_0^T(V\times\mathcal{V})}{\mathbf{m}}N_0))\nonumber\\
+c_2\frac{N_0^T(\pa_wV\times\mathcal{V})}{N_0^T(\mathcal{U}\times\mathcal{V})}
((N_0\times\mathcal{V})\times N_0+\frac{N_0^T(V\times\mathcal{V})}{\mathbf{m}}N_0)\nonumber\\
-c_4\frac{N_0^T(\pa_wV\times\mathcal{U})}{N_0^T(\mathcal{U}\times\mathcal{V})}
((N_0\times\mathcal{U})\times N_0+\frac{N_0^T(V\times\mathcal{U})}{\mathbf{m}}N_0)\nonumber\\
-\frac{\mathcal{M}N_0^T(\pa_wV\times\mathcal{V})
+\pa_vN_0^T\pa_wV+\frac{N_0^T(\pa_wV\times\mathcal{V})}
{\mathbf{m}}}{N_0^T(\mathcal{U}\times\mathcal{V})}(N_0\times\mathcal{U})\times N_0\nonumber\\
+\frac{\mathcal{M}N_0^T(\pa_wV\times\mathcal{U})+\pa_uN_0^T\pa_wV
+\frac{N_0^T(\pa_wV\times\mathcal{U})}
{\mathbf{m}}}{N_0^T(\mathcal{U}\times\mathcal{V})}(N_0\times\mathcal{V})\times N_0\nonumber\\
+(\mathcal{M}\frac{N_0^T(\pa_wV\times V)}{\mathbf{m}}+
\frac{\pa_uN_0^T\pa_wV}{N_0^T(\mathcal{U}\times\mathcal{V})}\frac{N_0^T(V\times\mathcal{V})}
{\mathbf{m}}-\frac{\pa_vN_0^T\pa_wV}{N_0^T(\mathcal{U}\times\mathcal{V})}
\frac{N_0^T(V\times\mathcal{U})}{\mathbf{m}})N_0.
\end{eqnarray}
These can be written for short as
$$\ti\om_u=c_1\ti\om_{uc_1}+c_2\ti\om_{uc_2}+\ti\om_{u1},$$
$$\ti\om_v=-c_1\ti\om_{uc_2}+c_4\ti\om_{uc_1}+\ti\om_{v1},$$
$$\ti\om_w=c_1(\mathcal{\ti V}\ti\om_{uc_1}+\mathcal{\ti U}\ti\om_{uc_2})
+c_2\mathcal{\ti V}\ti\om_{uc_2}-c_4\mathcal{\ti U}\ti\om_{uc_1}+\ti\om_{w1},$$
where $\mathcal{\ti U}:=\frac{N_0^T(\pa_wV\times\mathcal{U})}{N_0^T(\mathcal{U}\times\mathcal{V})},
\mathcal{\ti V}:=\frac{N_0^T(\pa_wV\times\mathcal{V})}{N_0^T(\mathcal{U}\times\mathcal{V})}$.

Note $\ti\om_{w1}=-(\mathcal{M}\mathcal{\ti V}+\frac{\pa_vN_0^T\pa_wV}
{N_0^T(\mathcal{U}\times\mathcal{V})}+\frac{\mathcal{\ti V}}{\mathbf{m}})\ti\om_{uc_1}
+(\mathcal{M}\mathcal{\ti U}+\frac{\pa_uN_0^T\pa_wV}{N_0^T(\mathcal{U}\times\mathcal{V})}
+\frac{\mathcal{\ti U}}{\mathbf{m}})\ti\om_{uc_2}
-\frac{N_0^T(\pa_wV\times V)}{\mathbf{m}^2}N_0$.

Now (\ref{eq:dR1}) becomes

\begin{eqnarray}\label{eq:dtiom}
\pa_u\ti\om_v-\pa_v\ti\om_u+\ti\om_u\times\ti\om_v=0,\nonumber\\
\pa_u\ti\om_w-\pa_w\ti\om_u+\ti\om_u\times\ti\om_w=0,\nonumber\\
\pa_w\ti\om_v-\pa_v\ti\om_w+\ti\om_w\times\ti\om_v=0,
\end{eqnarray}
or
$$(\pa_uc_4-\pa_vc_1)\ti\om_{uc_1}-(\pa_uc_1+\pa_vc_2)\ti\om_{uc_2}
+\pa_u\ti\om_{v1}-\pa_v\ti\om_{u1}+\ti\om_{u1}\times\ti\om_{v1}+
c_1(-\pa_u\ti\om_{uc_2}-\pa_v\ti\om_{uc_1}+
\ti\om_{uc_1}\times\ti\om_{v1}$$
$$-\ti\om_{u1}\times\ti\om_{uc_2})
+c_2(-\pa_v\ti\om_{uc_2}+\ti\om_{uc_2}\times\ti\om_{v1})
+c_4(\pa_u\ti\om_{uc_1}+\ti\om_{u1}\times\ti\om_{uc_1})
-(c_1^2+c_2c_4)\ti\om_{uc_1}\times\ti\om_{uc_2}=0,$$
$$\pa_uc_1(\mathcal{\ti V}\ti\om_{uc_1}+
\mathcal{\ti U}\ti\om_{uc_2})
-\pa_wc_1\ti\om_{uc_1}+(\pa_uc_2\mathcal{\ti V}-\pa_wc_2)\ti\om_{uc_2}
-\pa_uc_4\mathcal{\ti U}\ti\om_{uc_1}+\pa_u\ti\om_{w1}-\pa_w\ti\om_{u1}$$
$$+\ti\om_{u1}\times\ti\om_{w1}+c_1[\pa_u(\mathcal{\ti V}\ti\om_{uc_1}
+\mathcal{\ti U}\ti\om_{uc_2})-\pa_w\ti\om_{uc_1}+\ti\om_{uc_1}\times\ti\om_{w1}
+\ti\om_{u1}\times(\mathcal{\ti V}\ti\om_{uc_1}+\mathcal{\ti U}\ti\om_{uc_2})]$$
$$+c_2[\pa_u(\mathcal{\ti V}\ti\om_{uc_2})-\pa_w\ti\om_{uc_2}
+\ti\om_{uc_2}\times\ti\om_{w1}+\mathcal{\ti V}\ti\om_{u1}\times\ti\om_{uc_2}]$$
$$-c_4[\pa_u(\mathcal{\ti U}\ti\om_{uc_1})+\mathcal{\ti U}\ti\om_{u1}\times\ti\om_{uc_1}]
+(c_1^2+c_2c_4)\mathcal{\ti U}\ti\om_{uc_1}\times\ti\om_{uc_2}=0,$$
$$-\pa_wc_1\ti\om_{uc_2}-\pa_vc_1(\mathcal{\ti V}\ti\om_{uc_1}+
\mathcal{\ti U}\ti\om_{uc_2})-\pa_vc_2\mathcal{\ti V}\ti\om_{uc_2}
+(\pa_vc_4\mathcal{\ti U}+\pa_wc_4)\ti\om_{uc_1}$$
$$+\pa_w\ti\om_{v1}-\pa_v\ti\om_{w1}+\ti\om_{w1}\times\ti\om_{v1}
+c_1[-\pa_w\ti\om_{uc_2}-\pa_v(\mathcal{\ti
V}\ti\om_{uc_1}+\mathcal{\ti U}\ti\om_{uc_2})$$
$$+(\mathcal{\ti V}\ti\om_{uc_1}+
\mathcal{\ti U}\ti\om_{uc_2})\times\ti\om_{v1}-\ti\om_{w1}\times\ti\om_{uc_2}]
+c_2[-\pa_v(\mathcal{\ti V}\ti\om_{uc_2})+\mathcal{\ti V}\ti\om_{uc_2}\times\ti\om_{v1}]$$
$$+c_4[\pa_w\ti\om_{uc_1}+\pa_v(\mathcal{\ti U}\ti\om_{uc_1})
-\mathcal{\ti U}\ti\om_{uc_1}\times\ti\om_{v1}+\ti\om_{w1}\times\ti\om_{uc_1}]
-(c_1^2+c_2c_4)\mathcal{\ti V}\ti\om_{uc_1}\times\ti\om_{uc_2}=0,$$
or for short
$$(\pa_uc_4-\pa_vc_1)\ti\om_{uc_1}-(\pa_uc_1+\pa_vc_2)\ti\om_{uc_2}+C_0^1+c_1C_1^1+c_2C_2^1
+c_4C_4^1+(c_1^2+c_2c_4)C_{124}^1=0,$$
$$\pa_uc_1(\mathcal{\ti V}\ti\om_{uc_1}+\mathcal{\ti U}\ti\om_{uc_2})
-\pa_wc_1\ti\om_{uc_1}+(\pa_uc_2\mathcal{\ti V}-\pa_wc_2)\ti\om_{uc_2}
-\pa_uc_4\mathcal{\ti U}\ti\om_{uc_1}$$
$$+C_0^2+c_1C_1^2+c_2C_2^2+c_4C_4^2+(c_1^2+c_2c_4)C_{124}^2=0,$$
$$-\pa_wc_1\ti\om_{uc_2}-\pa_vc_1(\mathcal{\ti V}\ti\om_{uc_1}+\mathcal{\ti U}\ti\om_{uc_2})
-\pa_vc_2\mathcal{\ti V}\ti\om_{uc_2}
+(\pa_vc_4\mathcal{\ti U}+\pa_wc_4)\ti\om_{uc_1}$$
$$+C_0^3+c_1C_1^3+c_2C_2^3+c_4C_4^3+(c_1^2+c_2c_4)C_{124}^3=0.$$
Multiplying each of these $3$ vector equations on the left respectively with
$(\mathcal{U}\times N_0)^T,\ (\mathcal{V}\times N_0)^T,\ N_0^T$ we get $9$ scalar equations
linear in the $7$ variables $\pa_uc_1,\pa_vc_1,\pa_wc_1,
-\pa_uc_2\mathcal{\ti V}+\pa_wc_2,\pa_vc_2,\ \pa_uc_4,
\pa_vc_4\mathcal{\ti U}+\pa_wc_4$
with the rank of the matrix of the system being
$5$; dividing the rows $1,2,4,5,7,8$ with $N_0^T(\mathcal{U}\times\mathcal{V})$
the linear combinations of rows which are $0$ are as follows:

\begin{eqnarray}\label{eq:L3}
L_3+\frac{N_0^T(V\times\mathcal{V})}{\mathbf{m}}L_1
-\frac{N_0^T(V\times\mathcal{U})}{\mathbf{m}}L_2=0,\
L_6+\frac{N_0^T(V\times\mathcal{V})}{\mathbf{m}}L_4
-\frac{N_0^T(V\times\mathcal{U})}{\mathbf{m}}L_5=0,\nonumber\\
L_9+\frac{N_0^T(V\times\mathcal{V})}{\mathbf{m}}L_7
-\frac{N_0^T(V\times\mathcal{U})}{\mathbf{m}}L_8=0,\
L_7+L_5-\mathcal{\ti V}L_1+\mathcal{\ti U}L_2=0.
\end{eqnarray}
The condition that the augmented matrix of the linear system satisfies the same linear combinations
of rows being $0$ imposes $4$ equations in $c_1,c_2,c_4$; one solves for some of $c_1,c_2,c_4$ from
these equations and then one imposes compatibility conditions on the derivatives of $c_1,c_2,c_4$
and on the commuting of the mixed derivatives of these $3$ functions of the variables $u,v,w$
(using $\pa_uc_1,\pa_vc_1,\pa_wc_1,-\pa_uc_2\mathcal{\ti V}+\pa_wc_2,
\pa_vc_4\mathcal{\ti U}+\pa_wc_4$ of the
$7$ variables above given as solutions of the linear system above); this may impose
further functional relationships between $c_1,c_2,c_4$. At the end of this process we must get
precisely dependence on two functions of a variable, so all necessary conditions on $V$ (and
$\mathbf{m}$, but using (\ref{eq:dm}) this will become dependent on $V$) and on the
seed surface, adjoined to (\ref{eq:intrdst}) become the necessary and sufficient conditions to get
B\"{a}cklund transformation with
isometric correspondence of leaves of a general nature (independent of the shape of the seed).

Since we must get dependence on two functions of a variable, from the $4$ functional relationships
imposed on $c_1,c_2,c_4$ by the consistency requirement of the linear system above only at most two
of them are functionally independent; if we get $2$ functionally independent functional
relationships, then one can solve for two of $c_1,c_2,c_4$ as functions of the remaining third one
and compatibility conditions on derivatives and mixed derivatives do not allow a space of two
functions of a variable dependence of solutions; thus from the $4$ functional relationships
imposed on $c_1,c_2,c_4$ by the consistency requirement of the linear system above only at most one
of them is functionally independent; the conditions that arise from this premise will be conditions
on $V,\ \mathbf{m}$ and the seed surface.

The coefficient of the quadratic term $c_1^2+c_2c_4$ from the first equation of (\ref{eq:L3}) is
$-\frac{|m|^2}{\mathbf{m}^2}N_0^T(\mathcal{U}\times\mathcal{V})\neq 0$ and thus the equation is an
independent (non-vacuous) functional relationship between $c_1,c_2,c_4$.

The coefficient of the quadratic term $c_1^2+c_2c_4$ from the second equation of (\ref{eq:L3}) is
$\frac{|m|^2}{\mathbf{m}^2}N_0^T(\pa_wV\times\mathcal{U})$; thus the second equation is
$-\mathcal{\ti U}$ the first equation.

The coefficient of the quadratic term $c_1^2+c_2c_4$ from the third equation of (\ref{eq:L3}) is
$-\frac{|m|^2}{\mathbf{m}^2}N_0^T(\pa_wV\times\mathcal{V})$; thus the third equation is
$\mathcal{\ti V}$ the first equation.

The fourth equation of (\ref{eq:L3}) does not contain quadratic terms, so it is identically
satisfied.

We thus get the relations:

\begin{eqnarray}\label{eq:mtcj}
m^TC_j^2=-\mathcal{\ti U}m^TC_j^1,\
m^TC_j^3=\mathcal{\ti V}m^TC_j^1,\nonumber\\
-(\pa_wV\times N_0)^TC_j^1
+(\mathcal{U}\times N_0)^TC_j^3+(\mathcal{V}\times N_0)^TC_j^2=0,\ j=0,1,2,4,
\end{eqnarray}
which after removing two obvious relations become:
$$m^T[\mathcal{\ti V}\pa_u\ti\om_{uc_1}-\mathcal{\ti U}\pa_v\ti\om_{uc_1}-\pa_w\ti\om_{uc_1}
+\ti\om_{uc_1}\times\ti\om_{w1}$$
$$+\ti\om_{u1}\times(\mathcal{\ti V}\ti\om_{uc_1}+\mathcal{\ti U}\ti\om_{uc_2})
+\mathcal{\ti U}
(\ti\om_{uc_1}\times\ti\om_{v1}-\ti\om_{u1}\times\ti\om_{uc_2})]=0,$$
$$m^T[\mathcal{\ti V}\pa_u\ti\om_{uc_2}-\mathcal{\ti U}\pa_v\ti\om_{uc_2}
-\pa_w\ti\om_{uc_2}+\ti\om_{uc_2}\times\ti\om_{w1}
+\mathcal{\ti V}\ti\om_{u1}\times\ti\om_{uc_2}
+\mathcal{\ti U}\ti\om_{uc_2}\times\ti\om_{v1}]=0,$$
$$m^T[\pa_u\ti\om_{w1}-\pa_w\ti\om_{u1}+\ti\om_{u1}\times\ti\om_{w1}
+\mathcal{\ti U}(\pa_u\ti\om_{v1}-\pa_v\ti\om_{u1}+\ti\om_{u1}\times\ti\om_{v1})]=0,$$
$$m^T[\mathcal{\ti V}\pa_u\ti\om_{uc_2}-\mathcal{\ti U}\pa_v\ti\om_{uc_2}
-\pa_w\ti\om_{uc_2}+(\mathcal{\ti V}\ti\om_{uc_1}
+\mathcal{\ti U}\ti\om_{uc_2})\times\ti\om_{v1}$$
$$+\ti\om_{w1}\times\ti\om_{vc_1}-\mathcal{\ti V}(
\ti\om_{uc_1}\times\ti\om_{v1}-\ti\om_{u1}\times\ti\om_{uc_2})]=0,$$
$$m^T[\mathcal{\ti V}\pa_u\ti\om_{uc_1}-\mathcal{\ti U}\pa_v\ti\om_{uc_1}-\pa_w\ti\om_{uc_1}
-\mathcal{\ti U}\ti\om_{v1}\times\ti\om_{uc_1}+\ti\om_{uc_1}\times\ti\om_{w1}
+\mathcal{\ti V}\ti\om_{u1}\times\ti\om_{uc_1}]=0,$$
$$m^T[\pa_w\ti\om_{v1}-\pa_v\ti\om_{w1}+\ti\om_{w1}\times\ti\om_{v1}
-\mathcal{\ti V}(\pa_u\ti\om_{v1}-\pa_v\ti\om_{u1}+\ti\om_{u1}\times\ti\om_{v1})]=0,$$
$$(\pa_wV\times N_0)^T[\pa_u\ti\om_{uc_2}+\pa_v\ti\om_{uc_1}
-\ti\om_{uc_1}\times\ti\om_{v1}+\ti\om_{u1}\times\ti\om_{uc_2}]$$
$$+(\mathcal{V}\times N_0)^T[\pa_u(\mathcal{\ti V}\ti\om_{uc_1})
+\mathcal{\ti U}\pa_u\ti\om_{uc_2}-\pa_w\ti\om_{uc_1}+\ti\om_{uc_1}\times\ti\om_{w1}
+\ti\om_{u1}\times(\mathcal{\ti V}\ti\om_{uc_1}+\mathcal{\ti U}\ti\om_{uc_2})]$$
$$+(\mathcal{U}\times N_0)^T[-\pa_w\ti\om_{uc_2}-\mathcal{\ti V}\pa_v\ti\om_{uc_1}
-\pa_v(\mathcal{\ti U}\ti\om_{uc_2})+(\mathcal{\ti V}\ti\om_{uc_1}
+\mathcal{\ti U}\ti\om_{uc_2})
\times\ti\om_{v1}-\ti\om_{w1}\times\ti\om_{uc_2}]=0,$$
$$(\pa_wV\times N_0)^T[\pa_v\ti\om_{uc_2}-\ti\om_{uc_2}\times\ti\om_{v1}]$$
$$+(\mathcal{V}\times N_0)^T[\mathcal{\ti V}\pa_u\ti\om_{uc_2}-\pa_w\ti\om_{uc_2}
+\ti\om_{uc_2}\times\ti\om_{w1}+\mathcal{\ti V}\ti\om_{u1}\times\ti\om_{uc_2}]$$
$$+(\mathcal{U}\times N_0)^T[-\pa_v(\mathcal{\ti V}\ti\om_{uc_2})
+\mathcal{\ti V}\ti\om_{uc_2}\times\ti\om_{v1}]=0,$$
$$-(\pa_wV\times N_0)^T[\pa_u\ti\om_{uc_1}+\ti\om_{u1}\times\ti\om_{uc_1}]
+(\mathcal{V}\times N_0)^T[-\pa_u(\mathcal{\ti U}\ti\om_{uc_1})
-\mathcal{\ti U}\ti\om_{u1}\times\ti\om_{uc_1}]$$
$$+(\mathcal{U}\times N_0)^T[\pa_w\ti\om_{uc_1}+\mathcal{\ti U}\pa_v\ti\om_{uc_1}
-\mathcal{\ti U}\ti\om_{uc_1}\times\ti\om_{v1}+\ti\om_{w1}\times\ti\om_{uc_1}]=0,$$
$$-(\pa_wV\times N_0)^T(\pa_u\ti\om_{v1}-\pa_v\ti\om_{u1}+\ti\om_{u1}\times\ti\om_{v1})
+(\mathcal{V}\times N_0)^T(\pa_u\ti\om_{w1}-\pa_w\ti\om_{u1}+\ti\om_{u1}\times\ti\om_{w1})$$
$$+(\mathcal{U}\times N_0)^T[\pa_w\ti\om_{v1}-\pa_v\ti\om_{w1}+\ti\om_{w1}\times\ti\om_{v1}]=0.$$
The first and fifth relations are equivalent; so are the second and the fourth. We shall call
these relations R1; they depend on $V,\mathbf{m}$ (and their derivatives) and the first and
second fundamental forms of the seed surface $x_0$ (and their derivatives).

From the 9 scalar equations obtained by multiplying on the left respectively with
$(\mathcal{U}\times N_0)^T,\ (\mathcal{V}\times N_0)^T,\ N_0^T$ the 3 vector equations
of (\ref{eq:dtiom}), among which we have the linear
combinations of rows (\ref{eq:L3}), only the next 5 remain linearly independent:

\begin{eqnarray}\label{eq:pauc1}
\pa_uc_1=-\pa_vc_2-[C_0^1+c_1C_1^1+c_2C_2^1+c_4C_4^1+(c_1^2+c_2c_4)C_{124}^1]^T
\frac{(\mathcal{U}\times N_0)}{N_0^T(\mathcal{U}\times\mathcal{V})},\nonumber\\
\pa_vc_1=\pa_uc_4+[C_0^1+c_1C_1^1+c_2C_2^1+c_4C_4^1+(c_1^2+c_2c_4)C_{124}^1]^T
\frac{(\mathcal{V}\times N_0)}{N_0^T(\mathcal{U}\times\mathcal{V})},\nonumber\\
\pa_wc_2=(-\pa_vc_2-[C_0^1+c_1C_1^1+c_2C_2^1+c_4C_4^1+(c_1^2+c_2c_4)C_{124}^1]^T
\frac{(\mathcal{U}\times N_0)}{N_0^T(\mathcal{U}\times\mathcal{V})})\mathcal{\ti U}\nonumber\\
+\pa_uc_2\mathcal{\ti V}-[C_0^2+c_1C_1^2+c_2C_2^2+c_4C_4^2+(c_1^2+c_2c_4)C_{124}^2]^T
\frac{(\mathcal{U}\times N_0)}{N_0^T(\mathcal{U}\times\mathcal{V})},\nonumber\\
\pa_wc_1=-(\pa_vc_2+[C_0^1+c_1C_1^1+c_2C_2^1+c_4C_4^1+(c_1^2+c_2c_4)C_{124}^1]^T
\frac{(\mathcal{U}\times N_0)}{N_0^T(\mathcal{U}\times\mathcal{V})})\mathcal{\ti V}\nonumber\\
-\pa_uc_4\mathcal{\ti U}+[C_0^2+c_1C_1^2+c_2C_2^2+c_4C_4^2+(c_1^2+c_2c_4)C_{124}^2]^T
\frac{(\mathcal{V}\times N_0)}{N_0^T(\mathcal{U}\times\mathcal{V})},\nonumber\\
\pa_wc_4=(\pa_uc_4+[C_0^1+c_1C_1^1+c_2C_2^1+c_4C_4^1+(c_1^2+c_2c_4)C_{124}^1]^T
\frac{(\mathcal{V}\times N_0)}{N_0^T(\mathcal{U}\times\mathcal{V})})\mathcal{\ti V}\nonumber\\
-\pa_vc_4\mathcal{\ti U}-[C_0^3+c_1C_1^3+c_2C_2^3+c_4C_4^3+(c_1^2+c_2c_4)C_{124}^3]^T
\frac{(\mathcal{V}\times N_0)}{N_0^T(\mathcal{U}\times\mathcal{V})}.
\end{eqnarray}

From the first relation of (\ref{eq:L3}) we get:
$$c_4=-\frac{m^T(C_0^1+c_1C_1^1+c_2C_2^1+c_1^2C_{124}^1)}{m^T(C_4^1+c_2C_{124}^1)}.$$
Differentiating this with respect to $u,v$ and $w$ we get the derivatives $\pa_uc_4,\pa_vc_4$
and $\pa_wc_4$ of $c_4$.

Replacing these into (\ref{eq:pauc1}) we get the 5 equations
$$\pa_uc_1=-\pa_vc_2-[\frac{[(C_4^1+c_2C_{124}^1)\times(C_0^1+c_1C_1^1+c_2C_2^1
+c_1^2C_{124}^1)]\times m}{m^T(C_4^1+c_2C_{124}^1)}]^T
\frac{(\mathcal{U}\times N_0)}{N_0^T(\mathcal{U}\times\mathcal{V})},$$
$$\pa_vc_1=-\frac{\pa_u(m^TC_0^1)+c_1\pa_u(m^TC_1^1)+c_2\pa_u(m^TC_2^1)+c_1^2\pa_u(m^TC_{124}^1)}
{m^T(C_4^1+c_2C_{124}^1)}$$
$$+\frac{\pa_u(m^TC_4^1)+c_2\pa_u(m^TC_{124}^1)}{m^T(C_4^1+c_2C_{124}^1)}
\frac{m^T(C_0^1+c_1C_1^1+c_2C_2^1+c_1^2C_{124}^1)}{m^T(C_4^1+c_2C_{124}^1)}+(\pa_vc_2$$
$$+[\frac{[(C_4^1+c_2C_{124}^1)\times(C_0^1+c_1C_1^1+c_2C_2^1
+c_1^2C_{124}^1)]\times m}{m^T(C_4^1+c_2C_{124}^1)}]^T
\frac{(\mathcal{U}\times N_0)}{N_0^T(\mathcal{U}\times\mathcal{V})})
\frac{m^T(C_1^1+2c_1C_{124}^1)}{m^T(C_4^1+c_2C_{124}^1)}$$
$$+\pa_uc_2[-\frac{m^TC_2^1}{m^T(C_4^1+c_2C_{124}^1)}+\frac{m^TC_{124}^1}{m^T(C_4^1+c_2C_{124}^1)}
\frac{m^T(C_0^1+c_1C_1^1+c_2C_2^1+c_1^2C_{124}^1)}{m^T(C_4^1+c_2C_{124}^1)}]$$
$$+[\frac{[(C_4^1+c_2C_{124}^1)\times(C_0^1+c_1C_1^1+c_2C_2^1
+c_1^2C_{124}^1)]\times m}{m^T(C_4^1+c_2C_{124}^1)}]^T
\frac{(\mathcal{V}\times N_0)}{N_0^T(\mathcal{U}\times\mathcal{V})},$$
$$\pa_wc_2=-(\pa_vc_2+[\frac{[(C_4^1+c_2C_{124}^1)\times(C_0^1+c_1C_1^1+c_2C_2^1
+c_1^2C_{124}^1)]\times m}{m^T(C_4^1+c_2C_{124}^1)}]^T
\frac{(\mathcal{U}\times N_0)}{N_0^T(\mathcal{U}\times\mathcal{V})})\mathcal{\ti U}
+\pa_uc_2\mathcal{\ti V}$$
$$-[C_0^2+c_1C_1^2+c_2C_2^2+c_1^2C_{124}^2-(C_4^2+c_2C_{124}^2)
\frac{m^T(C_0^1+c_1C_1^1+c_2C_2^1+c_1^2C_{124}^1)}{m^T(C_4^1+c_2C_{124}^1)}]^T
\frac{(\mathcal{U}\times N_0)}{N_0^T(\mathcal{U}\times\mathcal{V})},$$
$$\pa_wc_1=-(\pa_vc_2+[\frac{[(C_4^1+c_2C_{124}^1)\times(C_0^1+c_1C_1^1+c_2C_2^1
+c_1^2C_{124}^1)]\times m}{m^T(C_4^1+c_2C_{124}^1)}]^T
\frac{(\mathcal{U}\times N_0)}{N_0^T(\mathcal{U}\times\mathcal{V})})\mathcal{\ti V}$$
$$-(-\frac{\pa_u(m^TC_0^1)+c_1\pa_u(m^TC_1^1)+c_2\pa_u(m^TC_2^1)+c_1^2\pa_u(m^TC_{124}^1)}
{m^T(C_4^1+c_2C_{124}^1)}$$
$$+\frac{\pa_u(m^TC_4^1)+c_2\pa_u(m^TC_{124}^1)}{m^T(C_4^1+c_2C_{124}^1)}
\frac{m^T(C_0^1+c_1C_1^1+c_2C_2^1+c_1^2C_{124}^1)}{m^T(C_4^1+c_2C_{124}^1)}+(\pa_vc_2$$
$$+[\frac{[(C_4^1+c_2C_{124}^1)\times(C_0^1+c_1C_1^1+c_2C_2^1
+c_1^2C_{124}^1)]\times m}{m^T(C_4^1+c_2C_{124}^1)}]^T
\frac{(\mathcal{U}\times N_0)}{N_0^T(\mathcal{U}\times\mathcal{V})})
\frac{m^T(C_1^1+2c_1C_{124}^1)}{m^T(C_4^1+c_2C_{124}^1)}$$
$$+\pa_uc_2[-\frac{m^TC_2^1}{m^T(C_4^1+c_2C_{124}^1)}+\frac{m^TC_{124}^1}{m^T(C_4^1+c_2C_{124}^1)}
\frac{m^T(C_0^1+c_1C_1^1+c_2C_2^1+c_1^2C_{124}^1)}{m^T(C_4^1+c_2C_{124}^1)}])\mathcal{\ti U}$$
$$+[C_0^2+c_1C_1^2+c_2C_2^2+c_1^2C_{124}^2-(C_4^2+c_2C_{124}^2)
\frac{m^T(C_0^1+c_1C_1^1+c_2C_2^1+c_1^2C_{124}^1)}{m^T(C_4^1+c_2C_{124}^1)}]^T
\frac{(\mathcal{V}\times N_0)}{N_0^T(\mathcal{U}\times\mathcal{V})},$$
$$-\frac{\pa_w(m^TC_0^1)+c_1\pa_w(m^TC_1^1)+c_2\pa_w(m^TC_2^1)+c_1^2\pa_w(m^TC_{124}^1)}
{m^T(C_4^1+c_2C_{124}^1)}$$
$$+\frac{\pa_w(m^TC_4^1)+c_2\pa_w(m^TC_{124}^1)}{m^T(C_4^1+c_2C_{124}^1)}
\frac{m^T(C_0^1+c_1C_1^1+c_2C_2^1+c_1^2C_{124}^1)}{m^T(C_4^1+c_2C_{124}^1)}$$
$$-(-(\pa_vc_2+[\frac{[(C_4^1+c_2C_{124}^1)\times(C_0^1+c_1C_1^1+c_2C_2^1
+c_1^2C_{124}^1)]\times m}{m^T(C_4^1+c_2C_{124}^1)}]^T
\frac{(\mathcal{U}\times N_0)}{N_0^T(\mathcal{U}\times\mathcal{V})})\mathcal{\ti V}$$
$$-(-\frac{\pa_u(m^TC_0^1)+c_1\pa_u(m^TC_1^1)+c_2\pa_u(m^TC_2^1)+c_1^2\pa_u(m^TC_{124}^1)}
{m^T(C_4^1+c_2C_{124}^1)}$$
$$+\frac{\pa_u(m^TC_4^1)+c_2\pa_u(m^TC_{124}^1)}{m^T(C_4^1+c_2C_{124}^1)}
\frac{m^T(C_0^1+c_1C_1^1+c_2C_2^1+c_1^2C_{124}^1)}{m^T(C_4^1+c_2C_{124}^1)}+(\pa_vc_2$$
$$+[\frac{[(C_4^1+c_2C_{124}^1)\times(C_0^1+c_1C_1^1+c_2C_2^1
+c_1^2C_{124}^1)]\times m}{m^T(C_4^1+c_2C_{124}^1)}]^T
\frac{(\mathcal{U}\times N_0)}{N_0^T(\mathcal{U}\times\mathcal{V})})
\frac{m^T(C_1^1+2c_1C_{124}^1)}{m^T(C_4^1+c_2C_{124}^1)}$$
$$+\pa_uc_2[-\frac{m^TC_2^1}{m^T(C_4^1+c_2C_{124}^1)}+\frac{m^TC_{124}^1}{m^T(C_4^1+c_2C_{124}^1)}
\frac{m^T(C_0^1+c_1C_1^1+c_2C_2^1+c_1^2C_{124}^1)}{m^T(C_4^1+c_2C_{124}^1)}])\mathcal{\ti U}$$
$$+[C_0^2+c_1C_1^2+c_2C_2^2+c_1^2C_{124}^2-(C_4^2+c_2C_{124}^2)
\frac{m^T(C_0^1+c_1C_1^1+c_2C_2^1+c_1^2C_{124}^1)}{m^T(C_4^1+c_2C_{124}^1)}]^T
\frac{(\mathcal{V}\times N_0)}{N_0^T(\mathcal{U}\times\mathcal{V})})$$
$$\frac{m^T(C_1^1+2c_1C_{124}^1)}{m^T(C_4^1+c_2C_{124}^1)}
+(-(\pa_vc_2+[\frac{[(C_4^1+c_2C_{124}^1)\times(C_0^1+c_1C_1^1+c_2C_2^1
+c_1^2C_{124}^1)]\times m}{m^T(C_4^1+c_2C_{124}^1)}]^T
\frac{(\mathcal{U}\times N_0)}{N_0^T(\mathcal{U}\times\mathcal{V})})$$
$$\mathcal{\ti U}+\pa_uc_2\mathcal{\ti V}
-[C_0^2+c_1C_1^2+c_2C_2^2+c_1^2C_{124}^2-(C_4^2$$
$$+c_2C_{124}^2)
\frac{m^T(C_0^1+c_1C_1^1+c_2C_2^1+c_1^2C_{124}^1)}{m^T(C_4^1+c_2C_{124}^1)}]^T
\frac{(\mathcal{U}\times N_0)}{N_0^T(\mathcal{U}\times\mathcal{V})})
[-\frac{m^TC_2^1}{m^T(C_4^1+c_2C_{124}^1)}$$
$$+\frac{m^TC_{124}^1}{m^T(C_4^1+c_2C_{124}^1)}
\frac{m^T(C_0^1+c_1C_1^1+c_2C_2^1+c_1^2C_{124}^1)}{m^T(C_4^1+c_2C_{124}^1)}]$$
$$=(-\frac{\pa_u(m^TC_0^1)+c_1\pa_u(m^TC_1^1)+c_2\pa_u(m^TC_2^1)+c_1^2\pa_u(m^TC_{124}^1)}
{m^T(C_4^1+c_2C_{124}^1)}$$
$$+\frac{\pa_u(m^TC_4^1)+c_2\pa_u(m^TC_{124}^1)}{m^T(C_4^1+c_2C_{124}^1)}
\frac{m^T(C_0^1+c_1C_1^1+c_2C_2^1+c_1^2C_{124}^1)}{m^T(C_4^1+c_2C_{124}^1)}+(\pa_vc_2$$
$$+[\frac{[(C_4^1+c_2C_{124}^1)\times(C_0^1+c_1C_1^1+c_2C_2^1
+c_1^2C_{124}^1)]\times m}{m^T(C_4^1+c_2C_{124}^1)}]^T
\frac{(\mathcal{U}\times N_0)}{N_0^T(\mathcal{U}\times\mathcal{V})})
\frac{m^T(C_1^1+2c_1C_{124}^1)}{m^T(C_4^1+c_2C_{124}^1)}$$
$$+\pa_uc_2[-\frac{m^TC_2^1}{m^T(C_4^1+c_2C_{124}^1)}+\frac{m^TC_{124}^1}{m^T(C_4^1+c_2C_{124}^1)}
\frac{m^T(C_0^1+c_1C_1^1+c_2C_2^1+c_1^2C_{124}^1)}{m^T(C_4^1+c_2C_{124}^1)}]$$
$$+[\frac{[(C_4^1+c_2C_{124}^1)\times(C_0^1+c_1C_1^1+c_2C_2^1
+c_1^2C_{124}^1)]\times m}{m^T(C_4^1+c_2C_{124}^1)}]^T
\frac{(\mathcal{V}\times N_0)}{N_0^T(\mathcal{U}\times\mathcal{V})})\mathcal{\ti V}$$
$$-(-\frac{\pa_v(m^TC_0^1)+c_1\pa_v(m^TC_1^1)+c_2\pa_v(m^TC_2^1)+c_1^2\pa_v(m^TC_{124}^1)}
{m^T(C_4^1+c_2C_{124}^1)}$$
$$+\frac{\pa_v(m^TC_4^1)+c_2\pa_v(m^TC_{124}^1)}{m^T(C_4^1+c_2C_{124}^1)}
\frac{m^T(C_0^1+c_1C_1^1+c_2C_2^1+c_1^2C_{124}^1)}{m^T(C_4^1+c_2C_{124}^1)}$$
$$-(-\frac{\pa_u(m^TC_0^1)+c_1\pa_u(m^TC_1^1)+c_2\pa_u(m^TC_2^1)+c_1^2\pa_u(m^TC_{124}^1)}
{m^T(C_4^1+c_2C_{124}^1)}$$
$$+\frac{\pa_u(m^TC_4^1)+c_2\pa_u(m^TC_{124}^1)}{m^T(C_4^1+c_2C_{124}^1)}
\frac{m^T(C_0^1+c_1C_1^1+c_2C_2^1+c_1^2C_{124}^1)}{m^T(C_4^1+c_2C_{124}^1)}+(\pa_vc_2$$
$$+[\frac{[(C_4^1+c_2C_{124}^1)\times(C_0^1+c_1C_1^1+c_2C_2^1
+c_1^2C_{124}^1)]\times m}{m^T(C_4^1+c_2C_{124}^1)}]^T
\frac{(\mathcal{U}\times N_0)}{N_0^T(\mathcal{U}\times\mathcal{V})})
\frac{m^T(C_1^1+2c_1C_{124}^1)}{m^T(C_4^1+c_2C_{124}^1)}$$
$$+\pa_uc_2[-\frac{m^TC_2^1}{m^T(C_4^1+c_2C_{124}^1)}+\frac{m^TC_{124}^1}{m^T(C_4^1+c_2C_{124}^1)}
\frac{m^T(C_0^1+c_1C_1^1+c_2C_2^1+c_1^2C_{124}^1)}{m^T(C_4^1+c_2C_{124}^1)}]$$
$$+[\frac{[(C_4^1+c_2C_{124}^1)\times(C_0^1+c_1C_1^1+c_2C_2^1
+c_1^2C_{124}^1)]\times m}{m^T(C_4^1+c_2C_{124}^1)}]^T
\frac{(\mathcal{V}\times N_0)}{N_0^T(\mathcal{U}\times\mathcal{V})})
\frac{m^T(C_1^1+2c_1C_{124}^1)}{m^T(C_4^1+c_2C_{124}^1)}$$
$$+\pa_vc_2[-\frac{m^TC_2^1}{m^T(C_4^1+c_2C_{124}^1)}+\frac{m^TC_{124}^1}{m^T(C_4^1+c_2C_{124}^1)}
\frac{m^T(C_0^1+c_1C_1^1+c_2C_2^1+c_1^2C_{124}^1)}{m^T(C_4^1+c_2C_{124}^1)}])\mathcal{\ti U}$$
$$
-[C_0^3+c_1C_1^3+c_2C_2^3+c_1^2C_{124}^3-(C_4^3+c_2C_{124}^3)
\frac{m^T(C_0^1+c_1C_1^1+c_2C_2^1+c_1^2C_{124}^1)}{m^T(C_4^1+c_2C_{124}^1)}]^T
\frac{(\mathcal{V}\times N_0)}
{N_0^T(\mathcal{U}\times\mathcal{V})}.$$
In the last relation the coefficients of $\pa_uc_2,\pa_vc_2$ are $0$; thus this relation is
a functional relationship (which is actually rational) between $c_1$ and $c_2$. Since $c_1$ and
$c_2$ are functionally independent, this last relation being identically $0$ gives relations
depending on $V,\mathbf{m}$ (and their derivatives) and the first and second fundamental form
(and their derivatives) of the seed surface $x_0$. We shall call these relations (and those
obtained by exchanging the r\^{o}le played by $c_2$ and $c_4$) R2.

The remaining 4 equations can be written as:
\begin{eqnarray}\label{eq:pauc1D}
\pa_uc_1=-\pa_vc_2+D_1,\nonumber\\
\pa_vc_1=D_2+\pa_uc_2D_3+\pa_vc_2D_4,\nonumber\\
\pa_wc_1=D_5+\pa_uc_2D_6+\pa_vc_2D_7,\nonumber\\
\pa_wc_2=D_8+\pa_uc_2\mathcal{\ti V}-\pa_vc_2\mathcal{\ti U},
\end{eqnarray}
where $D_j=D_j(c_1,c_2,u,v,w),\ j=1,...,8$ are rational functions of $c_1$ and $c_2$ and depend
otherwise on $V,\mathbf{m}$ (and their derivatives) and on the first and second fundamental form of the
seed surface $x_0$ (and their derivatives).

Imposing the three compatibility conditions $\pa_v(\pa_uc_1)=\pa_u(\pa_vc_1),\
\pa_w(\pa_uc_1)=\pa_u(\pa_wc_1),\ \pa_w(\pa_vc_1)\\=\pa_v(\pa_wc_1)$ and using the equations
(\ref{eq:pauc1D}) themselves we get the relations
$$-\pa_v^2c_2+\pa_{c_1}D_1(D_2+\pa_uc_2D_3+\pa_vc_2D_4)+\pa_{c_2}D_1\pa_vc_2+\pa_vD_1$$
$$-(\pa_{c_1}D_2+\pa_uc_2\pa_{c_1}D_3+\pa_vc_2\pa_{c_1}D_4)(-\pa_vc_2
+D_1)-(\pa_{c_2}D_2+\pa_uc_2\pa_{c_2}D_3+\pa_vc_2\pa_{c_2}D_4)\pa_uc_2$$
$$-(\pa_uD_2+\pa_uc_2\pa_uD_3+\pa_vc_2\pa_uD_4)-(\pa_u^2c_2D_3+\pa_{uv}^2c_2D_4)=0,$$
$$-[\pa_{c_1}D_8(D_2+\pa_uc_2D_3+\pa_vc_2D_4)+\pa_{c_2}D_8\pa_vc_2+\pa_vD_8+
\pa_{uv}^2c_2\mathcal{\ti V}+\pa_uc_2\pa_v\mathcal{\ti V}$$
$$-\pa_v^2c_2\mathcal{\ti U}-\pa_vc_2\pa_v\mathcal{\ti U}]
+\pa_{c_1}D_1(D_5+\pa_uc_2D_6+\pa_vc_2D_7)$$
$$+\pa_{c_2}D_1(D_8+\pa_uc_2\mathcal{\ti V}-\pa_vc_2\mathcal{\ti U})+\pa_wD_1$$
$$-(\pa_{c_1}D_5+\pa_uc_2\pa_{c_1}D_6+\pa_vc_2\pa_{c_1}D_7)(-\pa_vc_2+D_1)
-(\pa_{c_2}D_5+\pa_uc_2\pa_{c_2}D_6+\pa_vc_2\pa_{c_2}D_7)\pa_uc_2$$
$$-(\pa_uD_5+\pa_uc_2\pa_uD_6+\pa_vc_2\pa_uD_7)-(\pa_u^2c_2D_6+\pa_{uv}^2c_2D_7)=0,$$
$$(\pa_{c_1}D_2+\pa_uc_2\pa_{c_1}D_3+\pa_vc_2\pa_{c_1}D_4)(D_5+\pa_uc_2D_6+\pa_vc_2D_7)
+(\pa_{c_2}D_2+\pa_uc_2\pa_{c_2}D_3+\pa_vc_2\pa_{c_2}D_4)(D_8$$
$$+\pa_uc_2\mathcal{\ti V}-\pa_vc_2\mathcal{\ti U})
+\pa_wD_2+\pa_uc_2\pa_wD_3+\pa_vc_2\pa_wD_4+D_3[\pa_{c_1}D_8(-\pa_vc_2+D_1)$$
$$+\pa_{c_2}D_8\pa_uc_2+\pa_uD_8+\pa_u^2c_2\mathcal{\ti V}
+\pa_uc_2\pa_u\mathcal{\ti V}-\pa_{uv}^2c_2\mathcal{\ti U}$$
$$-\pa_vc_2\pa_u\mathcal{\ti U}]
+D_4[\pa_{c_1}D_8(D_2+\pa_uc_2D_3+\pa_vc_2D_4)+\pa_{c_2}D_8\pa_vc_2+\pa_vD_8$$
$$+\pa_{uv}^2c_2\mathcal{\ti V}+\pa_uc_2\pa_v\mathcal{\ti V}
-\pa_v^2c_2\mathcal{\ti U}-\pa_vc_2\pa_v\mathcal{\ti U}]$$
$$-(\pa_{c_1}D_5+\pa_uc_2\pa_{c_1}D_6+\pa_vc_2\pa_{c_1}D_7)(D_2+\pa_uc_2D_3+\pa_vc_2D_4)$$
$$-(\pa_{c_2}D_5+\pa_uc_2\pa_{c_2}D_6+\pa_vc_2\pa_{c_2}D_7)\pa_vc_2
-(\pa_vD_5+\pa_uc_2\pa_vD_6+\pa_vc_2\pa_vD_7)-(\pa_{uv}^2c_2D_6+\pa_v^2c_2D_7)=0.$$
These constitute a linear system in the variables $\pa_u^2c_2,\pa_{uv}^2c_2,\pa_v^2c_2$ with
the determinant of the matrix of the system being
$$\begin{vmatrix}-D_3& -D_4& -1\\ -D_6&-D_7+\mathcal{\ti V}&-\mathcal{\ti U}\\
\mathcal{\ti V}&-D_6+\mathcal{\ti V}-\mathcal{\ti U}&-D_7-\mathcal{\ti U}\end{vmatrix}=
-\begin{vmatrix}-D_6+D_3\mathcal{\ti U}&-D_7+\mathcal{\ti V}+D_4\mathcal{\ti U}\\
\mathcal{\ti V}+D_3(D_7+\mathcal{\ti U})&-D_6+\mathcal{\ti V}-\mathcal{\ti U}
+D_4(D_7+\mathcal{\ti U})\end{vmatrix}.$$
Using $D_3=-\frac{m^TC_2^1}{m^T(C_4^1+c_2C_{124}^1)}+\frac{m^TC_{124}^1}{m^T(C_4^1+c_2C_{124}^1)}
\frac{m^T(C_0^1+c_1C_1^1+c_2C_2^1+c_1^2C_{124}^1)}{m^T(C_4^1+c_2C_{124}^1)},\
D_4=\frac{m^T(C_1^1+2c_1C_{124}^1)}{m^T(C_4^1+c_2C_{124}^1)},\ D_6=-\mathcal{\ti U}D_3,\
D_7=-\mathcal{\ti V}-\mathcal{\ti U}D_4$ the coefficient of the highest order term $c_1^4$ in the
numerator (after bringing to a common denominator $[m^T(C_4^1+c_2C_{124}^1)]^4$) of the above
determinant is $2(m^TC_{124}^1)^4\mathcal{\ti U}^2\neq 0$, so this determinant is not $0$.

This allows us to solve for
\begin{eqnarray}\label{eq:pauuc2}
\pa_u^2c_2=F_1(\pa_uc_2,\pa_vc_2,c_1,c_2,u,v,w),\nonumber\\
\pa_{uv}^2c_2=F_2(\pa_uc_2,\pa_vc_2,c_1,c_2,u,v,w),\nonumber\\
\pa_v^2c_2=F_3(\pa_uc_2,\pa_vc_2,c_1,c_2,u,v,w),
\end{eqnarray}
where $F_j,\ j=1,2,3$ depend quadratically on $\pa_uc_2,\pa_vc_2$, rationally on $c_1,c_2$
and depend otherwise on $V,\mathbf{m}$ (and their derivatives) and on the first and second fundamental
form of the seed surface $x_0$ (and their derivatives).

Imposing the two compatibility conditions $\pa_v(\pa_u^2c_2)=\pa_u(\pa_{uv}^2c_2),\
\pa_u(\pa_v^2c_2)=\pa_v(\pa_{uv}^2c_2)$ and using the equations (\ref{eq:pauc1D}),
(\ref{eq:pauuc2}) themselves we get the relations
$$\pa_{\pa_uc_2}F_1F_2+\pa_{\pa_vc_2}F_1F_3+\pa_{c_1}F_1(D_2+\pa_uc_2D_3+\pa_vc_2D_4)
+\pa_{c_2}F_1\pa_vc_2+\pa_vF_1$$
$$-(\pa_{\pa_uc_2}F_2F_1+\pa_{\pa_vc_2}F_2F_2+\pa_{c_1}F_2(-\pa_vc_2+D_1)
+\pa_{c_2}F_2\pa_uc_2+\pa_uF_2)=0,$$
$$\pa_{\pa_uc_2}F_3F_1+\pa_{\pa_vc_2}F_3F_2+\pa_{c_1}F_3(-\pa_vc_2+D_1)
+\pa_{c_2}F_3\pa_uc_2+\pa_uF_3$$
$$-(\pa_{\pa_uc_2}F_2F_2+\pa_{\pa_vc_2}F_2F_3+\pa_{c_1}F_2(D_2+\pa_uc_2D_3+\pa_vc_2D_4)
+\pa_{c_2}F_2\pa_vc_2+\pa_vF_2)=0.$$
These can be written as
\begin{eqnarray}\label{eq:G}
G_1(\pa_uc_2,\pa_vc_2,c_1,c_2,u,v,w)=0,\ G_2(\pa_uc_2,\pa_vc_2,c_1,c_2,u,v,w)=0,
\end{eqnarray}
where $G_1,G_2$ depend cubically on $\pa_uc_2,\pa_vc_2$, rationally on $c_1,c_2$ and depend
otherwise on $V,\mathbf{m}$ (and their derivatives) and on the first and second fundamental
form of the seed surface $x_0$ (and their derivatives).

If $G_1,G_2$ are functionally independent in $\pa_uc_2,\pa_vc_2$, then one can solve them for
$\pa_uc_2=H_1(c_1,c_2,u,v,w),\ \pa_vc_2=H_2(c_1,c_2,u,v,w)$; replacing this into (\ref{eq:pauc1D})
we get all derivatives of $c_1$ and $c_2$ as functions of $c_1$ and $c_2$ and depend otherwise
on $V,\mathbf{m}$ (and their derivatives) and on the first and second fundamental form of the seed
surface $x_0$ (and their derivatives). Imposing the compatibility conditions of commuting of
mixed derivatives of $c_1$ and $c_2$ and using the derivatives of $c_1$ and $c_2$ themselves as
functions of $c_1$ and $c_2$ we get 6 compatibility conditions which must be identically satisfied
in $c_1$ and $c_2$, but in this case the space of solutions depends on constants and not on two
functions of a variable.

Since the coefficients of $\pa_uc_2^3,\ \pa_vc_2^3$ in $G_1$ are not $0$, $G_1$ is a
non-vacuous functional relationship between $\pa_uc_2,\ \pa_vc_2$, so $G_2$ must be a multiple
(namely $-\frac{\mathcal{\ti V}}{\mathcal{\ti U}}$) of it.

We thus have $G_2+\frac{\mathcal{\ti V}}{\mathcal{\ti U}}G_1=0$; after bringing to a common
denominator the numerator will be quadratic in $\pa_uc_2,\ \pa_vc_2$ and polynomial in
$c_1,c_2$ which must be identically $0$. We shall call these relations (and those obtained
by exchanging the r\^{o}le played by $c_2$ and $c_4$) R3.

We can solve $G_1=0$ for $\pa_uc_2=H(\pa_vc_2,c_1,c_2,u,v,w)$ (there are $3$ choices of $H$);
replacing this into (\ref{eq:pauc1D}) we get
\begin{eqnarray}\label{eq:pauc1DH}
\pa_uc_1=-\pa_vc_2+D_1,\nonumber\\
\pa_vc_1=D_2+H(\pa_vc_2,c_1,c_2,u,v,w)D_3+\pa_vc_2D_4,\nonumber\\
\pa_wc_1=D_5+H(\pa_vc_2,c_1,c_2,u,v,w)D_6+\pa_vc_2D_7,\nonumber\\
\pa_uc_2=H(\pa_vc_2,c_1,c_2,u,v,w),\nonumber\\
\pa_wc_2=D_8+H(\pa_vc_2,c_1,c_2,u,v,w)\mathcal{\ti V}-\pa_vc_2\mathcal{\ti U},
\end{eqnarray}
 and (\ref{eq:pauuc2}) becomes
 $$\pa_{\pa_vc_2}HF_2(H,\pa_vc_2,c_1,c_2,u,v,w)+\pa_{c_1}H(-\pa_vc_2+D_1)+\pa_{c_2}HH+\pa_uH
 =F_1(H,\pa_vc_2,c_1,c_2,u,v,w),$$
$$\pa_{\pa_vc_2}HF_3(H,\pa_vc_2,c_1,c_2,u,v,w)+\pa_{c_1}H(D_2+HD_3+\pa_vc_2D_4)
+\pa_{c_2}H\pa_vc_2+\pa_wH$$
$$=F_2(H,\pa_vc_2,c_1,c_2,u,v,w),$$
$$\pa_v^2c_2=F_3(H,\pa_vc_2,c_1,c_2,u,v,w).$$
The first two relations are valid for functionally independent $c_1,c_2,\pa_vc_2$, so impose
conditions on $V,\mathbf{m}$ (and their derivatives) and on the first and second fundamental
form of the seed surface $x_0$ (and their derivatives).

Imposing the compatibility condition $\pa_w(\pa_uc_2)=\pa_u(\pa_wc_2)$ on the last two equations
of (\ref{eq:pauc1DH}) and using the equations of (\ref{eq:pauc1DH}) themselves and
$\pa_v^2c_2=F_3(H,\pa_vc_2,c_1,c_2,u,v,w)$ we get
$$\pa_{\pa_vc_2}H[\pa_{c_1}D_8(D_2+HD_3+\pa_vc_2D_4)
+\pa_{c_2}D_8\pa_vc_2+\pa_vD_8+(\pa_{\pa_vc_2}HF_3(H,\pa_vc_2,c_1,c_2,u,v,w)$$
$$+\pa_{c_1}H(D_2+HD_3+\pa_vc_2D_4)+\pa_{c_2}H\pa_vc_2
+\pa_vH)\mathcal{\ti V}+H\pa_v\mathcal{\ti V}-F_3(H,\pa_vc_2,c_1,c_2,u,v,w)\mathcal{\ti U}
-\pa_vc_2\pa_v\mathcal{\ti U}]$$
$$+\pa_{c_1}H(D_5+HD_6+\pa_vc_2D_7)
+\pa_{c_2}H(D_8+H\mathcal{\ti V}-\pa_vc_2\mathcal{\ti U})+\pa_wH
-[\pa_{c_1}D_8(-\pa_vc_2+D_1)+\pa_{c_2}D_8H+\pa_uD_8$$
$$+(\pa_{\pa_vc_2}HF_2(H,\pa_vc_2,c_1,c_2,u,v,w)
+\pa_{c_1}H(-\pa_vc_2+D_1)+\pa_{c_2}HH+\pa_uH)\mathcal{\ti V}+H\pa_u\mathcal{\ti V}$$
$$-F_2(H,\pa_vc_2,c_1,c_2,u,v,w)\mathcal{\ti U}-\pa_vc_2\pa_u\mathcal{\ti U}]=0.$$
This is valid for functionally independent $c_1,c_2,\pa_vc_2$, so impose
conditions on $V,\mathbf{m}$ (and their derivatives) and on the first and second fundamental
form of the seed surface $x_0$ (and their derivatives).

Thus (\ref{eq:pauc1DH}) together with $\pa_v^2c_2=F_3(H,\pa_vc_2,c_1,c_2,u,v,w)$ is in
involution and we know the solution must depend on two functions of a variable.

The conditions on $V,
\mathbf{m}$ (and their derivatives) and on the first and second fundamental form of the seed
surface $x_0$ (and their derivatives) obtained at this last step will be called the relations R4.
Note that here we don't need to exchange the r\^{o}le played by $c_2$ and $c_4$, since there will
be $9$ cases to discuss instead of $3$.

For the relations R1,R2,R3,R4 we can replace the seed $x_0$ with any other isometric surface;
this will remove the dependence on the second fundamental form of the seed surface $x_0$
(and its derivatives). Then we can replace $\mathbf{m}$ and its derivatives from (\ref{eq:dm}),
thus giving a differential system on $V$ depending only on the first
fundamental form of the seed surface $x_0$, which adjoined to (\ref{eq:intrdst}) constitutes a
differential system on $V$ (here compatibility conditions may impose conditions on the metric
of the seed surface) which completely describes the B\"{a}cklund transformation with
isometric correspondence of leaves of a general nature (independent of the shape of the seed
surface).

\subsection{The relations R2}

After bringing to a common denominator $[m^T(C_4^1+c_2C_{124}^1)]^3
N_0^T(\mathcal{U}\times\mathcal{V})$ the last relation of (\ref{eq:pauc1}) becomes a
quartic polynomial $P_1$ in $c_1,c_2$ being identically $0$. We shall consider the coefficients of
decreasing powers of $c_1,c_2$, establishing recurrence relations.

From the coefficient of $c_1^4$ we get $(\mathcal{U}\times N_0)^TC_4^2=
-\mathcal{\ti U}(\mathcal{U}\times N_0)^TC_4^1$, which is obvious and thus removed;
from the coefficient of $c_2^4$ we get $(\mathcal{V}\times N_0)^TC_2^3=
\mathcal{\ti V}(\mathcal{V}\times N_0)^TC_2^1$, which is again obvious and thus removed.

From the coefficient of $c_1^3c_2$ we get $(\mathcal{V}\times N_0)^TC_4^2=
-\mathcal{\ti U}(\mathcal{V}\times N_0)^TC_4^1
+\frac{1}{2}(\mathcal{U}\times N_0)^T(C_1^2+\mathcal{\ti U}C_1^1)$;
from the coefficient of $c_1c_2^3$ we get $(\mathcal{V}\times N_0)^TC_1^3=
(\mathcal{V}\times N_0)^T(\mathcal{\ti V}C_1^1+2\mathcal{\ti U}C_2^1+2C_2^2)$.

From the coefficient of $c_1^2c_2^2$ we get $(\mathcal{V}\times N_0)^TC_4^3=
(\mathcal{V}\times N_0)^T(\mathcal{\ti V}C_4^1-2\mathcal{\ti U}C_1^1-2C_1^2)
-(\mathcal{U}\times N_0)^T(\mathcal{\ti U}C_2^1+C_2^2)$.

Using the relations thus far obtained, the coefficient of $c_1^3$ is $0$.

From the coefficient of $c_1^2c_2$ we get $(\mathcal{U}\times N_0)^TC_2^2=
\frac{1}{2m^TC_{124}^1m^TC_4^1}[2(m^TC_{124}^1)^2
(\mathcal{U}\times N_0)^T(\mathcal{\ti U}C_0^1+C_0^2)-2\mathcal{\ti U}m^TC_{124}^1m^TC_4^1
(\mathcal{U}\times N_0)^TC_2^1-m^TC_1^1m^TC_{124}^1
(\mathcal{U}\times N_0)^T(\mathcal{\ti U}C_1^1+C_1^2)\\
+2\mathcal{\ti V}N_0^T(\mathcal{U}\times\mathcal{V})[m^TC_{124}^1\pa_u(m^TC_4^1)
-\pa_u(m^TC_{124}^1)m^TC_4^1]
-2\mathcal{\ti U}N_0^T(\mathcal{U}\times\mathcal{V})[m^TC_{124}^1\pa_v(m^TC_4^1)
-\pa_v(m^TC_{124}^1)m^TC_4^1]
-2N_0^T(\mathcal{U}\times\mathcal{V})[m^TC_{124}^1\pa_w(m^TC_4^1)
-\pa_w(m^TC_{124}^1)m^TC_4^1]]$.

Using the relations thus far obtained, the coefficient of $c_1^2$ is $0$.

From the coefficient of $c_1c_2^2$ we get $(\mathcal{V}\times N_0)^TC_2^2=
\frac{-1}{2m^TC_{124}^1m^TC_4^1}[-2(m^TC_{124}^1)^2
(\mathcal{V}\times N_0)^T(\mathcal{\ti U}C_0^1+C_0^2)+2\mathcal{\ti U}m^TC_{124}^1m^TC_4^1
(\mathcal{V}\times N_0)^TC_2^1+m^TC_1^1m^TC_{124}^1
(\mathcal{V}\times N_0)^T(\mathcal{\ti U}C_1^1+C_1^2)\\
+\mathcal{\ti V}N_0^T(\mathcal{U}\times\mathcal{V})[m^TC_{124}^1\pa_u(m^TC_1^1)
-\pa_u(m^TC_{124}^1)m^TC_1^1]
-\mathcal{\ti U}N_0^T(\mathcal{U}\times\mathcal{V})[m^TC_{124}^1\pa_v(m^TC_1^1)\\
-\pa_v(m^TC_{124}^1)m^TC_1^1]
-N_0^T(\mathcal{U}\times\mathcal{V})[m^TC_{124}^1\pa_w(m^TC_1^1)
-\pa_w(m^TC_{124}^1)m^TC_1^1]]+m^TC_{124}^1m^TC_2^1(\mathcal{U}\times N_0)^TC_1^2$.

Using the relations thus far obtained, the coefficients of $c_1c_2$ and $c_1$ are $0$.

From the coefficient of $1$ we get $(\mathcal{V}\times N_0)^TC_0^3=
\frac{-1}{m^TC_{124}^1m^TC_4^1}[2m^TC_{124}^1m^TC_0^1
(\mathcal{V}\times N_0)^T(\mathcal{\ti U}C_1^1+C_1^2)
+m^TC_{124}^1[m^TC_2^1(\mathcal{U}\times N_0)^TC_0^2-m^TC_1^1(\mathcal{V}\times N_0)^TC_0^2]
+\mathcal{\ti U}m^TC_{124}^1[m^TC_2^1(\mathcal{U}\times N_0)^TC_0^1
-m^TC_1^1(\mathcal{V}\times N_0)^TC_0^1]
-\mathcal{\ti V}m^TC_{124}^1m^TC_4^1(\mathcal{V}\times N_0)^TC_0^1
+\mathcal{\ti V}N_0^T(\mathcal{U}\times\mathcal{V})[m^TC_{124}^1\pa_u(m^TC_0^1)
-\pa_u(m^TC_{124}^1)m^TC_0^1]
-\mathcal{\ti U}N_0^T(\mathcal{U}\times\mathcal{V})[m^TC_{124}^1\pa_v(m^TC_0^1)
-\pa_v(m^TC_{124}^1)m^TC_0^1]
-N_0^T(\mathcal{U}\times\mathcal{V})[m^TC_{124}^1\pa_w(m^TC_0^1)\\
-\pa_w(m^TC_{124}^1)m^TC_0^1]]$ and we only need to further consider the coefficients of $c_2,c_2^2,c_2^3$.

Since these satisfy $2m^TC_4^1(\mathrm{coeff\ of}\ c_2^3)=
m^TC_{124}^1(\mathrm{coeff\ of}\ c_2^2),\ 2m^TC_{124}^1(\mathrm{coeff\ of}\ c_2)\\=
m^TC_4^1(\mathrm{coeff\ of}\ c_2^2)$, it is enough to consider only the coefficient of $c_2^2$
and we get

$m^TC_{124}^1[(\mathcal{V}\times N_0)^TC_1^2+\mathcal{\ti U}(\mathcal{V}\times N_0)^TC_1^1]
[(m^TC_1^1)^2-4m^TC_0^1m^TC_{124}^1+4m^TC_2^1m^TC_4^1]\\
+2\mathcal{\ti U}N_0^T(\mathcal{U}\times\mathcal{V})m^TC_{124}^1
[m^TC_{124}^1\pa_v(m^TC_0^1)-\pa_v(m^TC_{124}^1)m^TC_0^1]\\
+2N_0^T(\mathcal{U}\times\mathcal{V})m^TC_4^1[\pa_w(m^TC_{124}^1)m^TC_2^1
-m^TC_{124}^1\pa_w(m^TC_2^1)]\\
+2N_0^T(\mathcal{U}\times\mathcal{V})m^TC_2^1[\pa_w(m^TC_{124}^1)m^TC_4^1
-m^TC_{124}^1\pa_w(m^TC_4^1)]\\
+2N_0^T(\mathcal{U}\times\mathcal{V})m^TC_{124}^1[m^TC_{124}^1\pa_w(m^TC_0^1)
-\pa_w(m^TC_{124}^1)m^TC_0^1]\\
+2\mathcal{\ti V}N_0^T(\mathcal{U}\times\mathcal{V})m^TC_4^1[m^TC_{124}^1\pa_u(m^TC_2^1)
-\pa_u(m^TC_{124}^1)m^TC_2^1]\\
+2\mathcal{\ti V}N_0^T(\mathcal{U}\times\mathcal{V})m^TC_2^1[m^TC_{124}^1\pa_u(m^TC_4^1)
-\pa_u(m^TC_{124}^1)m^TC_4^1]\\
-2\mathcal{\ti V}N_0^T(\mathcal{U}\times\mathcal{V})m^TC_{124}^1[m^TC_{124}^1\pa_u(m^TC_0^1)
-\pa_u(m^TC_{124}^1)m^TC_0^1]\\
-2\mathcal{\ti U}N_0^T(\mathcal{U}\times\mathcal{V})m^TC_4^1[m^TC_{124}^1\pa_v(m^TC_2^1)
-\pa_v(m^TC_{124}^1)m^TC_2^1]\\
-2\mathcal{\ti U}N_0^T(\mathcal{U}\times\mathcal{V})m^TC_2^1[m^TC_{124}^1\pa_v(m^TC_4^1)
-\pa_v(m^TC_{124}^1)m^TC_4^1]\\
-\mathcal{\ti U}N_0^T(\mathcal{U}\times\mathcal{V})m^TC_1^1[m^TC_{124}^1\pa_v(m^TC_1^1)
-\pa_v(m^TC_{124}^1)m^TC_1^1]\\
-N_0^T(\mathcal{U}\times\mathcal{V})m^TC_1^1[m^TC_{124}^1\pa_w(m^TC_1^1)
-\pa_w(m^TC_{124}^1)m^TC_1^1]\\
+\mathcal{\ti V}N_0^T(\mathcal{U}\times\mathcal{V})m^TC_1^1[m^TC_{124}^1\pa_u(m^TC_1^1)
-\pa_u(m^TC_{124}^1)m^TC_1^1]=0$.

By exchanging the r\^{o}le played by $c_2$ and $c_4$ and after bringing to a common denominator
$[m^T(C_2^1+c_4C_{124}^1)]^3N_0^T(\mathcal{U}\times\mathcal{V})$ the relation corresponding to
the last relation of (\ref{eq:pauc1}) by exchanging the r\^{o}le played by $c_2$ and $c_4$ becomes
a quartic polynomial $P_2$ in $c_1,c_4$ being identically $0$. Again we shall consider the
coefficients of decreasing powers of $c_1,c_4$ establishing recurrence relations.

It turns out that we don't get more relations than those obtained by the polynomial $P_1$ being
identically $0$: the coefficient of $c_1^4$ in $P_2$ is the coefficient
of $c_2^4$ in $P_1$; the coefficient of $c_4^4$ in $P_2$ is the coefficient
of $c_1^4$ in $P_1$; the coefficient of $c_1^3c_4$ in $P_2$ is $-$ the coefficient of $c_1c_2^3$ in
$P_1$; the coefficient of $c_1c_4^3$ in $P_2$ is $-$ the coefficient of $c_1^3c_2$ in $P_1$;
the coefficient of $c_1^2c_4^2$ in $P_2$ is the coefficient of $c_1^2c_2^2$ in $P_1$; the
coefficient of $c_1^3$ in $P_2$ is $0$; the coefficient of $c_1c_4^2$ in $P_2$ is the coefficient
of $c_1c_2^2$ in $P_1$. For the coefficient of $c_1^2c_4$ in $P_2$, by replacing
$(\mathcal{U}\times N_0)^TC_2^2$ and $(\mathcal{V}\times N_0)^TC_0^3$ with their values it becomes
$-$ the coefficient of $c_2^3$ in $P_1$. The coefficient of $c_4^3$ in $P_2$ is $0$; the
coefficient of $c_1^2$ in $P_2$ is $-\frac{m^TC_2^1}{m^TC_{124}^1}$ times the coefficient of
$c_2^3$ in $P_1$; the coefficient of $c_1c_4$ in $P_2$ is $-\frac{m^TC_1^1}{m^TC_{124}^1}$
times the coefficient of $c_2^3$ in $P_1$; the coefficient of $c_4^2$ in $P_2$ is
$-\frac{m^TC_4^1}{m^TC_{124}^1}$ times the coefficient of $c_2^3$ in $P_1$; the coefficient of
$c_1$ in $P_2$ is $-\frac{m^TC_2^1m^TC_1^1}{(m^TC_{124}^1)^2}$ times the coefficient of $c_2^3$
in $P_1$; the coefficient of $c_4$ in $P_2$ is
$-\frac{m^TC_{124}^1m^TC_0^1+m^TC_2^1m^TC_4^1}{(m^TC_{124}^1)^2}$ times  the coefficient of
$c_2^3$ in $P_1$; the coefficient of $1$ in $P_2$ is $-\frac{m^TC_2^1m^TC_0^1}{(m^TC_{124}^1)^2}$
times the coefficient of $c_2^3$ in $P_1$.

\subsection{The relations (\ref{eq:intrdst}), R1 and R2 revisited}

The relations that depend only on the second fundamental form of
$x_0$ (and not its derivatives) are of the form
$\pa_uN_0^TA+\pa_vN_0^TB+C=0$, where $A,B,C$ do not depend on the
second fundamental form of $x_0$. From each such relation we get
four ones: with $u^1:=u,\ u^2:=v,
g_{jk}=\pa_{u^j}x_0^T\pa_{u^k}x_0$ we have
$g^{1j}\pa_{u^j}x_0^TA=0,\ g^{2j}\pa_{u^j}x_0^TB=0,\
g^{2j}\pa_{u^j}x_0^TA+g^{1j}\pa_{u^j}x_0^TB=0,\ C=0$; if $A=0$ (or
$B=0$), then $B=0$ (or $A=0$) and we obtain only three relations.

For the fifth equation of R1 we have
$A=-N_0^T(\pa_wV\times\mathcal{V}){N_0^T(\mathcal{U}\times\mathcal{V})}N_0^T(V\times
\mathcal{U})N_0\times
V-N_0^T(\pa_wV\times\mathcal{V})[N_0^T(\mathcal{U}\times\mathcal{V})
\mathbf{m}^2+V^T\mathcal{U}N_0^T(V\times\mathcal{V})](N_0\times\mathcal{U})\times
N_0
+[(|V|^2+\mathbf{m}^2)N_0^T(\pa_wV\times\mathcal{U})N_0^T(\mathcal{U}\times\mathcal{V}))
+V^T\mathcal{U}N_0^T(\pa_wV\times\mathcal{V})N_0^T(V\times\mathcal{U})]
(N_0\times\mathcal{V})\times
N_0+(\mathbf{m}^2+|V|^2)(N_0^T(\mathcal{U}\times\mathcal{V}))^2
\pa_wV$, $B=N_0^T(\pa_wV\times\mathcal{U})N_0^T(V\times
\mathcal{U})[N_0^T(\mathcal{U}\times\mathcal{V})N_0\times V
+\mathcal{V}^TV(N_0\times\mathcal{U})\times N_0-V^T\mathcal{U}
(N_0\times\mathcal{V})\times N_0]=0,\\
C=[-\frac{N_0^T(\pa_wV\times\pa_vx_0)N_0^T(\pa_wV\times \pa_uV)}
{N_0^T(\mathcal{U}\times\mathcal{V})N_0^T(\pa_wV\times V)}
+\frac{N_0^T(\pa_wV\times\pa_ux_0)N_0^T(\pa_wV\times \pa_vV)}
{N_0^T(\mathcal{U}\times\mathcal{V})N_0^T(\pa_wV\times V)}]N_0^T(V\times\mathcal{U})\\
-\frac{N_0^T(\pa_wV\times\mathcal{V})}{N_0^T(\mathcal{U}\times\mathcal{V})}N_0^T(\pa_ux_0\times\pa_uV)
+\frac{N_0^T(\pa_wV\times\mathcal{U})}{N_0^T(\mathcal{U}\times\mathcal{V})}N_0^T(\mathcal{U}\times\pa_vV)
-(1+\mathbf{m}^2
\frac{KN_0^T(\pa_ux_0\times\pa_vx_0)}{N_0^T(\mathcal{U}\times\mathcal{V})})
N_0^T(\pa_wV\times\mathcal{U})
-\frac{KN_0^T(\pa_ux_0\times\pa_vx_0)}{N_0^T(\mathcal{U}\times\mathcal{V})}
\pa_wV^TVN_0^T(V\times\mathcal{U})
-V^T\mathcal{U}N_0^T(\pa_wV\times
V)\frac{KN_0^T(\pa_ux_0\times\pa_vx_0)}
{N_0^T(\mathcal{U}\times\mathcal{V})}$.

From the second equation of R1 we have
$A=N_0^T(\pa_wV\times\mathcal{V})N_0^T(V\times
\mathcal{V})[-N_0^T(\mathcal{U}\times\mathcal{V})N_0\times V
-\mathcal{V}^TV(N_0\times\mathcal{U})\times N_0+V^T\mathcal{U}
(N_0\times\mathcal{V})\times N_0]=0,$
$B=N_0^T(V\times\mathcal{V})N_0^T(\pa_wV\times\mathcal{U})N_0^T(\mathcal{U}\times\mathcal{V})N_0\times
V +[\mathbf{m}^2N_0^T(\mathcal{U}\times\mathcal{V})
-V^T\mathcal{V}N_0^T(V\times\mathcal{U})]N_0^T(\pa_wV\times\mathcal{U})
(N_0\times\mathcal{V})\times N_0+
[-(|V|^2+\mathbf{m}^2)N_0^T(\pa_wV\times\mathcal{V})N_0^T(\mathcal{U}\times\mathcal{V})
+V^T\mathcal{V}N_0^T(\pa_wV\times\mathcal{U})N_0^T(V\times\mathcal{V})]
(N_0\times\mathcal{U})\times N_0+
(\mathbf{m}^2+|V|^2)(N_0^T(\mathcal{U}\times\mathcal{V}))^2\pa_wV,$
$C=[-\frac{N_0^T(\pa_wV\times\pa_vx_0)N_0^T(\pa_wV\times\pa_uV)}
{N_0^T(\mathcal{U}\times\mathcal{V})N_0^T(\pa_wV\times V)}
+\frac{N_0^T(\pa_wV\times\pa_ux_0)N_0^T(\pa_wV\times\pa_vV)}
{N_0^T(\mathcal{U}\times\mathcal{V})N_0^T(\pa_wV\times
V)}]N_0^T(V\times\mathcal{V})
-\frac{N_0^T(\pa_wV\times\mathcal{V})}{N_0^T(\mathcal{U}\times\mathcal{V})}
N_0^T(\mathcal{V}\times\pa_uV)+\frac{N_0^T(\pa_wV\times\mathcal{U})}
{N_0^T(\mathcal{U}\times\mathcal{V})}N_0^T(\pa_vx_0\times\pa_vV)
-(1+\mathbf{m}^2
\frac{KN_0^T(\pa_ux_0\times\pa_vx_0)}{N_0^T(\mathcal{U}\times\mathcal{V})})
N_0^T(\pa_wV\times\mathcal{V})\\
-\frac{KN_0^T(\pa_ux_0\times\pa_vx_0)}{N_0^T(\mathcal{U}\times\mathcal{V})}
\pa_wV^TVN_0^T(V\times\mathcal{V})
-KN_0^T(\pa_ux_0\times\pa_vx_0)V^T\mathcal{V}\frac{N_0^T(\pa_wV\times
V)} {N_0^T(\mathcal{U}\times\mathcal{V})}$.

From the seventh equation of R1 we have
$A=(N_0^T(\mathcal{U}\times\mathcal{V}))^2
N_0^T(V\times\mathcal{V})\pa_wV\times N_0
-[\pa_wV^T\mathcal{V}N_0^T(\mathcal{U}\times\mathcal{V})
+(\mathcal{V}\times N_0)^T(\mathcal{U}\times
N_0)N_0^T(\pa_wV\times\mathcal{V}) +|\mathcal{V}\times
N_0|^2N_0^T(\pa_wV\times\mathcal{U})]
N_0^T(V\times\mathcal{V})(N_0\times \mathcal{U})\times N_0
+(N_0^T(\mathcal{U}\times\mathcal{V}))^2\mathcal{V}^TV\pa_wV+
[N_0^T(V\times\mathcal{U})N_0^T(\pa_wV\times\mathcal{V})+
N_0^T(V\times\mathcal{V})N_0^T(\pa_wV\times\mathcal{U})]
N_0^T(\mathcal{U}\times\mathcal{V})\mathcal{V}\times N_0+
[\pa_wV^T\mathcal{V}N_0^T(\mathcal{U}\times\mathcal{V})N_0^T(V\times\mathcal{U})
+(\mathcal{V}\times N_0)^T(\mathcal{U}\times
N_0)N_0^T(\pa_wV\times\mathcal{U})
N_0^T(V\times\mathcal{V})+|\mathcal{U}\times
N_0|^2N_0^T(\pa_wV\times\mathcal{V})
N_0^T(V\times\mathcal{V})](N_0\times\mathcal{V})\times N_0$,
$B=(N_0^T(\mathcal{U}\times\mathcal{V}))^2
N_0^T(V\times\mathcal{U})\pa_wV\times N_0
+(N_0^T(\mathcal{U}\times\mathcal{V}))^2\mathcal{U}^TV\pa_wV
-[N_0^T(\pa_wV\times\mathcal{V})N_0^T(V\times\mathcal{U})+
N_0^T(\pa_wV\times\mathcal{U})N_0^T(V\times\mathcal{V})]
N_0^T(\mathcal{U}\times\mathcal{V})\mathcal{U}\times N_0
+[N_0^T(\mathcal{U}\times\mathcal{V})\pa_wV^T\mathcal{U}
-|\mathcal{U}\times N_0|^2N_0^T(\pa_wV\times\mathcal{V})
-(\mathcal{U}\times N_0)^T(\mathcal{V}\times
N_0)N_0^T(\pa_wV\times\mathcal{U})]
N_0^T(V\times\mathcal{U})(N_0\times\mathcal{V})\times N_0
+[-\pa_wV^T\mathcal{U}N_0^T(\mathcal{U}\times\mathcal{V})N_0^T(V\times\mathcal{V})
+|\mathcal{V}\times
N_0|^2N_0^T(V\times\mathcal{U})N_0^T(\pa_wV\times\mathcal{U})
+(\mathcal{U}\times N_0)^T(\mathcal{V}\times
N_0)N_0^T(V\times\mathcal{U})
N_0^T(\pa_wV\times\mathcal{V})](N_0\times\mathcal{U})\times N_0$,
$C=(\pa_wV\times N_0)^T[\pa_u[(N_0\times\mathcal{V})\times N_0]+
\pa_v[(N_0\times\mathcal{U})\times N_0]] +(\mathcal{V}\times
N_0)^T[\pa_u[\mathcal{\ti V}(N_0\times\mathcal{U})\times N_0]
+\mathcal{\ti U}\pa_u[(N_0\times\mathcal{V})\times
N_0]-(N_0\times\pa_w\mathcal{U})\times N_0] -(\mathcal{U}\times
N_0)^T[(N_0\times\pa_w\mathcal{V})\times N_0 +\mathcal{\ti
V}\pa_v[(N_0\times\mathcal{U})\times N_0] +\pa_v[\mathcal{\ti
U}(N_0\times\mathcal{V})\times N_0]]
-\frac{1}{\mathbf{m}^2}\pa_wV^T
[\mathcal{U}N_0^T(V\times\mathcal{V})+\mathcal{V}N_0^T(V\times\mathcal{U})]
+\frac{KN_0^T(\pa_ux_0\times\pa_vx_0)}{N_0^T(\mathcal{U}\times\mathcal{V})}
[-(\mathcal{U}\times N_0)^T(\mathcal{V}\times N_0) (\pa_wV\times
V)^TN_0 +|\mathcal{V}\times N_0|^2\mathcal{U}^T(N_0\times
V)\mathcal{\ti U} -|\mathcal{U}\times
N_0|^2\mathcal{V}^T(N_0\times V)\mathcal{\ti V}]
-(\mathcal{U}\times N_0)^T(\mathcal{V}\times
N_0)\frac{(\pa_wV\times V)^TN_0}{\mathbf{m}^2} -|\mathcal{V}\times
N_0|^2\mathcal{\ti
U}\frac{N_0^T(V\times\mathcal{U})}{\mathbf{m}^2}
+|\mathcal{U}\times N_0|^2\mathcal{\ti
V}\frac{N_0^T(V\times\mathcal{V})}{\mathbf{m}^2}$.

From the eighth equation of R1 we have $A=
N_0^T(\pa_wV\times\mathcal{V})N_0^T(V\times\mathcal{V})[N_0^T(\mathcal{U}\times\mathcal{V})
\mathcal{V}\times N_0-|\mathcal{V}\times N_0|^2(N_0\times\mathcal{U})\times N_0
+(\mathcal{V}\times N_0)^T(\mathcal{U}\times N_0)(N_0\times\mathcal{V})\times N_0]=0$,
$B=(N_0^T(\mathcal{U}\times\mathcal{V}))^2N_0^T(V\times\mathcal{V})\pa_wV\times N_0
-N_0^T(\mathcal{U}\times\mathcal{V})N_0^T(\pa_wV\times\mathcal{V})
N_0^T(V\times \mathcal{V})\mathcal{U}\times N_0+
[-\pa_wV^T\mathcal{V}N_0^T(V\times\mathcal{V})N_0^T(\mathcal{U}\times\mathcal{V})
+|\mathcal{V}\times N_0|^2N_0^T(\pa_wV\times\mathcal{V})N_0^T(V\times\mathcal{U})
](N_0\times\mathcal{U})\times N_0
+[\pa_wV^T\mathcal{V}N_0^T(\mathcal{U}\times\mathcal{V})
-(\mathcal{U}\times N_0)^T(\mathcal{V}\times N_0)
N_0^T(\pa_wV\times\mathcal{V})]N_0^T(V\times\mathcal{U})(N_0\times\mathcal{V})\times N_0
+V^T\mathcal{\mathcal{V}}(N_0^T(\mathcal{U}\times\mathcal{V}))^2\pa_wV$,
$C=(\pa_wV\times N_0)^T\pa_v[(N_0\times\mathcal{V})\times N_0]
+(\mathcal{V}\times N_0)^T[\mathcal{\ti V}\pa_u[(N_0\times\mathcal{V})\times N_0]
-(N_0\times\pa_w\mathcal{V})\times N_0]
-(\mathcal{U}\times N_0)^T\pa_v[\mathcal{\ti V}(N_0\times\mathcal{V})\times N_0]
-\frac{KN_0^T(\pa_ux_0\times\pa_vx_0)}{N_0^T(\mathcal{U}\times\mathcal{V})}
N_0^T(\pa_wV\times\mathcal{V})V^T\mathcal{V}$.

From the ninth equation of R1 we have
$A=-(N_0^T(\mathcal{U}\times\mathcal{V}))^2N_0^T(V\times\mathcal{U})\pa_wV\times N_0
-N_0^T(\mathcal{U}\times\mathcal{V})N_0^T(\pa_wV\times\mathcal{U})N_0^T(V\times\mathcal{U})
\mathcal{V}\times N_0-[\pa_wV^T\mathcal{U}N_0^T(V\times\mathcal{U})
N_0^T(\mathcal{U}\times\mathcal{V})
+|\mathcal{U}\times N_0|^2N_0^T(\pa_wV\times\mathcal{U})N_0^T(V\times\mathcal{V})]
(N_0\times\mathcal{V})\times N_0
+[\pa_wV^T\mathcal{U}N_0^T(\mathcal{U}\times\mathcal{V})
+(\mathcal{U}\times N_0)^T(\mathcal{V}\times N_0)
N_0^T(\pa_wV\times\mathcal{U})]N_0^T(V\times\mathcal{V})(N_0\times\mathcal{U})\times N_0
-V^T\mathcal{U}(N_0^T(\mathcal{U}\times\mathcal{V}))^2\pa_wV$,
$B=
N_0^T(\pa_wV\times\mathcal{U})N_0^T(V\times\mathcal{U})[
N_0^T(\mathcal{U}\times\mathcal{V})\mathcal{U}\times N_0
+|\mathcal{U}\times N_0|^2(N_0\times\mathcal{V})\times N_0
-(\mathcal{V}\times N_0)^T(\mathcal{U}\times N_0)(N_0\times\mathcal{U})\times N_0]=0$,
$C=-(\pa_wV\times N_0)^T\pa_u[(N_0\times\mathcal{U})\times N_0]
-(\mathcal{V}\times N_0)^T\pa_u[\mathcal{\ti U}(N_0\times\mathcal{U})\times N_0]
+(\mathcal{U}\times N_0)^T[(N_0\times\pa_w\mathcal{U})\times N_0]
+\mathcal{\ti U}(\mathcal{U}\times N_0)^T\pa_v[(N_0\times\mathcal{U})\times N_0]
+\frac{KN_0^T(\pa_ux_0\times\pa_vx_0)}{N_0^T(\mathcal{U}\times\mathcal{V})}
N_0^T(\pa_wV\times\mathcal{U})\mathcal{U}^TV$.

From the first equation of R2 we have
$A=\frac{1}{2}N_0^T(V\times\mathcal{U})N_0^T(\pa_wV\times\mathcal{V})
N_0^T(\mathcal{U}\times\mathcal{V})\mathcal{U}\times N_0
+[\frac{1}{2}|\mathcal{U}\times N_0|^2N_0^T(V\times\mathcal{U})N_0^T(\pa_wV\times\mathcal{V})
+\frac{1}{2}N_0^T(\pa_wV\times\mathcal{U})V^T\mathcal{U}N_0^T(\mathcal{U}\times\mathcal{V})]
(N_0\times\mathcal{V})\times N_0
-\frac{1}{2}|\mathcal{U}\times N_0|^2N_0^T(V\times\mathcal{V})N_0^T(\pa_wV\times\mathcal{V})
(N_0\times\mathcal{U})\times N_0+\frac{1}{2}V^T\mathcal{U}
(N_0^T(\mathcal{U}\times\mathcal{V}))^2\pa_wV$,
$B=\frac{1}{2}N_0^T(V\times\mathcal{U})N_0^T(\pa_wV\times\mathcal{U})
[-N_0^T(\mathcal{U}\times\mathcal{V})\mathcal{U}\times N_0+
(\mathcal{U}\times N_0)^T(\mathcal{V}\times N_0)(N_0\times\mathcal{U})\times N_0
-|\mathcal{U}\times N_0|^2(N_0\times\mathcal{V})\times N_0]=0$,
$C=\frac{1}{2}\pa_u\mathcal{\ti U}N_0^T(\mathcal{U}\times\mathcal{V})
+\frac{1}{2}\mathcal{\ti V}(\mathcal{U}\times N_0)^T\pa_u[(N_0\times\mathcal{U})\times N_0]
-\frac{1}{2}(\mathcal{U}\times N_0)^T\pa_w\mathcal{U}
-\frac{1}{2}\mathcal{\ti U}(\mathcal{U}\times N_0)^T\pa_v[(N_0\times\mathcal{U})\times N_0]
-\frac{KN_0^T(\pa_ux_0\times\pa_vx_0)}{2N_0^T(\mathcal{U}\times\mathcal{V})}
N_0^T(\pa_wV\times\mathcal{U})V^T\mathcal{U}$.

From the second equation of R2 we have
$A=N_0^T(V\times\mathcal{V})N_0^T(\pa_wV\times\mathcal{V})
[N_0^T(\mathcal{U}\times\mathcal{V})\mathcal{V}\times N_0+
(\mathcal{U}\times N_0)^T(\mathcal{V}\times
N_0)(N_0\times\mathcal{V})\times N_0 -|\mathcal{V}\times
N_0|^2(N_0\times\mathcal{U})\times N_0]=0$,
$B=\frac{N_0^T(\pa_wV\times\mathcal{U})N_0^T(V\times\mathcal{V})}
{\mathbf{m}N_0^T(\mathcal{U}\times\mathcal{V})}\mathcal{V}\times
N_0+(\frac{|\mathcal{V}\times
N_0|^2N_0^T(V\times\mathcal{U})}{\mathbf{m}N_0^T(\mathcal{U}\times\mathcal{V})}
-\frac{(\mathcal{V}\times N_0)^T(\mathcal{U}\times
N_0)N_0^T(V\times\mathcal{V})}{\mathbf{m}N_0^T(\mathcal{U}\times\mathcal{V})})\pa_wV+
\frac{|\mathcal{V}\times
N_0|^2N_0^T(\pa_wV\times\mathcal{U})N_0^T(V\times\mathcal{U})}
{\mathbf{m}N_0^T(\mathcal{U}\times\mathcal{V})^2}(N_0\times\mathcal{V})\times
N_0+[-\frac{3}{2}\frac{|\mathcal{V}\times
N_0|^2N_0^T(\pa_wV\times\mathcal{U})N_0^T(V\times\mathcal{V})}
{\mathbf{m}N_0^T(\mathbf{U}\times\mathcal{V})^2}+|\mathcal{V}\times
N_0|^2\frac{N_0^T(\pa_wV\times
V)}{2\mathbf{m}N_0^T(\mathcal{U}\times\mathcal{V})}+\frac{(\mathcal{V}\times
N_0)^T(\mathcal{U}\times
N_0)N_0^T(V\times\mathcal{V})N_0^T(\pa_wV\times
\mathcal{V})}{\mathbf{m}N_0^T(\mathcal{U}\times\mathcal{V})^2}\\
-\frac{|\mathcal{V}\times N_0|^2N_0^T(\pa_wV\times
\mathcal{V})N_0^T(V\times\mathcal{U})}{2\mathbf{m}N_0^T(\mathcal{U}\times\mathcal{V})^2}]
(N_0\times\mathcal{U})\times N_0$, $C=(\mathcal{V}\times
N_0)^T\pa_w\mathcal{V}-\pa_v\ti{\mathcal{V}}N_0^T(\mathcal{U}\times\mathcal{V})
+\ti{\mathcal{U}}(\mathcal{V}\times
N_0)^T\pa_v[(N_0\times\mathcal{V})\times
N_0]+\ti{\mathcal{V}}KN_0^T(\pa_ux_0\times\pa_vx_0)\mathcal{V}^TV
-\ti{\mathcal{V}}(\mathcal{V}\times
N_0)^T\pa_u[(N_0\times\mathcal{V})\times N_0]=0$.

From the third equation of R2 we have
$A=-N_0^T(\mathcal{U}\times\mathcal{V})N_0^T(\pa_wV\times\mathcal{V})
[N_0^T(V\times\mathcal{U})\mathcal{V}\times N_0+
N_0^T(V\times\mathcal{V})\mathcal{U}\times N_0]
+2N_0^T(\pa_wV\times\mathcal{V})N_0^T(V\times\mathcal{V})
(\mathcal{V}\times N_0)^T(\mathcal{U}\times
N_0)(N_0\times\mathcal{U})\times N_0
-[N_0^T(V\times\mathcal{U})N_0^T(\mathcal{U}\times\mathcal{V})\pa_wV^T\mathcal{V}
+N_0^T(\pa_wV\times\mathcal{U})N_0^T(V\times\mathcal{V})
(\mathcal{V}\times N_0)^T(\mathcal{U}\times N_0)
+N_0^T(\pa_wV\times\mathcal{V})N_0^T(V\times\mathcal{V})
|\mathcal{U}\times N_0|^2](N_0\times\mathcal{V})\times N_0
-V^T\mathcal{V}(N_0^T(\mathcal{U}\times\mathcal{V}))^2\pa_wV$,
$B=N_0^T(\pa_wV\times\mathcal{U})N_0^T(\mathcal{U}\times\mathcal{V})[N_0^T(V\times\mathcal{V})
\mathcal{U}\times N_0+N_0^T(V\times\mathcal{U}) \mathcal{V}\times
N_0]
+[N_0^T(V\times\mathcal{V})N_0^T(\mathcal{U}\times\mathcal{V})\pa_wV^T\mathcal{U}
-N_0^T(\pa_wV\times\mathcal{V})N_0^T(V\times\mathcal{U})
(\mathcal{V}\times N_0)^T(\mathcal{U}\times N_0)
-N_0^T(\pa_wV\times\mathcal{U})N_0^T(V\times\mathcal{U})
|\mathcal{V}\times N_0|^2](N_0\times\mathcal{U})\times N_0
+2N_0^T(\pa_wV\times\mathcal{U})N_0^T(V\times\mathcal{U})
(\mathcal{V}\times N_0)^T(\mathcal{U}\times
N_0)(N_0\times\mathcal{V})\times N_0
-(N_0^T(\mathcal{U}\times\mathcal{V}))^2V^T\mathcal{U}\pa_wV$,
$C=(\mathcal{V}\times N_0)^T[(N_0\times\pa_w\mathcal{U})\times N_0
-\pa_v\mathcal{\ti U}(N_0\times\mathcal{U})\times N_0
+\mathcal{\ti U}\pa_v[(N_0\times\mathcal{U})\times N_0]
-\mathcal{\ti V}\pa_u[(N_0\times\mathcal{U})\times N_0]
-2\pa_u\mathcal{\ti V}(N_0\times\mathcal{U})\times N_0]
+(\mathcal{U}\times N_0)^T[\mathcal{\ti
U}\pa_v[(N_0\times\mathcal{V})\times N_0] -\pa_u[\mathcal{\ti
V}(N_0\times\mathcal{V})\times N_0]
+(N_0\times\pa_w\mathcal{V})\times N_0]
+\frac{KN_0^T(\pa_ux_0\times\pa_vx_0)}{N_0^T(\mathcal{U}\times\mathcal{V})}
[N_0^T(\pa_wV\times V) (\mathcal{U}\times N_0)^T(\mathcal{V}\times
N_0) -|\mathcal{V}\times N_0|^2\mathcal{\ti
U}N_0^T(V\times\mathcal{U}) +|\mathcal{U}\times N_0|^2\mathcal{\ti
V}N_0^T(V\times\mathcal{V})]$.

Note that by the exchange of coordinates $u\leftrightarrow v$ we
have $\mathcal{U}\leftrightarrow\mathcal{V}$ and $C=0$ for the
fifth equation of R1 is exchanged with $C=0$ for the second
equation of R1, $A=0$ of the fifth equation of R1 is exchanged
with $B=0$ of the second equation of R1, $A$ of the seventh
equation of R1 is exchanged with $B$ of the seventh equation of
R1, $C=0$ for the seventh equation of R1 is exchanged with itself,
$B=0$ for the eighth equation of R1 is exchanged with $-A=0$ for
the ninth equation of R1, $C=0$ for the eighth equation of R1 is
exchanged with $-C=0$ for the ninth equation of R1, $A=0$ for the
first equation of R2 is exchanged with $\hat A=0$, where $\hat
A=-\frac{1}{2}N_0^T(V\times\mathcal{V})N_0^T(\pa_wV\times\mathcal{U})
N_0^T(\mathcal{U}\times\mathcal{V})\mathcal{V}\times N_0
+[\frac{1}{2}|\mathcal{V}\times
N_0|^2N_0^T(V\times\mathcal{V})N_0^T(\pa_wV\times\mathcal{U})
-\frac{1}{2}N_0^T(\pa_wV\times\mathcal{V})V^T\mathcal{V}N_0^T(\mathcal{U}\times\mathcal{V})]
(N_0\times\mathcal{U})\times N_0 -\frac{1}{2}|\mathcal{V}\times
N_0|^2N_0^T(V\times\mathcal{U})N_0^T(\pa_wV\times\mathcal{U})
(N_0\times\mathcal{V})\times N_0+\frac{1}{2}V^T\mathcal{V}
(N_0^T(\mathcal{U}\times\mathcal{V}))^2\pa_wV$, $C=0$ for the
first equation of R2 is exchanged with $\hat C=0$, where $\hat
C=\frac{1}{2}\pa_v\mathcal{\ti
V}N_0^T(\mathcal{U}\times\mathcal{V}) -\frac{1}{2}\mathcal{\ti
U}(\mathcal{V}\times N_0)^T\pa_v[(N_0\times\mathcal{V})\times N_0]
-\frac{1}{2}(\mathcal{V}\times N_0)^T\pa_w\mathcal{V}
+\frac{1}{2}\mathcal{\ti V}(\mathcal{V}\times
N_0)^T\pa_u[(N_0\times\mathcal{V})\times N_0]
-\frac{KN_0^T(\pa_ux_0\times\pa_vx_0)}{2N_0^T(\mathcal{U}\times\mathcal{V})}
N_0^T(\pa_wV\times\mathcal{V})V^T\mathcal{V}$, $B=0$ for the
second equation of R2 is exchanged with $\hat B=0$, where $\hat
B=-\frac{N_0^T(\pa_wV\times\mathcal{V})N_0^T(V\times\mathcal{U})}
{\mathbf{m}N_0^T(\mathcal{U}\times\mathcal{V})}\mathcal{U}\times
N_0-(\frac{|\mathcal{U}\times
N_0|^2N_0^T(V\times\mathcal{V})}{\mathbf{m}N_0^T(\mathcal{U}\times\mathcal{V})}
-\frac{(\mathcal{V}\times N_0)^T(\mathcal{U}\times
N_0)N_0^T(V\times\mathcal{U})}{\mathbf{m}N_0^T(\mathcal{U}\times\mathcal{V})})\pa_wV+
\frac{|\mathcal{U}\times
N_0|^2N_0^T(\pa_wV\times\mathcal{V})N_0^T(V\times\mathcal{V})}
{\mathbf{m}N_0^T(\mathcal{U}\times\mathcal{V})^2}(N_0\times\mathcal{U})\times
N_0+[-\frac{3}{2}\frac{|\mathcal{U}\times
N_0|^2N_0^T(\pa_wV\times\mathcal{V})N_0^T(V\times\mathcal{U})}
{\mathbf{m}N_0^T(\mathbf{U}\times\mathcal{V})^2}-|\mathcal{U}\times
N_0|^2\frac{N_0^T(\pa_wV\times
V)}{2\mathbf{m}N_0^T(\mathcal{U}\times\mathcal{V})}+\frac{(\mathcal{V}\times
N_0)^T(\mathcal{U}\times
N_0)N_0^T(V\times\mathcal{U})N_0^T(\pa_wV\times
\mathcal{U})}{\mathbf{m}N_0^T(\mathcal{U}\times\mathcal{V})^2}
-\frac{|\mathcal{U}\times N_0|^2N_0^T(\pa_wV\times
\mathcal{U})N_0^T(V\times\mathcal{V})}{2\mathbf{m}N_0^T(\mathcal{U}\times\mathcal{V})^2}]
(N_0\times\mathcal{V})\times N_0$, $C=0$ for the second equation
of R2 is exchanged with $\hat C=0$, where $\hat
C=(\mathcal{U}\times
N_0)^T\pa_w\mathcal{U}-\pa_u\ti{\mathcal{U}}N_0^T(\mathcal{U}\times\mathcal{V})
-\ti{\mathcal{V}}(\mathcal{U}\times
N_0)^T\pa_u[(N_0\times\mathcal{U})\times
N_0]+\ti{\mathcal{U}}KN_0^T(\pa_ux_0\times\pa_vx_0)\mathcal{U}^TV
+\ti{\mathcal{U}}(\mathcal{U}\times
N_0)^T\pa_v[(N_0\times\mathcal{U})\times N_0]$, $A$ of the third
equation of R2 is exchanged with $B$ of the third equation of R2
and $C=0$ for the third equation of R3 is exchanged with $\hat
C=0$, where $\hat C=(\mathcal{U}\times
N_0)^T[(N_0\times\pa_w\mathcal{V})\times N_0 +\pa_u\mathcal{\ti
V}(N_0\times\mathcal{V})\times N_0 -\mathcal{\ti
V}\pa_u[(N_0\times\mathcal{V})\times N_0] +\mathcal{\ti
U}\pa_v[(N_0\times\mathcal{V})\times N_0] +2\pa_v\mathcal{\ti
U}(N_0\times\mathcal{V})\times N_0] +(\mathcal{V}\times
N_0)^T[-\mathcal{\ti V}\pa_u[(N_0\times\mathcal{U})\times N_0]
+\pa_v[\mathcal{\ti U}(N_0\times\mathcal{U})\times N_0]
+(N_0\times\pa_w\mathcal{U})\times N_0]
+\frac{KN_0^T(\pa_ux_0\times\pa_vx_0)}{N_0^T(\mathcal{U}\times\mathcal{V})}
[N_0^T(\pa_wV\times V) (\mathcal{U}\times N_0)^T(\mathcal{V}\times
N_0) +|\mathcal{U}\times N_0|^2\mathcal{\ti
V}N_0^T(V\times\mathcal{V}) -|\mathcal{V}\times N_0|^2\mathcal{\ti
U}N_0^T(V\times\mathcal{U})]$.

Using
$$W=\frac{N_0^T(W\times\mathcal{V})}{N_0^T(\mathcal{U}\times\mathcal{V})}
(N_0\times\mathcal{U})\times N_0
-\frac{N_0^T(W\times\mathcal{U})}{N_0^T(\mathcal{U}\times\mathcal{V})}
(N_0\times\mathcal{V})\times N_0+W^TN_0N_0$$ for any vector field $W$
we get
\begin{eqnarray}\label{eq:tang}
N_0\times V=-\frac{\mathcal{V}^TV}{N_0^T(\mathcal{U}\times\mathcal{V})}
(N_0\times\mathcal{U})\times N_0+\frac{V^T\mathcal{U}}{N_0^T(\mathcal{U}\times\mathcal{V})}
(N_0\times\mathcal{V})\times N_0,\nonumber\\
\mathcal{V}\times N_0=\frac{|\mathcal{V}\times N_0|^2}{N_0^T(\mathcal{U}\times\mathcal{V})}
(N_0\times\mathcal{U})\times N_0
-\frac{(\mathcal{V}\times N_0)^T(\mathcal{U}\times N_0)}{N_0^T(\mathcal{U}\times\mathcal{V})}
(N_0\times\mathcal{V})\times N_0,\nonumber\\
\mathcal{U}\times N_0=
\frac{(\mathcal{V}\times N_0)^T(\mathcal{U}\times N_0)}{N_0^T(\mathcal{U}\times\mathcal{V})}
(N_0\times\mathcal{U})\times N_0
-\frac{|\mathcal{U}\times N_0|^2}{N_0^T(\mathcal{U}\times\mathcal{V})}
(N_0\times\mathcal{V})\times N_0,\nonumber\\
\pa_wV=\frac{N_0^T(\pa_wV\times\mathcal{V})}{N_0^T(\mathcal{U}\times\mathcal{V})}
(N_0\times\mathcal{U})\times N_0
-\frac{N_0^T(\pa_wV\times\mathcal{U})}{N_0^T(\mathcal{U}\times\mathcal{V})}
(N_0\times\mathcal{V})\times N_0,\nonumber\\
\pa_wV\times N_0=\frac{\pa_wV^T\mathcal{V}}{N_0^T(\mathcal{U}\times\mathcal{V})}
(N_0\times\mathcal{U})\times N_0
-\frac{\pa_wV^T\mathcal{U}}{N_0^T(\mathcal{U}\times\mathcal{V})}
(N_0\times\mathcal{V})\times N_0
\end{eqnarray}
and the relations $A=0$ for the fifth equation of R1, $B=0$ for the second equation of R1,
$B=0$ for the eighth equation of R1, $A=0$ for the ninth equation of R1, $A=0$ and $\hat A=0$
for the first equation of R2 and $B=0$ and $\hat B=0$ for the second equation of R2 become
linear combinations of the linearly independent vector fields
$(N_0\times\mathcal{U})\times N_0,\ (N_0\times\mathcal{V})\times N_0$ being $0$.

Using (\ref{eq:tang}) the relations $A=0$ for the fifth equation
of R1, $B=0$ for the second equation of R1, $B=0$ for the eighth
equation of R1, $A=0$ for the ninth equation of R1, $A=0$ and
$\hat A=0$ for the first equation of R2 and $B=0$ and $\hat B=0$
for the second equation of R2 are identically satisfied.

Using again (\ref{eq:tang}) for the seventh equation of R1 and for the third equation of
R2 we have $A=0,\ B=0$.

With $V=V_1\pa_ux_0+V_2\pa_vx_0$ the relations $C=0$ for the fifth
and second equation of R1 depend on $V_1,V_2$ and their first
derivatives; the remaining relations $C=0$ and $\hat C=0$ depend
a-priori linearly on the second derivatives of $V_1,V_2$ and on
their first derivatives and on $V_1,V_2$; we shall see that the
dependence on the second derivatives of $V_1,V_2$ disappears.

From $C=0$ for the fifth relation of R1 we get

$-[\pa_wV_1[\pa_wV_1(\pa_uV_2+V_1\Ga_{11}^2+V_2\Ga_{12}^2)
-\pa_wV_2(\pa_uV_1+V_1\Ga_{11}^1+V_2\Ga_{12}^1)]
+\pa_wV_2[\pa_wV_1(\pa_vV_2+V_1\Ga_{12}^2+V_2\Ga_{22}^2)
-\pa_wV_2(\pa_vV_1+V_1\Ga_{12}^1+V_2\Ga_{22}^1)]][V_1(\pa_uV_2+V_1\Ga_{11}^2+V_2\Ga_{12}^2)
-V_2(\pa_uV_1+V_1\Ga_{11}^1+V_2\Ga_{12}^1+1)]-[\pa_wV_1(\pa_vV_2+V_1\Ga_{12}^2+V_2\Ga_{22}^2+1)
-\pa_wV_2(\pa_vV_1+V_1\Ga_{12}^1+V_2\Ga_{22}^1)](\pa_uV_2+V_1\Ga_{11}^2+V_2\Ga_{12}^2)
(\pa_wV_1V_2-\pa_wV_2V_1)+[\pa_wV_1(\pa_uV_2+V_1\Ga_{11}^2+V_2\Ga_{12}^2)
-\pa_wV_2(\pa_uV_1+V_1\Ga_{11}^1+V_2\Ga_{12}^1+1)](\pa_wV_1V_2-\pa_wV_2V_1)
[(\pa_uV_1+V_1\Ga_{11}^1+V_2\Ga_{12}^1+1)(\pa_vV_2+V_1\Ga_{12}^2+V_2\Ga_{22}^2)
-(\pa_uV_2+V_1\Ga_{11}^2+V_2\Ga_{12}^2)(\pa_vV_1+V_1\Ga_{12}^1+V_2\Ga_{22}^1)]
-[(\pa_wV_1V_2-\pa_wV_2V_1)[(\pa_uV_1+V_1\Ga_{11}^1+V_2\Ga_{12}^1+1)
(\pa_vV_2+V_1\Ga_{12}^2+V_2\Ga_{22}^2+1)-(\pa_uV_2+V_1\Ga_{11}^2+V_2\Ga_{12}^2)
(\pa_vV_1+V_1\Ga_{12}^1+V_2\Ga_{22}^1)]-|V|^2K(\pa_wV_1V_2-\pa_wV_2V_1)
-V_1[\pa_wV_1(\pa_uV_2+V_1\Ga_{11}^2+V_2\Ga_{12}^2)
-\pa_wV_2(\pa_uV_1+V_1\Ga_{11}^1+V_2\Ga_{12}^1+1)]
-V_2[\pa_wV_1(\pa_vV_2+V_1\Ga_{12}^2+V_2\Ga_{22}^2+1)
-\pa_wV_2(\pa_vV_1+V_1\Ga_{12}^1+V_2\Ga_{22}^1)]]
[\pa_wV_1(\pa_uV_2+V_1\Ga_{11}^2+V_2\Ga_{12}^2)
-\pa_wV_2(\pa_uV_1+V_1\Ga_{11}^1+V_2\Ga_{12}^1+1)]
-K(\pa_wV_1V_2-\pa_wV_2V_1)\pa_wV^TV[V_1(\pa_uV_2+V_1\Ga_{11}^2+V_2\Ga_{12}^2)
-V_2(\pa_uV_1+V_1\Ga_{11}^1+V_2\Ga_{12}^1+1)]-K(\pa_wV_1V_2-\pa_wV_2V_1)^2
[V^T\pa_ux_0(\pa_uV_1+V_1\Ga_{11}^1+V_2\Ga_{12}^1+1)
+V^T\pa_vx_0(\pa_uV_2+V_1\Ga_{11}^2+V_2\Ga_{12}^2)]=0.$ This is
identically satisfied (does not introduce any conditions).

From $C=0$ for the second relation of R1 we get

$-[\pa_wV_1[\pa_wV_1(\pa_uV_2+V_1\Ga_{11}^2+V_2\Ga_{12}^2)
-\pa_wV_2(\pa_uV_1+V_1\Ga_{11}^1+V_2\Ga_{12}^1)]
+\pa_wV_2[\pa_wV_1(\pa_vV_2+V_1\Ga_{12}^2+V_2\Ga_{22}^2)
-\pa_wV_2(\pa_vV_1+V_1\Ga_{12}^1+V_2\Ga_{22}^1)]]
[V_1(\pa_vV_2+V_1\Ga_{12}^2+V_2\Ga_{22}^2+1)-V_2(\pa_vV_1+V_1\Ga_{12}^1+V_2\Ga_{22}^1)]
-(\pa_wV_1V_2-\pa_wV_2V_1)[\pa_wV_1(\pa_vV_2+V_1\Ga_{12}^2+V_2\Ga_{22}^2+1)
-\pa_wV_2(\pa_vV_1+V_1\Ga_{12}^1+V_2\Ga_{22}^1)]
[(\pa_vV_1+V_1\Ga_{12}^1+V_2\Ga_{22}^1)(\pa_uV_2+V_1\Ga_{11}^2+V_2\Ga_{12}^2)-
(\pa_vV_2+V_1\Ga_{12}^2+V_2\Ga_{22}^2+1)(\pa_uV_1+V_1\Ga_{11}^1+V_2\Ga_{12}^1)]
-(\pa_wV_1V_2-\pa_wV_2V_1)(\pa_vV_1+V_1\Ga_{12}^1+V_2\Ga_{22}^1)
[\pa_wV_1(\pa_uV_2+V_1\Ga_{11}^2+V_2\Ga_{12}^2)-\pa_wV_2(\pa_uV_1+V_1\Ga_{11}^1+V_2\Ga_{12}^1+1)]
-[(\pa_wV_1V_2-\pa_wV_2V_1)[(\pa_uV_1+V_1\Ga_{11}^1+V_2\Ga_{12}^1+1)
(\pa_vV_2+V_1\Ga_{12}^2+V_2\Ga_{22}^2+1)-(\pa_uV_2+V_1\Ga_{11}^2+V_2\Ga_{12}^2)
(\pa_vV_1+V_1\Ga_{12}^1+V_2\Ga_{22}^1)]-|V|^2K(\pa_wV_1V_2-\pa_wV_2V_1)
-V_1[\pa_wV_1(\pa_uV_2+V_1\Ga_{11}^2+V_2\Ga_{12}^2)
-\pa_wV_2(\pa_uV_1+V_1\Ga_{11}^1+V_2\Ga_{12}^1+1)]
-V_2[\pa_wV_1(\pa_vV_2+V_1\Ga_{12}^2+V_2\Ga_{22}^2+1)
-\pa_wV_2(\pa_vV_1+V_1\Ga_{12}^1+V_2\Ga_{22}^1)]]
[\pa_wV_1(\pa_vV_2+V_1\Ga_{12}^2+V_2\Ga_{22}^2+1)-\pa_wV_2(\pa_vV_1+V_1\Ga_{12}^1+V_2\Ga_{22}^1)]
-K(\pa_wV_1V_2\\-\pa_wV_2V_1)\pa_wV^TV[V_1(\pa_vV_2+V_1\Ga_{12}^2+V_2\Ga_{22}^2+1)
-V_2(\pa_vV_1+V_1\Ga_{12}^1+V_2\Ga_{22}^1)]-K(\pa_wV_1V_2-\pa_wV_2V_1)^2
[V^T\pa_ux_0(\pa_vV_1+V_1\Ga_{12}^1+V_2\Ga_{22}^1)
+V^T\pa_vx_0(\pa_vV_2+V_1\Ga_{12}^2+V_2\Ga_{22}^2+1)]=0.$ This is
identically satisfied (does not introduce any conditions).

From $C=0$ for the seventh relation of R1 we get

$\pa_wV_2[\Ga_{11}^1(\pa_vV_1+V_1\Ga_{12}^1+V_2\Ga_{22}^1)
+\Ga_{12}^1(\pa_vV_2+V_1\Ga_{12}^2+V_2\Ga_{22}^2+1)]
-\pa_wV_1[\Ga_{11}^2(\pa_vV_1+V_1\Ga_{12}^1+V_2\Ga_{22}^1)
+\Ga_{12}^2(\pa_vV_2+V_1\Ga_{12}^2+V_2\Ga_{22}^2+1)]
+\pa_wV_2[\Ga_{12}^1(\pa_uV_1+V_1\Ga_{11}^1+V_2\Ga_{12}^1+1)
+\Ga_{22}^1(\pa_uV_2+V_1\Ga_{11}^2+V_2\Ga_{12}^2)]
-\pa_wV_1[\Ga_{12}^2(\pa_uV_1+V_1\Ga_{11}^1+V_2\Ga_{12}^1+1)
+\Ga_{22}^2(\pa_uV_2+V_1\Ga_{11}^2+V_2\Ga_{12}^2)]\\
+\frac{\pa_wV_1(\pa_vV_2+V_1\Ga_{12}^2+V_2\Ga_{22}^2+1)
-\pa_wV_2(\pa_vV_1+V_1\Ga_{12}^1+V_2\Ga_{22}^1)}{(\pa_uV_1+V_1\Ga_{11}^1+V_2\Ga_{12}^1+1)
(\pa_vV_2+V_1\Ga_{12}^2+V_2\Ga_{22}^2+1)-(\pa_uV_2+V_1\Ga_{11}^2+V_2\Ga_{12}^2)
(\pa_vV_1+V_1\Ga_{12}^1+V_2\Ga_{22}^1)}[-\pa_u(V_1\Ga_{12}^2\\+V_2\Ga_{22}^2)
(\pa_uV_1+V_1\Ga_{11}^1+V_2\Ga_{12}^1+1)+\pa_u(V_1\Ga_{12}^1+V_2\Ga_{22}^1)
(\pa_uV_2+V_1\Ga_{11}^2+V_2\Ga_{12}^2)+[\Ga_{11}^1(\pa_uV_1+V_1\Ga_{11}^1+V_2\Ga_{12}^1+1)
+\Ga_{12}^1(\pa_uV_2+V_1\Ga_{11}^2+V_2\Ga_{12}^2)](\pa_vV_2+V_1\Ga_{12}^2+V_2\Ga_{22}^2+1)
-[\Ga_{11}^2(\pa_uV_1+V_1\Ga_{11}^1+V_2\Ga_{12}^1+1)
+\Ga_{12}^2(\pa_uV_2+V_1\Ga_{11}^2+V_2\Ga_{12}^2)](\pa_vV_1+V_1\Ga_{12}^1+V_2\Ga_{22}^1)
-[\Ga_{12}^1(\pa_uV_1+V_1\Ga_{11}^1+V_2\Ga_{12}^1+1)
+\Ga_{22}^1(\pa_uV_2+V_1\Ga_{11}^2+V_2\Ga_{12}^2)](\pa_uV_2+V_1\Ga_{11}^2+V_2\Ga_{12}^2)
+[\Ga_{12}^2(\pa_uV_1+V_1\Ga_{11}^1+V_2\Ga_{12}^1+1)
+\Ga_{22}^2(\pa_uV_2+V_1\Ga_{11}^2+V_2\Ga_{12}^2)](\pa_uV_1+V_1\Ga_{11}^1+V_2\Ga_{12}^1+1)
-\pa_v(V_1\Ga_{11}^1+V_2\Ga_{12}^1)(\pa_uV_2+V_1\Ga_{11}^2+V_2\Ga_{12}^2)
+\pa_v(V_1\Ga_{11}^2+V_2\Ga_{12}^2)(\pa_uV_1+V_1\Ga_{11}^1+V_2\Ga_{12}^1+1)]\\
+\frac{\pa_wV_1(\pa_uV_2+V_1\Ga_{11}^2+V_2\Ga_{12}^2)
-\pa_wV_2(\pa_uV_1+V_1\Ga_{11}^1+V_2\Ga_{12}^1+1)}{(\pa_uV_1+V_1\Ga_{11}^1+V_2\Ga_{12}^1+1)
(\pa_vV_2+V_1\Ga_{12}^2+V_2\Ga_{22}^2+1)-(\pa_uV_2+V_1\Ga_{11}^2+V_2\Ga_{12}^2)
(\pa_vV_1+V_1\Ga_{12}^1+V_2\Ga_{22}^1)}[-\pa_v(V_1\Ga_{11}^1\\+V_2\Ga_{12}^1)
(\pa_vV_2+V_1\Ga_{12}^2+V_2\Ga_{22}^2+1)+\pa_v(V_1\Ga_{11}^2+V_2\Ga_{12}^2)
(\pa_vV_1+V_1\Ga_{12}^1+V_2\Ga_{22}^1)-[\Ga_{12}^1(\pa_vV_1+V_1\Ga_{12}^1+V_2\Ga_{22}^1)
+\Ga_{22}^1(\pa_vV_2+V_1\Ga_{12}^2+V_2\Ga_{22}^2+1)](\pa_uV_2+V_1\Ga_{11}^2+V_2\Ga_{12}^2)
+[\Ga_{12}^2(\pa_vV_1+V_1\Ga_{12}^1+V_2\Ga_{22}^1)
+\Ga_{22}^2(\pa_vV_2+V_1\Ga_{12}^2+V_2\Ga_{22}^2+1)](\pa_uV_1+V_1\Ga_{11}^1+V_2\Ga_{12}^1+1)
+[\Ga_{11}^1(\pa_vV_1+V_1\Ga_{12}^1+V_2\Ga_{22}^1)
+\Ga_{12}^1(\pa_vV_2+V_1\Ga_{12}^2+V_2\Ga_{22}^2+1)](\pa_vV_2+V_1\Ga_{12}^2+V_2\Ga_{22}^2+1)
-[\Ga_{11}^2(\pa_vV_1+V_1\Ga_{12}^1+V_2\Ga_{22}^1)
+\Ga_{12}^2(\pa_vV_2+V_1\Ga_{12}^2+V_2\Ga_{22}^2+1)](\pa_vV_1+V_1\Ga_{12}^1+V_2\Ga_{22}^1)
+\pa_u(V_1\Ga_{12}^1+V_2\Ga_{22}^1)(\pa_vV_2+V_1\Ga_{12}^2+V_2\Ga_{22}^2+1)
-\pa_u(V_1\Ga_{12}^2+V_2\Ga_{22}^2)(\pa_vV_1+V_1\Ga_{12}^1+V_2\Ga_{22}^1)]
-(\pa_wV_1\Ga_{11}^1+\pa_wV_2\Ga_{12}^1)(\pa_vV_2+V_1\Ga_{12}^2+V_2\Ga_{22}^2+1)
+(\pa_wV_1\Ga_{11}^2+\pa_wV_2\Ga_{12}^2)(\pa_vV_1+V_1\Ga_{12}^1+V_2\Ga_{22}^1)
-(\pa_wV_1\Ga_{12}^1+\pa_wV_2\Ga_{22}^1)(\pa_uV_2+V_1\Ga_{11}^2+V_2\Ga_{12}^2)
+(\pa_wV_1\Ga_{12}^2+\pa_wV_2\Ga_{22}^2)(\pa_uV_1+V_1\Ga_{11}^1+V_2\Ga_{12}^1+1)
+\frac{K}{(\pa_uV_1+V_1\Ga_{11}^1+V_2\Ga_{12}^1+1)
(\pa_vV_2+V_1\Ga_{12}^2+V_2\Ga_{22}^2+1)-(\pa_uV_2+V_1\Ga_{11}^2+V_2\Ga_{12}^2)
(\pa_vV_1+V_1\Ga_{12}^1+V_2\Ga_{22}^1)}[-(\pa_wV_1V_2\\-\pa_wV_2V_1)
[(\pa_uV_1+V_1\Ga_{11}^1+V_2\Ga_{12}^1+1)\pa_ux_0
+(\pa_uV_2+V_1\Ga_{11}^2+V_2\Ga_{12}^2)\pa_vx_0]^T
[(\pa_vV_1+V_1\Ga_{12}^1+V_2\Ga_{22}^1)\pa_ux_0+(\pa_vV_2+V_1\Ga_{12}^2+V_2\Ga_{22}^2+1)\pa_vx_0]
\\+\frac{\pa_wV_1(\pa_uV_2+V_1\Ga_{11}^2+V_2\Ga_{12}^2)
-\pa_wV_2(\pa_uV_1+V_1\Ga_{11}^1+V_2\Ga_{12}^1+1)}{(\pa_uV_1+V_1\Ga_{11}^1+V_2\Ga_{12}^1+1)
(\pa_vV_2+V_1\Ga_{12}^2+V_2\Ga_{22}^2+1)-(\pa_uV_2+V_1\Ga_{11}^2+V_2\Ga_{12}^2)
(\pa_vV_1+V_1\Ga_{12}^1+V_2\Ga_{22}^1)}|(\pa_vV_1+V_1\Ga_{12}^1+V_2\Ga_{22}^1)\pa_ux_0
+(\pa_vV_2+V_1\Ga_{12}^2+V_2\Ga_{22}^2+1)\pa_vx_0|^2[V_1(\pa_uV_2+V_1\Ga_{11}^2+V_2\Ga_{12}^2)
-V_2(\pa_uV_1+V_1\Ga_{11}^1+V_2\Ga_{12}^1+1)]
-\frac{\pa_wV_1(\pa_vV_2+V_1\Ga_{12}^2+V_2\Ga_{22}^2+1)
-\pa_wV_2(\pa_vV_1+V_1\Ga_{12}^1+V_2\Ga_{22}^1)}{(\pa_uV_1+V_1\Ga_{11}^1+V_2\Ga_{12}^1+1)
(\pa_vV_2+V_1\Ga_{12}^2+V_2\Ga_{22}^2+1)-(\pa_uV_2+V_1\Ga_{11}^2+V_2\Ga_{12}^2)
(\pa_vV_1+V_1\Ga_{12}^1+V_2\Ga_{22}^1)}|(\pa_uV_1+V_1\Ga_{11}^1+V_2\Ga_{12}^1+1)\pa_ux_0
+(\pa_uV_2+V_1\Ga_{11}^2+V_2\Ga_{12}^2)\pa_vx_0|^2[V_1(\pa_vV_2+V_1\Ga_{12}^2+V_2\Ga_{22}^2+1)
-V_2(\pa_vV_1+V_1\Ga_{12}^1+V_2\Ga_{22}^1)]]
+K(\pa_wV_1V_2-\pa_wV_2V_1)/[V_1[\pa_wV_1(\pa_uV_2+V_1\Ga_{11}^2+V_2\Ga_{12}^2)
-\pa_wV_2(\pa_uV_1+V_1\Ga_{11}^1+V_2\Ga_{12}^1+1)]
+V_2[\pa_wV_1(\pa_vV_2+V_1\Ga_{12}^2+V_2\Ga_{22}^2+1)
-\pa_wV_2(\pa_vV_1+V_1\Ga_{12}^1+V_2\Ga_{22}^1)]+|V|^2K(\pa_wV_1V_2-\pa_wV_2V_1)]
\{[(\pa_wV_1\pa_ux_0+\pa_wV_2\pa_vx_0)^T[(\pa_uV_1+V_1\Ga_{11}^1+V_2\Ga_{12}^1+1)\pa_ux_0
+(\pa_uV_2+V_1\Ga_{11}^2\\+V_2\Ga_{12}^2)\pa_vx_0]
[V_1(\pa_vV_2+V_1\Ga_{12}^2+V_2\Ga_{22}^2+1)
-V_2(\pa_vV_1+V_1\Ga_{12}^1+V_2\Ga_{22}^1)]+(\pa_wV_1\pa_ux_0\\+\pa_wV_2\pa_vx_0)^T
[(\pa_vV_1+V_1\Ga_{12}^1+V_2\Ga_{22}^1)\pa_ux_0
+(\pa_vV_2+V_1\Ga_{12}^2+V_2\Ga_{22}^2+1)\pa_vx_0][V_1(\pa_uV_2+V_1\Ga_{11}^2+V_2\Ga_{12}^2)
-V_2(\pa_uV_1+V_1\Ga_{11}^1+V_2\Ga_{12}^1+1)]
+[(\pa_uV_1+V_1\Ga_{11}^1+V_2\Ga_{12}^1+1)\pa_ux_0
+(\pa_uV_2+V_1\Ga_{11}^2+V_2\Ga_{12}^2)\pa_vx_0]^T
[(\pa_vV_1+V_1\Ga_{12}^1+V_2\Ga_{22}^1)\pa_ux_0
+(\pa_vV_2+V_1\Ga_{12}^2+V_2\Ga_{22}^2+1)\pa_vx_0](\pa_wV_1V_2-\pa_wV_2V_1)
+|(\pa_vV_1+V_1\Ga_{12}^1+V_2\Ga_{22}^1)\pa_ux_0
+(\pa_vV_2+V_1\Ga_{12}^2+V_2\Ga_{22}^2+1)\pa_vx_0|^2[V_1(\pa_uV_2+V_1\Ga_{11}^2+V_2\Ga_{12}^2)
-V_2(\pa_uV_1+V_1\Ga_{11}^1+V_2\Ga_{12}^1+1)]
\frac{\pa_wV_1(\pa_uV_2+V_1\Ga_{11}^2+V_2\Ga_{12}^2)
-\pa_wV_2(\pa_uV_1+V_1\Ga_{11}^1+V_2\Ga_{12}^1+1)}{(\pa_uV_1+V_1\Ga_{11}^1+V_2\Ga_{12}^1+1)
(\pa_vV_2+V_1\Ga_{12}^2+V_2\Ga_{22}^2+1)-(\pa_uV_2+V_1\Ga_{11}^2+V_2\Ga_{12}^2)
(\pa_vV_1+V_1\Ga_{12}^1+V_2\Ga_{22}^1)}-|(\pa_uV_1+V_1\Ga_{11}^1+V_2\Ga_{12}^1+1)\pa_ux_0
+(\pa_uV_2+V_1\Ga_{11}^2+V_2\Ga_{12}^2)\pa_vx_0|^2[V_1(\pa_vV_2+V_1\Ga_{12}^2+V_2\Ga_{22}^2+1)
-V_2(\pa_vV_1+V_1\Ga_{12}^1+V_2\Ga_{22}^1)]\\\frac{\pa_wV_1(\pa_vV_2+V_1\Ga_{12}^2+V_2\Ga_{22}^2+1)
-\pa_wV_2(\pa_vV_1+V_1\Ga_{12}^1+V_2\Ga_{22}^1)}{(\pa_uV_1+V_1\Ga_{11}^1+V_2\Ga_{12}^1+1)
(\pa_vV_2+V_1\Ga_{12}^2+V_2\Ga_{22}^2+1)-(\pa_uV_2+V_1\Ga_{11}^2+V_2\Ga_{12}^2)
(\pa_vV_1+V_1\Ga_{12}^1+V_2\Ga_{22}^1)}]\}=0,$

From $C=0$ for the eighth relation of R1 we get

$\pa_wV_2[\Ga_{12}^1(\pa_vV_1+V_1\Ga_{12}^1+V_2\Ga_{22}^1)
+\Ga_{22}^1(\pa_vV_2+V_1\Ga_{12}^2+V_2\Ga_{22}^2+1)]
-\pa_wV_1[\Ga_{12}^2(\pa_vV_1+V_1\Ga_{12}^1+V_2\Ga_{22}^1)
+\Ga_{22}^2(\pa_vV_2+V_1\Ga_{12}^2+V_2\Ga_{22}^2+1)]\\
+\frac{\pa_wV_1(\pa_vV_2+V_1\Ga_{12}^2+V_2\Ga_{22}^2+1)
-\pa_wV_2(\pa_vV_1+V_1\Ga_{12}^1+V_2\Ga_{22}^1)}{(\pa_uV_1+V_1\Ga_{11}^1+V_2\Ga_{12}^1+1)
(\pa_vV_2+V_1\Ga_{12}^2+V_2\Ga_{22}^2+1)-(\pa_uV_2+V_1\Ga_{11}^2+V_2\Ga_{12}^2)
(\pa_vV_1+V_1\Ga_{12}^1+V_2\Ga_{22}^1)}\{\pa_u(V_1\Ga_{12}^1\\+V_2\Ga_{22}^1)
(\pa_vV_2+V_1\Ga_{12}^2+V_2\Ga_{22}^2+1)-\pa_u(V_1\Ga_{12}^2+V_2\Ga_{22}^2)
(\pa_vV_1+V_1\Ga_{12}^1+V_2\Ga_{22}^1)+[\Ga_{11}^1(\pa_vV_1+V_1\Ga_{12}^1+V_2\Ga_{22}^1)
+\Ga_{12}^1(\pa_vV_2+V_1\Ga_{12}^2+V_2\Ga_{22}^2+1)](\pa_vV_2+V_1\Ga_{12}^2+V_2\Ga_{22}^2+1)
-[\Ga_{11}^2(\pa_vV_1+V_1\Ga_{12}^1+V_2\Ga_{22}^1)
+\Ga_{12}^2(\pa_vV_2+V_1\Ga_{12}^2+V_2\Ga_{22}^2+1)](\pa_vV_1+V_1\Ga_{12}^1+V_2\Ga_{22}^1)
-[\Ga_{12}^1(\pa_vV_1+V_1\Ga_{12}^1+V_2\Ga_{22}^1)
+\Ga_{22}^1(\pa_vV_2+V_1\Ga_{12}^2+V_2\Ga_{22}^2+1)](\pa_uV_2+V_1\Ga_{11}^2+V_2\Ga_{12}^2)
+[\Ga_{12}^2(\pa_vV_1+V_1\Ga_{12}^1+V_2\Ga_{22}^1)
+\Ga_{22}^2(\pa_vV_2+V_1\Ga_{12}^2+V_2\Ga_{22}^2+1)](\pa_uV_1+V_1\Ga_{11}^1+V_2\Ga_{12}^1+1)
-\pa_v(V_1\Ga_{11}^1+V_2\Ga_{12}^1)(\pa_vV_2+V_1\Ga_{12}^2+V_2\Ga_{22}^2+1)
+\pa_v(V_1\Ga_{11}^2+V_2\Ga_{12}^2)(\pa_vV_1+V_1\Ga_{12}^1+V_2\Ga_{22}^1)\}
-(\pa_wV_1\Ga_{12}^1+\pa_wV_2\Ga_{22}^1)(\pa_vV_2+V_1\Ga_{12}^2+V_2\Ga_{22}^2+1)
+(\pa_wV_1\Ga_{12}^2+\pa_wV_2\Ga_{22}^2)(\pa_vV_1+V_1\Ga_{12}^1+V_2\Ga_{22}^1)
-\frac{K}{(\pa_uV_1+V_1\Ga_{11}^1+V_2\Ga_{12}^1+1)
(\pa_vV_2+V_1\Ga_{12}^2+V_2\Ga_{22}^2+1)-(\pa_uV_2+V_1\Ga_{11}^2+V_2\Ga_{12}^2)
(\pa_vV_1+V_1\Ga_{12}^1+V_2\Ga_{22}^1)}[\pa_wV_1(\pa_vV_2+V_1\Ga_{12}^2+V_2\Ga_{22}^2+1)
-\pa_wV_2(\pa_vV_1+V_1\Ga_{12}^1+V_2\Ga_{22}^1)]
[V^T\pa_ux_0(\pa_vV_1+V_1\Ga_{12}^1+V_2\Ga_{22}^1)
+V^T\pa_vx_0(\pa_vV_2+V_1\Ga_{12}^2+V_2\Ga_{22}^2+1)]=0,$

From $C=0$ for the ninth relation of R1 we get

$-\pa_wV_2[\Ga_{11}^1(\pa_uV_1+V_1\Ga_{11}^1+V_2\Ga_{12}^1+1)
+\Ga_{12}^1(\pa_uV_2+V_1\Ga_{11}^2+V_2\Ga_{12}^2)]
+\pa_wV_1[\Ga_{11}^2(\pa_uV_1+V_1\Ga_{11}^1+V_2\Ga_{12}^1+1)
+\Ga_{12}^2(\pa_uV_2+V_1\Ga_{11}^2+V_2\Ga_{12}^2)]\\
+\frac{\pa_wV_1(\pa_uV_2+V_1\Ga_{11}^2+V_2\Ga_{12}^2)
-\pa_wV_2(\pa_uV_1+V_1\Ga_{11}^1+V_2\Ga_{12}^1+1)}{(\pa_uV_1+V_1\Ga_{11}^1+V_2\Ga_{12}^1+1)
(\pa_vV_2+V_1\Ga_{12}^2+V_2\Ga_{22}^2+1)-(\pa_uV_2+V_1\Ga_{11}^2+V_2\Ga_{12}^2)
(\pa_vV_1+V_1\Ga_{12}^1+V_2\Ga_{22}^1)}\{-[\Ga_{11}^1(\pa_uV_1\\+V_1\Ga_{11}^1+V_2\Ga_{12}^1+1)
+\Ga_{12}^1(\pa_uV_2+V_1\Ga_{11}^2+V_2\Ga_{12}^2)](\pa_vV_2+V_1\Ga_{12}^2+V_2\Ga_{22}^2+1)
+[\Ga_{11}^2(\pa_uV_1+V_1\Ga_{11}^1+V_2\Ga_{12}^1+1)
+\Ga_{12}^2(\pa_uV_2+V_1\Ga_{11}^2+V_2\Ga_{12}^2)](\pa_vV_1+V_1\Ga_{12}^1+V_2\Ga_{22}^1)
+\pa_u(V_1\Ga_{12}^2+V_2\Ga_{22}^2)(\pa_uV_1+V_1\Ga_{11}^1+V_2\Ga_{12}^1+1)
-\pa_u(V_1\Ga_{12}^1+V_2\Ga_{22}^1)(\pa_uV_2+V_1\Ga_{11}^2+V_2\Ga_{12}^2)
+\pa_v(V_1\Ga_{11}^1+V_2\Ga_{12}^1)(\pa_uV_2+V_1\Ga_{11}^2+V_2\Ga_{12}^2)
-\pa_v(V_1\Ga_{11}^2+V_2\Ga_{12}^2)(\pa_uV_1+V_1\Ga_{11}^1+V_2\Ga_{12}^1+1)
+[\Ga_{12}^1(\pa_uV_1+V_1\Ga_{11}^1+V_2\Ga_{12}^1+1)
+\Ga_{22}^1(\pa_uV_2+V_1\Ga_{11}^2+V_2\Ga_{12}^2)](\pa_uV_2+V_1\Ga_{11}^2+V_2\Ga_{12}^2)
-[\Ga_{12}^2(\pa_uV_1+V_1\Ga_{11}^1+V_2\Ga_{12}^1+1)
+\Ga_{22}^2(\pa_uV_2+V_1\Ga_{11}^2+V_2\Ga_{12}^2)](\pa_uV_1+V_1\Ga_{11}^1+V_2\Ga_{12}^1+1)\}
+(\pa_wV_1\Ga_{11}^1+\pa_wV_2\Ga_{12}^1)(\pa_uV_2+V_1\Ga_{11}^2+V_2\Ga_{12}^2)
-(\pa_wV_1\Ga_{11}^2+\pa_wV_2\Ga_{12}^2)(\pa_uV_1+V_1\Ga_{11}^1+V_2\Ga_{12}^1+1)
+\frac{K}{(\pa_uV_1+V_1\Ga_{11}^1+V_2\Ga_{12}^1+1)
(\pa_vV_2+V_1\Ga_{12}^2+V_2\Ga_{22}^2+1)-(\pa_uV_2+V_1\Ga_{11}^2+V_2\Ga_{12}^2)
(\pa_vV_1+V_1\Ga_{12}^1+V_2\Ga_{22}^1)}\\\{\pa_wV_1(\pa_uV_2+V_1\Ga_{11}^2+V_2\Ga_{12}^2)
-\pa_wV_2(\pa_uV_1+V_1\Ga_{11}^1+V_2\Ga_{12}^1+1)\}
[V^T\pa_ux_0(\pa_uV_1+V_1\Ga_{11}^1+V_2\Ga_{12}^1+1)
+V^T\pa_vx_0(\pa_uV_2+V_1\Ga_{11}^2+V_2\Ga_{12}^2)]=0,$

From $C=0$ for the first relation of R2 we get

$-\pa_wV_1\pa_u(V_1\Ga_{11}^2+V_2\Ga_{12}^2)
+\pa_wV_2\pa_u(V_1\Ga_{11}^1+V_2\Ga_{12}^1)\\
+\frac{\pa_wV_1(\pa_uV_2+V_1\Ga_{11}^2+V_2\Ga_{12}^2)
-\pa_wV_2(\pa_uV_1+V_1\Ga_{11}^1+V_2\Ga_{12}^1+1)}{(\pa_uV_1+V_1\Ga_{11}^1+V_2\Ga_{12}^1+1)
(\pa_vV_2+V_1\Ga_{12}^2+V_2\Ga_{22}^2+1)-(\pa_uV_2+V_1\Ga_{11}^2+V_2\Ga_{12}^2)
(\pa_vV_1+V_1\Ga_{12}^1+V_2\Ga_{22}^1)}\{\pa_u(V_1\Ga_{11}^1\\+V_2\Ga_{12}^1)
(\pa_vV_2+V_1\Ga_{12}^2+V_2\Ga_{22}^2+1)-\pa_u(V_1\Ga_{11}^2+V_2\Ga_{12}^2)
(\pa_vV_1+V_1\Ga_{12}^1+V_2\Ga_{22}^1)+\pa_u(V_1\Ga_{12}^2+V_2\Ga_{22}^2)
(\pa_uV_1+V_1\Ga_{11}^1+V_2\Ga_{12}^1+1)-\pa_u(V_1\Ga_{12}^1+V_2\Ga_{22}^1)
(\pa_uV_2+V_1\Ga_{11}^2+V_2\Ga_{12}^2)+\pa_v(V_1\Ga_{11}^1+V_2\Ga_{12}^1)
(\pa_uV_2+V_1\Ga_{11}^2+V_2\Ga_{12}^2)-\pa_v(V_1\Ga_{11}^2+V_2\Ga_{12}^2)
(\pa_uV_1+V_1\Ga_{11}^1+V_2\Ga_{12}^1+1)+
[\Ga_{12}^1(\pa_uV_1+V_1\Ga_{11}^1+V_2\Ga_{12}^1+1)
+\Ga_{22}^1(\pa_uV_2+V_1\Ga_{11}^2+V_2\Ga_{12}^2)](\pa_uV_2+V_1\Ga_{11}^2+V_2\Ga_{12}^2)
-[\Ga_{12}^2(\pa_uV_1+V_1\Ga_{11}^1+V_2\Ga_{12}^1+1)
+\Ga_{22}^2(\pa_uV_2+V_1\Ga_{11}^2+V_2\Ga_{12}^2)](\pa_uV_1+V_1\Ga_{11}^1+V_2\Ga_{12}^1+1)\}\\
+\frac{\pa_wV_1(\pa_vV_2+V_1\Ga_{12}^2+V_2\Ga_{22}^2+1)
-\pa_wV_2(\pa_vV_1+V_1\Ga_{12}^1+V_2\Ga_{22}^1)}{(\pa_uV_1+V_1\Ga_{11}^1+V_2\Ga_{12}^1+1)
(\pa_vV_2+V_1\Ga_{12}^2+V_2\Ga_{22}^2+1)-(\pa_uV_2+V_1\Ga_{11}^2+V_2\Ga_{12}^2)
(\pa_vV_1+V_1\Ga_{12}^1+V_2\Ga_{22}^1)}\{-\pa_u(V_1\Ga_{11}^1\\+V_2\Ga_{12}^1)
(\pa_uV_2+V_1\Ga_{11}^2+V_2\Ga_{12}^2)+\pa_u(V_1\Ga_{11}^2+V_2\Ga_{12}^2)
(\pa_uV_1+V_1\Ga_{11}^1+V_2\Ga_{12}^1+1)-[\Ga_{11}^1(\pa_uV_1+V_1\Ga_{11}^1+V_2\Ga_{12}^1+1)
+\Ga_{12}^1(\pa_uV_2+V_1\Ga_{11}^2+V_2\Ga_{12}^2)](\pa_uV_2+V_1\Ga_{11}^2+V_2\Ga_{12}^2)
+[\Ga_{11}^2(\pa_uV_1+V_1\Ga_{11}^1+V_2\Ga_{12}^1+1)
+\Ga_{12}^2(\pa_uV_2+V_1\Ga_{11}^2+V_2\Ga_{12}^2)](\pa_uV_1+V_1\Ga_{11}^1+V_2\Ga_{12}^1+1)\}
+(\pa_wV_1\Ga_{11}^1+\pa_wV_2\Ga_{12}^1)(\pa_uV_2+V_1\Ga_{11}^2+V_2\Ga_{12}^2)
-(\pa_wV_1\Ga_{11}^2+\pa_wV_2\Ga_{12}^2)(\pa_uV_1+V_1\Ga_{11}^1+V_2\Ga_{12}^1+1)\\
+\frac{K}{(\pa_uV_1+V_1\Ga_{11}^1+V_2\Ga_{12}^1+1)
(\pa_vV_2+V_1\Ga_{12}^2+V_2\Ga_{22}^2+1)-(\pa_uV_2+V_1\Ga_{11}^2+V_2\Ga_{12}^2)
(\pa_vV_1+V_1\Ga_{12}^1+V_2\Ga_{22}^1)}[\pa_wV_1(\pa_uV_2+V_1\Ga_{11}^2+V_2\Ga_{12}^2)
-\pa_wV_2(\pa_uV_1+V_1\Ga_{11}^1+V_2\Ga_{12}^1+1)]
[V^T\pa_ux_0(\pa_uV_1+V_1\Ga_{11}^1+V_2\Ga_{12}^1+1)
+V^T\pa_vx_0(\pa_uV_2+V_1\Ga_{11}^2+V_2\Ga_{12}^2)]=0$. This is
the same as the relation $C=0$ for the ninth relation of R1.

From $\hat C=0$ for the first relation of R2 we get

$\pa_wV_1\pa_v(V_1\Ga_{12}^2+V_2\Ga_{22}^2)-\pa_wV_2\pa_v(V_1\Ga_{12}^1+V_2\Ga_{22}^1)
+\\\frac{\pa_wV_1(\pa_vV_2+V_1\Ga_{12}^2+V_2\Ga_{22}^2+1)
-\pa_wV_2(\pa_vV_1+V_1\Ga_{12}^1+V_2\Ga_{22}^1)}{(\pa_uV_1+V_1\Ga_{11}^1+V_2\Ga_{12}^1+1)
(\pa_vV_2+V_1\Ga_{12}^2+V_2\Ga_{22}^2+1)-(\pa_uV_2+V_1\Ga_{11}^2+V_2\Ga_{12}^2)
(\pa_vV_1+V_1\Ga_{12}^1+V_2\Ga_{22}^1)}\{-\pa_v(V_1\Ga_{11}^1\\+V_2\Ga_{12}^1)
(\pa_vV_2+V_1\Ga_{12}^2+V_2\Ga_{22}^2+1)+\pa_v(V_1\Ga_{11}^2+V_2\Ga_{12}^2)
(\pa_vV_1+V_1\Ga_{12}^1+V_2\Ga_{22}^1)-\pa_v(V_1\Ga_{12}^2+V_2\Ga_{22}^2)
(\pa_uV_1+V_1\Ga_{11}^1+V_2\Ga_{12}^1+1)+\pa_v(V_1\Ga_{12}^1+V_2\Ga_{22}^1)
(\pa_uV_2+V_1\Ga_{11}^2+V_2\Ga_{12}^2)+\pa_u(V_1\Ga_{12}^1+V_2\Ga_{22}^1)
(\pa_vV_2+V_1\Ga_{12}^2+V_2\Ga_{22}^2+1)-\pa_u(V_1\Ga_{12}^2+V_2\Ga_{22}^2)
(\pa_vV_1+V_1\Ga_{12}^1+V_2\Ga_{22}^1)+[\Ga_{11}^1(\pa_vV_1+V_1\Ga_{12}^1+V_2\Ga_{22}^1)
+\Ga_{12}^1(\pa_vV_2+V_1\Ga_{12}^2+V_2\Ga_{22}^2+1)](\pa_vV_2+V_1\Ga_{12}^2+V_2\Ga_{22}^2+1)
-[\Ga_{11}^2(\pa_vV_1+V_1\Ga_{12}^1+V_2\Ga_{22}^1)
+\Ga_{12}^2(\pa_vV_2+V_1\Ga_{12}^2+V_2\Ga_{22}^2+1)](\pa_vV_1+V_1\Ga_{12}^1+V_2\Ga_{22}^1)\}\\
+\frac{\pa_wV_1(\pa_uV_2+V_1\Ga_{11}^2+V_2\Ga_{12}^2)
-\pa_wV_2(\pa_uV_1+V_1\Ga_{11}^1+V_2\Ga_{12}^1+1)}{(\pa_uV_1+V_1\Ga_{11}^1+V_2\Ga_{12}^1+1)
(\pa_vV_2+V_1\Ga_{12}^2+V_2\Ga_{22}^2+1)-(\pa_uV_2+V_1\Ga_{11}^2+V_2\Ga_{12}^2)
(\pa_vV_1+V_1\Ga_{12}^1+V_2\Ga_{22}^1)}\{-\pa_v(V_1\Ga_{12}^1\\+V_2\Ga_{22}^1)
(\pa_vV_2+V_1\Ga_{12}^2+V_2\Ga_{22}^2+1)+\pa_v(V_1\Ga_{12}^2+V_2\Ga_{22}^2)
(\pa_vV_1+V_1\Ga_{12}^1+V_2\Ga_{22}^1)-[\Ga_{12}^1(\pa_vV_1+V_1\Ga_{12}^1+V_2\Ga_{22}^1)
+\Ga_{22}^1(\pa_vV_2+V_1\Ga_{12}^2+V_2\Ga_{22}^2+1)](\pa_vV_2+V_1\Ga_{12}^2+V_2\Ga_{22}^2+1)
+[\Ga_{12}^2(\pa_vV_1+V_1\Ga_{12}^1+V_2\Ga_{22}^1)
+\Ga_{22}^2(\pa_vV_2+V_1\Ga_{12}^2+V_2\Ga_{22}^2+1)](\pa_vV_1+V_1\Ga_{12}^1+V_2\Ga_{22}^1)\}
-(\pa_wV_1\Ga_{12}^1+\pa_wV_2\Ga_{22}^1)(\pa_vV_2+V_1\Ga_{12}^2+V_2\Ga_{22}^2+1)
+(\pa_wV_1\Ga_{12}^2+\pa_wV_2\Ga_{22}^2)(\pa_vV_1+V_1\Ga_{12}^1+V_2\Ga_{22}^1)\\
-\frac{K}{(\pa_uV_1+V_1\Ga_{11}^1+V_2\Ga_{12}^1+1)
(\pa_vV_2+V_1\Ga_{12}^2+V_2\Ga_{22}^2+1)-(\pa_uV_2+V_1\Ga_{11}^2+V_2\Ga_{12}^2)
(\pa_vV_1+V_1\Ga_{12}^1+V_2\Ga_{22}^1)}[\pa_wV_1(\pa_vV_2+V_1\Ga_{12}^2+V_2\Ga_{22}^2+1)
-\pa_wV_2(\pa_vV_1+V_1\Ga_{12}^1+V_2\Ga_{22}^1)]
[V^T\pa_ux_0(\pa_vV_1+V_1\Ga_{12}^1+V_2\Ga_{22}^1)
+V^T\pa_vx_0(\pa_vV_2+V_1\Ga_{12}^2+V_2\Ga_{22}^2+1)]=0.$ This is
the same as the relation $C=0$ for the eighth relation od R1.

From $C=0$ for the second relation of R2 we get

$(\pa_wV_1\Ga_{12}^1+\pa_wV_2\Ga_{22}^1)(\pa_vV_2+V_1\Ga_{12}^2+V_2\Ga_{22}^2+1)
-(\pa_wV_1\Ga_{12}^2+\pa_wV_2\Ga_{22}^2)(\pa_vV_1+V_1\Ga_{12}^1+V_2\Ga_{22}^1)
+\frac{\pa_wV_1(\pa_vV_2+V_1\Ga_{12}^2+V_2\Ga_{22}^2+1)
-\pa_wV_2(\pa_vV_1+V_1\Ga_{12}^1+V_2\Ga_{22}^1)}{(\pa_uV_1+V_1\Ga_{11}^1+V_2\Ga_{12}^1+1)
(\pa_vV_2+V_1\Ga_{12}^2+V_2\Ga_{22}^2+1)-(\pa_uV_2+V_1\Ga_{11}^2+V_2\Ga_{12}^2)
(\pa_vV_1+V_1\Ga_{12}^1+V_2\Ga_{22}^1)}[\pa_v(V_1\Ga_{11}^1+V_2\Ga_{12}^1)
(\pa_vV_2+V_1\Ga_{12}^2+V_2\Ga_{22}^2+1)-\pa_v(V_1\Ga_{11}^2+V_2\Ga_{12}^2)
(\pa_vV_1+V_1\Ga_{12}^1+V_2\Ga_{22}^1)]\\
+\frac{\pa_wV_1(\pa_uV_2+V_1\Ga_{11}^2+V_2\Ga_{12}^2)
-\pa_wV_2(\pa_uV_1+V_1\Ga_{11}^1+V_2\Ga_{12}^1+1)}{(\pa_uV_1+V_1\Ga_{11}^1+V_2\Ga_{12}^1+1)
(\pa_vV_2+V_1\Ga_{12}^2+V_2\Ga_{22}^2+1)-(\pa_uV_2+V_1\Ga_{11}^2+V_2\Ga_{12}^2)
(\pa_vV_1+V_1\Ga_{12}^1+V_2\Ga_{22}^1)}\{[\Ga_{12}^1(\pa_vV_1+V_1\Ga_{12}^1+V_2\Ga_{22}^1)
+\Ga_{22}^1(\pa_vV_2+V_1\Ga_{12}^2+V_2\Ga_{22}^2+1)](\pa_vV_2+V_1\Ga_{12}^2+V_2\Ga_{22}^2+1)
-[\Ga_{12}^2(\pa_vV_1+V_1\Ga_{12}^1+V_2\Ga_{22}^1)
+\Ga_{22}^2(\pa_vV_2+V_1\Ga_{12}^2+V_2\Ga_{22}^2+1)](\pa_vV_1+V_1\Ga_{12}^1+V_2\Ga_{22}^1)\}\\
-\frac{\pa_wV_1(\pa_vV_2+V_1\Ga_{12}^2+V_2\Ga_{22}^2+1)
-\pa_wV_2(\pa_vV_1+V_1\Ga_{12}^1+V_2\Ga_{22}^1)}{(\pa_uV_1+V_1\Ga_{11}^1+V_2\Ga_{12}^1+1)
(\pa_vV_2+V_1\Ga_{12}^2+V_2\Ga_{22}^2+1)-(\pa_uV_2+V_1\Ga_{11}^2+V_2\Ga_{12}^2)
(\pa_vV_1+V_1\Ga_{12}^1+V_2\Ga_{22}^1)}\{[\pa_u(V_1\Ga_{12}^1+V_2\Ga_{22}^1)
+\Ga_{11}^1(\pa_vV_1+V_1\Ga_{12}^1+V_2\Ga_{22}^1)
+\Ga_{12}^1(\pa_vV_2+V_1\Ga_{12}^1+V_2\Ga_{22}^2+1)](\pa_vV_2+V_1\Ga_{12}^2+V_2\Ga_{22}^2+1)
-[\pa_u(V_1\Ga_{12}^2+V_2\Ga_{22}^2)+\Ga_{11}^2(\pa_vV_1+V_1\Ga_{12}^1+V_2\Ga_{22}^1)
+\Ga_{12}^2(\pa_vV_2+V_1\Ga_{12}^2+V_2\Ga_{22}^2+1)](\pa_vV_1+V_1\Ga_{12}^1+V_2\Ga_{22}^1)\}
+\frac{\pa_wV_1(\pa_vV_2+V_1\Ga_{12}^2+V_2\Ga_{22}^2+1)
-\pa_wV_2(\pa_vV_1+V_1\Ga_{12}^1+V_2\Ga_{22}^1)}{(\pa_uV_1+V_1\Ga_{11}^1+V_2\Ga_{12}^1+1)
(\pa_vV_2+V_1\Ga_{12}^2+V_2\Ga_{22}^2+1)-(\pa_uV_2+V_1\Ga_{11}^2+V_2\Ga_{12}^2)
(\pa_vV_1+V_1\Ga_{12}^1+V_2\Ga_{22}^1)}K[(\pa_vV_1+V_1\Ga_{12}^1+V_2\Ga_{22}^1)V_1g_{11}
+[(\pa_vV_1+V_1\Ga_{12}^1+V_2\Ga_{22}^1)V_2+(\pa_vV_2+V_1\Ga_{12}^2+V_2\Ga_{22}^2+1)V_1]g_{12}
+(\pa_vV_2+V_1\Ga_{12}^2+V_2\Ga_{22}^2+1)V_2g_{22}]=0.$ This is
minus the relation $C=0$ for the eighth relation of R1.

From $\hat C=0$ for the second relation of R2 we get

$(\pa_wV_1\Ga_{11}^1+\pa_wV_2\Ga_{12}^1)(\pa_uV_2+V_1\Ga_{11}^2+V_2\Ga_{12}^2)
-(\pa_wV_1\Ga_{11}^2+\pa_wV_2\Ga_{12}^2)(\pa_uV_1+V_1\Ga_{11}^1+V_2\Ga_{12}^1+1)
+\frac{\pa_wV_1(\pa_uV_2+V_1\Ga_{11}^2+V_2\Ga_{12}^2)
-\pa_wV_2(\pa_uV_1+V_1\Ga_{11}^1+V_2\Ga_{12}^1+1)}{(\pa_uV_1+V_1\Ga_{11}^1+V_2\Ga_{12}^1+1)
(\pa_vV_2+V_1\Ga_{12}^2+V_2\Ga_{22}^2+1)-(\pa_uV_2+V_1\Ga_{11}^2+V_2\Ga_{12}^2)
(\pa_vV_1+V_1\Ga_{12}^1+V_2\Ga_{22}^1)}[\pa_u(V_1\Ga_{12}^2+V_2\Ga_{22}^2)
(\pa_uV_1+V_1\Ga_{11}^1+V_2\Ga_{12}^1+1)-\pa_u(V_1\Ga_{12}^1+V_2\Ga_{22}^1)
(\pa_uV_2+V_1\Ga_{11}^2+V_2\Ga_{12}^2)]\\-
\frac{\pa_wV_1(\pa_vV_2+V_1\Ga_{12}^2+V_2\Ga_{22}^2+1)
-\pa_wV_2(\pa_vV_1+V_1\Ga_{12}^1+V_2\Ga_{22}^1)}{(\pa_uV_1+V_1\Ga_{11}^1+V_2\Ga_{12}^1+1)
(\pa_vV_2+V_1\Ga_{12}^2+V_2\Ga_{22}^2+1)-(\pa_uV_2+V_1\Ga_{11}^2+V_2\Ga_{12}^2)
(\pa_vV_1+V_1\Ga_{12}^1+V_2\Ga_{22}^1)}\{[\Ga_{11}^1(\pa_uV_1+V_1\Ga_{11}^1+V_2\Ga_{12}^1+1)
+\Ga_{12}^1(\pa_uV_2+V_1\Ga_{11}^2+V_2\Ga_{12}^2)](\pa_uV_2+V_1\Ga_{11}^2+V_2\Ga_{12}^2)
-[\Ga_{11}^2(\pa_uV_1+V_1\Ga_{11}^1+V_2\Ga_{12}^1+1)
+\Ga_{12}^2(\pa_uV_2+V_1\Ga_{11}^2+V_2\Ga_{12}^2)](\pa_uV_1+V_1\Ga_{11}^1+V_2\Ga_{12}^1+1)\}
\\+\frac{\pa_wV_1(\pa_uV_2+V_1\Ga_{11}^2+V_2\Ga_{12}^2)
-\pa_wV_2(\pa_uV_1+V_1\Ga_{11}^1+V_2\Ga_{12}^1+1)}{(\pa_uV_1+V_1\Ga_{11}^1+V_2\Ga_{12}^1+1)
(\pa_vV_2+V_1\Ga_{12}^2+V_2\Ga_{22}^2+1)-(\pa_uV_2+V_1\Ga_{11}^2+V_2\Ga_{12}^2)
(\pa_vV_1+V_1\Ga_{12}^1+V_2\Ga_{22}^1)}\{[\pa_v(V_1\Ga_{11}^1+V_2\Ga_{12}^1)
+\Ga_{12}^1(\pa_uV_1+V_1\Ga_{11}^1+V_2\Ga_{12}^1+1)+
\Ga_{22}^1(\pa_uV_2+V_1\Ga_{11}^2+V_2\Ga_{12}^2)](\pa_uV_2+V_1\Ga_{11}^2+V_2\Ga_{12}^2)-
[\pa_v(V_1\Ga_{11}^2+V_2\Ga_{12}^2)+\Ga_{12}^2(\pa_uV_1+V_1\Ga_{11}^1+V_2\Ga_{12}^1+1)+
\Ga_{22}^2(\pa_uV_2+V_1\Ga_{11}^2+V_2\Ga_{12}^2)](\pa_uV_1+V_1\Ga_{11}^1+V_2\Ga_{12}^1+1)\}+
\frac{\pa_wV_1(\pa_uV_2+V_1\Ga_{11}^2+V_2\Ga_{12}^2)
-\pa_wV_2(\pa_uV_1+V_1\Ga_{11}^1+V_2\Ga_{12}^1+1)}{(\pa_uV_1+V_1\Ga_{11}^1+V_2\Ga_{12}^1+1)
(\pa_vV_2+V_1\Ga_{12}^2+V_2\Ga_{22}^2+1)-(\pa_uV_2+V_1\Ga_{11}^2+V_2\Ga_{12}^2)
(\pa_vV_1+V_1\Ga_{12}^1+V_2\Ga_{22}^1)}K[V_1(\pa_uV_1+V_1\Ga_{11}^1+V_2\Ga_{12}^1+1)g_{11}+
[V_1(\pa_uV_2+V_1\Ga_{11}^2+V_2\Ga_{12}^2)+V_2(\pa_uV_1+V_1\Ga_{11}^1+V_2\Ga_{12}^1+1)]g_{12}+
V_2(\pa_uV_2+V_1\Ga_{11}^2+V_2\Ga_{12}^2)g_{22}]=0.$ This is the
same as the relation $C=0$ for the ninth relation of R1.

From $C=0$ for the third relation of R2 we get

$(\pa_wV_1\Ga_{11}^1+\pa_wV_2\Ga_{12}^1)(\pa_vV_2+V_1\Ga_{12}^2+V_2\Ga_{22}^2+1)
-(\pa_wV_1\Ga_{11}^2+\pa_wV_2\Ga_{12}^2)(\pa_vV_1+V_1\Ga_{12}^1+V_2\Ga_{22}^1)
+(\pa_wV_1\Ga_{12}^1+\pa_wV_2\Ga_{22}^1)(\pa_uV_2+V_1\Ga_{11}^2+V_2\Ga_{12}^2)
-(\pa_wV_1\Ga_{12}^2+\pa_wV_2\Ga_{22}^2)(\pa_uV_1+V_1\Ga_{11}^1+V_2\Ga_{12}^1+1)
-\pa_wV_1\pa_v(V_1\Ga_{11}^2+V_2\Ga_{12}^2)+\pa_wV_2\pa_v(V_1\Ga_{11}^1+V_2\Ga_{12}^1)
-\pa_wV_1\pa_u(V_1\Ga_{12}^2+V_2\Ga_{22}^2)+\pa_wV_2\pa_u(V_1\Ga_{12}^1+V_2\Ga_{22}^1)
+\frac{\pa_wV_1(\pa_uV_2+V_1\Ga_{11}^2+V_2\Ga_{12}^2)
-\pa_wV_2(\pa_uV_1+V_1\Ga_{11}^1+V_2\Ga_{12}^1+1)}{(\pa_uV_1+V_1\Ga_{11}^1+V_2\Ga_{12}^1+1)
(\pa_vV_2+V_1\Ga_{12}^2+V_2\Ga_{22}^2+1)-(\pa_uV_2+V_1\Ga_{11}^2+V_2\Ga_{12}^2)
(\pa_vV_1+V_1\Ga_{12}^1+V_2\Ga_{22}^1)}\{2\pa_v(V_1\Ga_{11}^1\\+V_2\Ga_{12}^1)
(\pa_vV_2+V_1\Ga_{12}^2+V_2\Ga_{22}^2+1)-2\pa_v(V_1\Ga_{11}^2+V_2\Ga_{12}^2)
(\pa_vV_1+V_1\Ga_{12}^1+V_2\Ga_{22}^1)+[\Ga_{12}^1(\pa_uV_1+V_1\Ga_{11}^1+V_2\Ga_{12}^1+1)
+\Ga_{22}^1(\pa_uV_2+V_1\Ga_{11}^2+V_2\Ga_{12}^2)](\pa_vV_2+V_1\Ga_{12}^2+V_2\Ga_{22}^2+1)
-[\Ga_{12}^2(\pa_uV_1+V_1\Ga_{11}^1+V_2\Ga_{12}^1+1)
+\Ga_{22}^2(\pa_uV_2+V_1\Ga_{11}^2+V_2\Ga_{12}^2)](\pa_vV_1+V_1\Ga_{12}^1+V_2\Ga_{22}^1)
+[\Ga_{12}^1(\pa_vV_1+V_1\Ga_{12}^1+V_2\Ga_{22}^1)
+\Ga_{22}^1(\pa_vV_2+V_1\Ga_{12}^2+V_2\Ga_{22}^2+1)](\pa_uV_2+V_1\Ga_{11}^2+V_2\Ga_{12}^2)
-[\Ga_{12}^2(\pa_vV_1+V_1\Ga_{12}^1+V_2\Ga_{22}^1)
+\Ga_{22}^2(\pa_vV_2+V_1\Ga_{12}^2+V_2\Ga_{22}^2+1)](\pa_uV_1+V_1\Ga_{11}^1+V_2\Ga_{12}^1+1)\}
+\frac{\pa_wV_1(\pa_vV_2+V_1\Ga_{12}^2+V_2\Ga_{22}^2+1)
-\pa_wV_2(\pa_vV_1+V_1\Ga_{12}^1+V_2\Ga_{22}^1)}{(\pa_uV_1+V_1\Ga_{11}^1+V_2\Ga_{12}^1+1)
(\pa_vV_2+V_1\Ga_{12}^2+V_2\Ga_{22}^2+1)-(\pa_uV_2+V_1\Ga_{11}^2+V_2\Ga_{12}^2)
(\pa_vV_1+V_1\Ga_{12}^1+V_2\Ga_{22}^1)}\\\{-2\pa_u(V_1\Ga_{12}^1+V_2\Ga_{22}^1)
(\pa_uV_2+V_1\Ga_{11}^2+V_2\Ga_{12}^2)+2\pa_u(V_1\Ga_{12}^2+V_2\Ga_{22}^2)
(\pa_uV_1+V_1\Ga_{11}^1+V_2\Ga_{12}^1+1)-[\Ga_{11}^1(\pa_uV_1+V_1\Ga_{11}^1+V_2\Ga_{12}^1+1)
+\Ga_{12}^1(\pa_uV_2+V_1\Ga_{11}^2+V_2\Ga_{12}^2)](\pa_vV_2+V_1\Ga_{12}^2+V_2\Ga_{22}^2+1)
+[\Ga_{11}^2(\pa_uV_1+V_1\Ga_{11}^1+V_2\Ga_{12}^1+1)
+\Ga_{12}^2(\pa_uV_2+V_1\Ga_{11}^2+V_2\Ga_{12}^2)](\pa_vV_1+V_1\Ga_{12}^1+V_2\Ga_{22}^1)\}
-[\Ga_{11}^1(\pa_vV_1+V_1\Ga_{12}^1+V_2\Ga_{22}^1)
+\Ga_{12}^1(\pa_vV_2+V_1\Ga_{12}^2+V_2\Ga_{22}^2+1)](\pa_uV_2+V_1\Ga_{11}^2+V_2\Ga_{12}^2)
+[\Ga_{11}^2(\pa_vV_1+V_1\Ga_{12}^1+V_2\Ga_{22}^1)
+\Ga_{12}^2(\pa_vV_2+V_1\Ga_{12}^2+V_2\Ga_{22}^2+1)](\pa_uV_1+V_1\Ga_{11}^1+V_2\Ga_{12}^1+1)\}\\
+\frac{K}{(\pa_uV_1+V_1\Ga_{11}^1+V_2\Ga_{12}^1+1)
(\pa_vV_2+V_1\Ga_{12}^2+V_2\Ga_{22}^2+1)-(\pa_uV_2+V_1\Ga_{11}^2+V_2\Ga_{12}^2)
(\pa_vV_1+V_1\Ga_{12}^1+V_2\Ga_{22}^1)}\{(\pa_wV_1V_2\\-\pa_wV_2V_1)
[(\pa_uV_1+V_1\Ga_{11}^1+V_2\Ga_{12}^1+1)\pa_ux_0
+(\pa_uV_2+V_1\Ga_{11}^2+V_2\Ga_{12}^2)\pa_vx_0]^T[(\pa_vV_1+V_1\Ga_{12}^1+V_2\Ga_{22}^1)\pa_ux_0
+(\pa_vV_2+V_1\Ga_{12}^2+V_2\Ga_{22}^2+1)\pa_vx_0]\\
-\frac{\pa_wV_1(\pa_uV_2+V_1\Ga_{11}^2+V_2\Ga_{12}^2)
-\pa_wV_2(\pa_uV_1+V_1\Ga_{11}^1+V_2\Ga_{12}^1+1)}{(\pa_uV_1+V_1\Ga_{11}^1+V_2\Ga_{12}^1+1)
(\pa_vV_2+V_1\Ga_{12}^2+V_2\Ga_{22}^2+1)-(\pa_uV_2+V_1\Ga_{11}^2+V_2\Ga_{12}^2)
(\pa_vV_1+V_1\Ga_{12}^1+V_2\Ga_{22}^1)}|(\pa_vV_1+V_1\Ga_{12}^1+V_2\Ga_{22}^1)\pa_ux_0
+(\pa_vV_2+V_1\Ga_{12}^2+V_2\Ga_{22}^2+1)\pa_vx_0|^2[V_1(\pa_uV_2+V_1\Ga_{11}^2+V_2\Ga_{12}^2)
-V_2(\pa_uV_1+V_1\Ga_{11}^1+V_2\Ga_{12}^1+1)]
+\frac{\pa_wV_1(\pa_vV_2+V_1\Ga_{12}^2+V_2\Ga_{22}^2+1)
-\pa_wV_2(\pa_vV_1+V_1\Ga_{12}^1+V_2\Ga_{22}^1)}{(\pa_uV_1+V_1\Ga_{11}^1+V_2\Ga_{12}^1+1)
(\pa_vV_2+V_1\Ga_{12}^2+V_2\Ga_{22}^2+1)-(\pa_uV_2+V_1\Ga_{11}^2+V_2\Ga_{12}^2)
(\pa_vV_1+V_1\Ga_{12}^1+V_2\Ga_{22}^1)}|(\pa_uV_1+V_1\Ga_{11}^1+V_2\Ga_{12}^1+1)\pa_ux_0
+(\pa_uV_2+V_1\Ga_{11}^2+V_2\Ga_{12}^2)\pa_vx_0|^2[V_1(\pa_vV_2+V_1\Ga_{12}^2+V_2\Ga_{22}^2+1)
-V_2(\pa_vV_1+V_1\Ga_{12}^1+V_2\Ga_{22}^1)]\}=0,$

From $\hat C=0$ for the third relation of R2 we get

$(\pa_wV_1\Ga_{12}^1+\pa_wV_2\Ga_{22}^1)(\pa_uV_2+V_1\Ga_{11}^2+V_2\Ga_{12}^2)
-(\pa_wV_1\Ga_{12}^2+\pa_wV_2\Ga_{22}^2)(\pa_uV_1+V_1\Ga_{11}^1+V_2\Ga_{12}^1+1)
+(\pa_wV_1\Ga_{11}^1+\pa_wV_2\Ga_{12}^1)(\pa_vV_2+V_1\Ga_{12}^2+V_2\Ga_{22}^2+1)
-(\pa_wV_1\Ga_{11}^2+\pa_wV_2\Ga_{12}^2)(\pa_vV_1+V_1\Ga_{12}^1+V_2\Ga_{22}^1)
-\pa_wV_1\pa_u(V_1\Ga_{12}^2+V_2\Ga_{22}^2)+\pa_wV_2\pa_u(V_1\Ga_{12}^1+V_2\Ga_{22}^1)
-\pa_wV_1\pa_v(V_1\Ga_{11}^2+V_2\Ga_{12}^2)+\pa_wV_2\pa_v(V_1\Ga_{11}^1+V_2\Ga_{12}^1)
+\frac{\pa_wV_1(\pa_uV_2+V_1\Ga_{11}^2+V_2\Ga_{12}^2)
-\pa_wV_2(\pa_uV_1+V_1\Ga_{11}^1+V_2\Ga_{12}^1+1)}{(\pa_uV_1+V_1\Ga_{11}^1+V_2\Ga_{12}^1+1)
(\pa_vV_2+V_1\Ga_{12}^2+V_2\Ga_{22}^2+1)-(\pa_uV_2+V_1\Ga_{11}^2+V_2\Ga_{12}^2)
(\pa_vV_1+V_1\Ga_{12}^1+V_2\Ga_{22}^1)}\{2\pa_v(V_1\Ga_{11}^1+V_2\Ga_{12}^1)
(\pa_vV_2+V_1\Ga_{12}^2+V_2\Ga_{22}^2+1)-2\pa_v(V_1\Ga_{11}^2+V_2\Ga_{12}^2)
(\pa_vV_1+V_1\Ga_{12}^1+V_2\Ga_{22}^1)+[\Ga_{12}^1(\pa_vV_1+V_1\Ga_{12}^1+V_2\Ga_{22}^1)
+\Ga_{22}^1(\pa_vV_2+V_1\Ga_{12}^2+V_2\Ga_{22}^2+1)](\pa_uV_2+V_1\Ga_{11}^2+V_2\Ga_{12}^2)
-[\Ga_{12}^2(\pa_vV_1+V_1\Ga_{12}^1+V_2\Ga_{22}^1)
+\Ga_{22}^2(\pa_vV_2+V_1\Ga_{12}^2+V_2\Ga_{22}^2+1)](\pa_uV_1+V_1\Ga_{11}^1+V_2\Ga_{12}^1+1)
+[\Ga_{12}^1(\pa_uV_1+V_1\Ga_{11}^1+V_2\Ga_{12}^1+1)
+\Ga_{22}^1(\pa_uV_2+V_1\Ga_{11}^2+V_2\Ga_{12}^2)](\pa_vV_2+V_1\Ga_{12}^2+V_2\Ga_{22}^2+1)
-[\Ga_{12}^2(\pa_uV_1+V_1\Ga_{11}^1+V_2\Ga_{12}^1+1)
+\Ga_{22}^2(\pa_uV_2+V_1\Ga_{11}^2+V_2\Ga_{12}^2)](\pa_vV_1+V_1\Ga_{12}^1+V_2\Ga_{22}^1)\}
+\frac{\pa_wV_1(\pa_vV_2+V_1\Ga_{12}^2+V_2\Ga_{22}^2+1)
-\pa_wV_2(\pa_vV_1+V_1\Ga_{12}^1+V_2\Ga_{22}^1)}{(\pa_uV_1+V_1\Ga_{11}^1+V_2\Ga_{12}^1+1)
(\pa_vV_2+V_1\Ga_{12}^2+V_2\Ga_{22}^2+1)-(\pa_uV_2+V_1\Ga_{11}^2+V_2\Ga_{12}^2)
(\pa_vV_1+V_1\Ga_{12}^1+V_2\Ga_{22}^1)}\{-2\pa_u(V_1\Ga_{12}^1+V_2\Ga_{22}^1)
(\pa_uV_2+V_1\Ga_{11}^2+V_2\Ga_{12}^2)+2\pa_u(V_1\Ga_{12}^2+V_2\Ga_{22}^2)
(\pa_uV_1+V_1\Ga_{11}^1+V_2\Ga_{12}^1+1)-[\Ga_{11}^1(\pa_vV_1+V_1\Ga_{12}^1+V_2\Ga_{22}^1)
+\Ga_{12}^1(\pa_vV_2+V_1\Ga_{12}^2+V_2\Ga_{22}^2+1)](\pa_uV_2+V_1\Ga_{11}^2+V_2\Ga_{12}^2)
+[\Ga_{11}^2(\pa_vV_1+V_1\Ga_{12}^1+V_2\Ga_{22}^1)
+\Ga_{12}^2(\pa_vV_2+V_1\Ga_{12}^2+V_2\Ga_{22}^2+1)](\pa_uV_1+V_1\Ga_{11}^1+V_2\Ga_{12}^1+1)
-[\Ga_{11}^1(\pa_uV_1+V_1\Ga_{11}^1+V_2\Ga_{12}^1+1)
+\Ga_{12}^1(\pa_uV_2+V_1\Ga_{11}^2+V_2\Ga_{12}^2)](\pa_vV_2+V_1\Ga_{12}^2+V_2\Ga_{22}^2+1)
+[\Ga_{11}^2(\pa_uV_1+V_1\Ga_{11}^1+V_2\Ga_{12}^1+1)
+\Ga_{12}^2(\pa_uV_2+V_1\Ga_{11}^2+V_2\Ga_{12}^2)](\pa_vV_1+V_1\Ga_{12}^1+V_2\Ga_{22}^1)\}
+\frac{K}{(\pa_uV_1+V_1\Ga_{11}^1+V_2\Ga_{12}^1+1)
(\pa_vV_2+V_1\Ga_{12}^2+V_2\Ga_{22}^2+1)-(\pa_uV_2+V_1\Ga_{11}^2+V_2\Ga_{12}^2)
(\pa_vV_1+V_1\Ga_{12}^1+V_2\Ga_{22}^1)}\{(\pa_wV_1V_2\\-\pa_wV_2V_1)
[(\pa_uV_1+V_1\Ga_{11}^1+V_2\Ga_{12}^1+1)\pa_ux_0
+(\pa_uV_2+V_1\Ga_{11}^2+V_2\Ga_{12}^2)\pa_vx_0]^T[(\pa_vV_1+V_1\Ga_{12}^1+V_2\Ga_{22}^1)\}\pa_ux_0
+(\pa_vV_2+V_1\Ga_{12}^2+V_2\Ga_{22}^2+1)\pa_vx_0]\\
-\frac{\pa_wV_1(\pa_uV_2+V_1\Ga_{11}^2+V_2\Ga_{12}^2)
-\pa_wV_2(\pa_uV_1+V_1\Ga_{11}^1+V_2\Ga_{12}^1+1)}{(\pa_uV_1+V_1\Ga_{11}^1+V_2\Ga_{12}^1+1)
(\pa_vV_2+V_1\Ga_{12}^2+V_2\Ga_{22}^2+1)-(\pa_uV_2+V_1\Ga_{11}^2+V_2\Ga_{12}^2)
(\pa_vV_1+V_1\Ga_{12}^1+V_2\Ga_{22}^1)}|(\pa_vV_1+V_1\Ga_{12}^1+V_2\Ga_{22}^1)\pa_ux_0
+(\pa_vV_2+V_1\Ga_{12}^2+V_2\Ga_{22}^2+1)\pa_vx_0|^2[V_1(\pa_uV_2+V_1\Ga_{11}^2+V_2\Ga_{12}^2)
-V_2(\pa_uV_1+V_1\Ga_{11}^1+V_2\Ga_{12}^1+1)]
+\frac{\pa_wV_1(\pa_vV_2+V_1\Ga_{12}^2+V_2\Ga_{22}^2+1)
-\pa_wV_2(\pa_vV_1+V_1\Ga_{12}^1+V_2\Ga_{22}^1)}{(\pa_uV_1+V_1\Ga_{11}^1+V_2\Ga_{12}^1+1)
(\pa_vV_2+V_1\Ga_{12}^2+V_2\Ga_{22}^2+1)-(\pa_uV_2+V_1\Ga_{11}^2+V_2\Ga_{12}^2)
(\pa_vV_1+V_1\Ga_{12}^1+V_2\Ga_{22}^1)}|(\pa_uV_1+V_1\Ga_{11}^1+V_2\Ga_{12}^1+1)\pa_ux_0
+(\pa_uV_2+V_1\Ga_{11}^2+V_2\Ga_{12}^2)\pa_vx_0|^2[V_1(\pa_vV_2+V_1\Ga_{12}^2+V_2\Ga_{22}^2+1)
-V_2(\pa_vV_1+V_1\Ga_{12}^1+V_2\Ga_{22}^1)]\}=0.$

Using the Gau\ss\ equations (see Eisehhart  \cite{E1} pg 155)

$Kg_{11}=\pa_v\Ga_{11}^2-\pa_u\Ga_{12}^2+\Ga_{11}^1\Ga_{12}^2+\Ga_{11}^2\Ga_{22}^2
-\Ga_{12}^1\Ga_{11}^2-{\Ga_{12}^2}^2,\\
Kg_{12}=\pa_u\Ga_{12}^1-\pa_v\Ga_{11}^1+\Ga_{12}^2\Ga_{12}^1-\Ga_{11}^2\Ga_{22}^1,\\
Kg_{12}=\pa_v\Ga_{12}^2-\pa_u\Ga_{22}^2+\Ga_{12}^1\Ga_{12}^2-\Ga_{22}^1\Ga_{11}^2,\\
Kg_{22}=\pa_u\Ga_{22}^1-\pa_v\Ga_{12}^1+\Ga_{22}^2\Ga_{12}^1+\Ga_{22}^1\Ga_{11}^1
-\Ga_{12}^2\Ga_{22}^1-{\Ga_{12}^1}^2$

the relations $C=0$ for the seventh, eighth and ninth relations of
R1, $C=0$ and $\hat C=0$ for the third relation of R2 are
identically satisfied (do not impose any condition).

The remaining equations of R1 (the third, sixth and tenth) depend
a-priori also linearly on the first derivatives of the second
fundamental form of $x_0$ and quadratically on the second
fundamental form of $x_0$; using (\ref{eq:tang}) and the Gau\ss\ equation these boil down
as before to $\pa_uN_0^TA+\pa_vN_0^TB+C=0$, where $A,B,C$ do not
depend on the second fundamental form of $x_0$.

For the third relation of R1 we have $C=-K\mathbf{m}\mathcal{\ti
U}(\pa_ux_0\times\pa_vx_0)^TN_0
+\frac{N_0^T(V\times\mathcal{U})}{\mathbf{m}N_0^T(\mathcal{U}\times\mathcal{V})}
[V^T\pa_wVKN_0^T(\pa_ux_0\times\pa_vx_0)+
\frac{(\pa_wV\times\pa_uV)^TN_0N_0^T(\pa_wV\times\pa_vx_0)}{N_0^T(\pa_wV\times
V)}-\frac{(\pa_wV\times\pa_vV)^TN_0N_0^T(\pa_wV\times\pa_ux_0)}{N_0^T(\pa_wV\times
V)}]+...$

From the relations (\ref{eq:intrdst}), applying $\pa_w$ to the
second equation of (\ref{eq:intrdst}) and using the first one we
obtain the linear system of 5 equations in
$\pa_u^2V_1,\pa_v^2V_1,\pa_w^2V_1,\pa_{vw}^2V_1,
\pa_{uw}^2V_1,\pa_{uv}^2V_1,\\
\pa_u^2V_2,\pa_v^2V_2,\pa_w^2V_2,\pa_{vw}^2V_2,\pa_{uw}^2V_2,\pa_{uv}^2V_2$:

$\pa_w^2V_1[\frac{N_0^T(\pa_ux_0\times\mathcal{U})}{N_0^T(\pa_wV\times
V)}N_0^T(V\times\pa_vx_0)-\frac{N_0^T(\pa_ux_0\times\mathcal{V})}{N_0^T(\pa_wV\times
V)}N_0^T(V\times\pa_ux_0)-\frac{N_0^T(\pa_ux_0\times
V)N_0^T(\pa_wV\times\mathcal{U})}{N_0^T(\pa_wV\times
V)^2}N_0^T(V\times\pa_vx_0)+\frac{N_0^T(\pa_ux_0\times
V)N_0^T(\pa_wV\times\mathcal{V})}{N_0^T(\pa_wV\times
V)^2}N_0^T(V\times\pa_ux_0)]+
\pa_{uw}^2V_1\frac{N_0^T(\pa_wV\times\pa_ux_0)}{N_0^T(\pa_wV\times
V)}N_0^T(V\times\pa_vx_0)\\
-\pa_{vw}^2V_1\frac{N_0^T(\pa_wV\times\pa_ux_0)}{N_0^T(\pa_wV\times
V)}N_0^T(V\times\pa_ux_0)
+\pa_{uw}^2V_2\frac{N_0^T(\pa_wV\times\pa_vx_0)}{N_0^T(\pa_wV\times
V)}N_0^T(V\times\pa_vx_0)\\
-\pa_{vw}^2V_2\frac{N_0^T(\pa_wV\times\pa_vx_0)}{N_0^T(\pa_wV\times
V)}N_0^T(V\times\pa_ux_0)
+\pa_w^2V_2[\frac{N_0^T(\pa_vx_0\times\mathcal{U})}{N_0^T(\pa_wV\times
V)}N_0^T(V\times\pa_vx_0)-\frac{N_0^T(\pa_vx_0\times\mathcal{V})}{N_0^T(\pa_wV\times
V)}N_0^T(V\times\pa_ux_0)-\frac{N_0^T(\pa_vx_0\times
V)N_0^T(\pa_wV\times\mathcal{U})}{N_0^T(\pa_wV\times
V)^2}N_0^T(V\times\pa_vx_0)+\frac{N_0^T(\pa_vx_0\times
V)N_0^T(\pa_wV\times\mathcal{V})}{N_0^T(\pa_wV\times
V)^2}N_0^T(V\times\pa_ux_0)]=\\
\frac{N_0^T(\pa_wV\times\pa_uV)}{N_0^T(\pa_wV\times
V)}N_0^T(\pa_wV\times\pa_vx_0)-\frac{N_0^T(\pa_wV\times\pa_vV)}{N_0^T(\pa_wV\times
V)}N_0^T(\pa_wV\times\pa_ux_0)\\-\frac{N_0^T(\pa_wV\times[\pa_wV_1(\Ga_{11}^1\pa_ux_0+
\Ga_{11}^2\pa_vx_0)+\pa_wV_2(\Ga_{12}^1\pa_ux_0+\Ga_{12}^2\pa_vx_0)])}{N_0^T(\pa_wV\times
V)}N_0^T(V\times\pa_vx_0)\\+\frac{N_0^T(\pa_wV\times[\pa_wV_1(\Ga_{12}^1\pa_ux_0+
\Ga_{12}^2\pa_vx_0)+\pa_wV_2(\Ga_{22}^1\pa_ux_0+\Ga_{22}^2\pa_vx_0)])}{N_0^T(\pa_wV\times
V)}N_0^T(V\times\pa_ux_0),$

$-\pa_u^2V_1\frac{V_1\pa_wV_2}{2K(\pa_wV_1V_2-\pa_wV_2V_1)}
+\pa_{uw}^2V_1[\frac{V_1(\pa_uV_2+V_1\Ga_{11}^2+V_2\Ga_{12}^2)
+V_2(\pa_vV_2+V_1\Ga_{12}^2+V_2\Ga_{22}^2+1)}{2K(\pa_wV_1V_2-\pa_wV_2V_1)}\\
-V_2\frac{V_1[\pa_wV_1(\pa_uV_2+V_1\Ga_{11}^2+V_2\Ga_{12}^2)
-\pa_wV_2(\pa_uV_1+V_1\Ga_{11}^1+V_2\Ga_{12}^1+1)]}{2K(\pa_wV_1V_2-\pa_wV_2V_1)^2}\\
-V_2\frac{V_2[\pa_wV_1(\pa_vV_2+V_1\Ga_{12}^2+V_2\Ga_{22}^2+1)
-\pa_wV_2(\pa_vV_1+V_1\Ga_{12}^1+V_2\Ga_{22}^1)]}{2K(\pa_wV_1V_2-\pa_wV_2V_1)^2}]
-\pa_{uv}V_1\frac{V_2\pa_wV_2}{2K(\pa_wV_1V_2-\pa_wV_2V_1)}\\
+\pa_u^2V_2\frac{V_1\pa_wV_1}{2K(\pa_wV_1V_2-\pa_wV_2V_1)}
-\pa_{uw}^2V_2[\frac{V_1(\pa_uV_1+V_1\Ga_{11}^1+V_2\Ga_{12}^1+1)
+V_2(\pa_vV_1+V_1\Ga_{12}^1+V_2\Ga_{22}^1)}{2K(\pa_wV_1V_2-\pa_wV_2V_1)}\\
-V_1\frac{V_1[\pa_wV_1(\pa_uV_2+V_1\Ga_{11}^2+V_2\Ga_{12}^2)
-\pa_wV_2(\pa_uV_1+V_1\Ga_{11}^1+V_2\Ga_{12}^1+1)]}{2K(\pa_wV_1V_2-\pa_wV_2V_1)^2}\\
-V_1\frac{V_2[\pa_wV_1(\pa_vV_2+V_1\Ga_{12}^2+V_2\Ga_{22}^2+1)
-\pa_wV_2(\pa_vV_1+V_1\Ga_{12}^1+V_2\Ga_{22}^1)]}{2K(\pa_wV_1V_2-\pa_wV_2V_1)^2}]
+\pa_{uv}^2V_2\frac{V_2\pa_wV_1}{2K(\pa_wV_1V_2-\pa_wV_2V_1)}=\\
-\frac{\pa_uV_1[\pa_wV_1(\pa_uV_2+V_1\Ga_{11}^2+V_2\Ga_{12}^2)
-\pa_wV_2(\pa_uV_1+V_1\Ga_{11}^1+V_2\Ga_{12}^1+1)]}{2K(\pa_wV_1V_2-\pa_wV_2V_1)}\\
-\frac{\pa_uV_2[\pa_wV_1(\pa_vV_2+V_1\Ga_{12}^2+V_2\Ga_{22}^2+1)
-\pa_wV_2(\pa_vV_1+V_1\Ga_{12}^1+V_2\Ga_{22}^1)]}{2K(\pa_wV_1V_2-\pa_wV_2V_1)}\\
-\frac{V_1[\pa_wV_1\pa_u(V_1\Ga_{11}^2+V_2\Ga_{12}^2)
-\pa_wV_2\pa_u(V_1\Ga_{11}^1+V_2\Ga_{12}^1)]}{2K(\pa_wV_1V_2-\pa_wV_2V_1)}
-\frac{V_2[\pa_wV_1\pa_u(V_1\Ga_{12}^2+V_2\Ga_{22}^2)
-\pa_wV_2\pa_u(V_1\Ga_{12}^1+V_2\Ga_{22}^1)]}{2K(\pa_wV_1V_2-\pa_wV_2V_1)}\\
+[\frac{V_1[\pa_wV_1(\pa_uV_2+V_1\Ga_{11}^2+V_2\Ga_{12}^2)
-\pa_wV_2(\pa_uV_1+V_1\Ga_{11}^1+V_2\Ga_{12}^1+1)]}{2K(\pa_wV_1V_2-\pa_wV_2V_1)^2}\\
+\frac{V_2[\pa_wV_1(\pa_vV_2+V_1\Ga_{12}^2+V_2\Ga_{22}^2+1)
-\pa_wV_2(\pa_vV_1+V_1\Ga_{12}^1+V_2\Ga_{22}^1)]}{2K(\pa_wV_1V_2-\pa_wV_2V_1)^2}]
(\pa_wV_1\pa_uV_2-\pa_wV_2\pa_uV_1)\\
+[\frac{V_1[\pa_wV_1(\pa_uV_2+V_1\Ga_{11}^2+V_2\Ga_{12}^2)
-\pa_wV_2(\pa_uV_1+V_1\Ga_{11}^1+V_2\Ga_{12}^1+1)]}{2K^2(\pa_wV_1V_2-\pa_wV_2V_1)}\\
+\frac{V_2[\pa_wV_1(\pa_vV_2+V_1\Ga_{12}^2+V_2\Ga_{22}^2+1)
-\pa_wV_2(\pa_vV_1+V_1\Ga_{12}^1+V_2\Ga_{22}^1)]}{2K^2(\pa_wV_1V_2-\pa_wV_2V_1)}]\pa_uK
-V^T\pa_uV-(2\frac{[\pa_wV\times dV]^T\wedge(\pa_wV\times
dx_0)}{K(\pa_wV\times V)^T(dx_0\times\wedge dx_0)} +V^T\pa_wV)
\frac{N_0^T(V\times\mathcal{U})}{N_0^T(\pa_wV\times V)}+
(2\frac{[\pa_wV\times d(V+x_0)]^T\wedge(V\times
dx_0)}{K(\pa_wV\times V)^T(dx_0\times\wedge dx_0)}
+|V|^2)\frac{N_0^T(\pa_wV\times\pa_uV)}{N_0^T(\pa_wV\times V)},$

$-\pa_v^2V_1\frac{V_2\pa_wV_2}{2K(\pa_wV_1V_2-\pa_wV_2V_1)}
+\pa_{vw}^2V_1[\frac{V_1(\pa_uV_2+V_1\Ga_{11}^2+V_2\Ga_{12}^2)
+V_2(\pa_vV_2+V_1\Ga_{12}^2+V_2\Ga_{12}^2+1)}{2K(\pa_wV_1V_2-\pa_wV_2V_1)}\\
-V_2\frac{V_1[\pa_wV_1(\pa_uV_2+V_1\Ga_{11}^2+V_2\Ga_{12}^2)
-\pa_wV_2(\pa_uV_1+V_1\Ga_{11}^1+V_2\Ga_{12}^1+1)]}{2K(\pa_wV_1V_2-\pa_wV_2V_1)^2}\\
-V_2\frac{V_2[\pa_wV_1(\pa_vV_2+V_1\Ga_{12}^2+V_2\Ga_{22}^2+1)
-\pa_wV_2(\pa_vV_1+V_1\Ga_{12}^1+V_2\Ga_{22}^1)]}{2K(\pa_wV_1V_2-\pa_wV_2V_1)^2}]
-\pa_{uv}^2V_1\frac{V_1\pa_wV_2}{2K(\pa_wV_1V_2-\pa_wV_2V_1)}\\
+\pa_v^2V_2\frac{V_2\pa_wV_1}{2K(\pa_wV_1V_2-\pa_wV_2V_1)}
-\pa_{vw}^2V_2[\frac{V_1(\pa_uV_1+V_1\Ga_{11}^1+V_2\Ga_{12}^1+1)
+V_2(\pa_vV_1+V_1\Ga_{12}^1+V_2\Ga_{22}^1)}{2K(\pa_wV_1V_2-\pa_wV_2V_1)}\\
-V_1\frac{V_1[\pa_wV_1(\pa_uV_2+V_1\Ga_{11}^2+V_2\Ga_{12}^2)
-\pa_wV_2(\pa_uV_1+V_1\Ga_{11}^1+V_2\Ga_{12}^1+1)]}{2K(\pa_wV_1V_2-\pa_wV_2V_1)^2}\\
-V_1\frac{V_2[\pa_wV_1(\pa_vV_2+V_1\Ga_{12}^2+V_2\Ga_{22}^2+1)
-\pa_wV_2(\pa_vV_1+V_1\Ga_{12}^1+V_2\Ga_{22}^1)]}{2K(\pa_wV_1V_2-\pa_wV_2V_1)^2}]
+\pa_{uv}^2V_2\frac{V_1\pa_wV_1}{2K(\pa_wV_1V_2-\pa_wV_2V_1)}=\\
-\frac{\pa_vV_1[\pa_wV_1(\pa_uV_2+V_1\Ga_{11}^2+V_2\Ga_{12}^2)
-\pa_wV_2(\pa_uV_1+V_1\Ga_{11}^1+V_2\Ga_{12}^1+1)]}{2K(\pa_wV_1V_2-\pa_wV_2V_1)}\\
-\frac{\pa_vV_2[\pa_wV_1(\pa_vV_2+V_1\Ga_{12}^2+V_2\Ga_{22}^2+1)
-\pa_wV_2(\pa_vV_1+V_1\Ga_{12}^1+V_2\Ga_{22}^1)]}{2K(\pa_wV_1V_2-\pa_wV_2V_1)}\\
-\frac{V_1[\pa_wV_1\pa_v(V_1\Ga_{11}^2+V_2\Ga_{12}^2)
-\pa_wV_2\pa_v(V_1\Ga_{11}^1+V_2\Ga_{12}^1)]}{2K(\pa_wV_1V_2-\pa_wV_2V_1)}
-\frac{V_2[\pa_wV_1\pa_v(V_1\Ga_{12}^2+V_2\Ga_{22}^2)
-\pa_wV_2\pa_v(V_1\Ga_{12}^1+V_2\Ga_{22}^1)]}{2K(\pa_wV_1V_2-\pa_wV_2V_1)}\\
+[\frac{V_1[\pa_wV_1(\pa_uV_2+V_1\Ga_{11}^2+V_2\Ga_{12}^2)
-\pa_wV_2(\pa_uV_1+V_1\Ga_{11}^1+V_2\Ga_{12}^1+1)]}{2K(\pa_wV_1V_2-\pa_wV_2V_1)^2}\\
+\frac{V_2[\pa_wV_1(\pa_vV_2+V_1\Ga_{12}^2+V_2\Ga_{22}^2+1)
-\pa_wV_2(\pa_vV_1+V_1\Ga_{12}^1+V_2\Ga_{22}^1)]}{2K(\pa_wV_1V_2-\pa_wV_2V_1)^2}]
(\pa_wV_1\pa_vV_2-\pa_wV_2\pa_vV_1)\\
+[\frac{V_1[\pa_wV_1(\pa_uV_2+V_1\Ga_{11}^2+V_2\Ga_{12}^2)
-\pa_wV_2(\pa_uV_1+V_1\Ga_{11}^1+V_2\Ga_{12}^1+1)]}{2K^2(\pa_wV_1V_2-\pa_wV_2V_1)}\\
+\frac{V_2[\pa_wV_1(\pa_vV_2+V_1\Ga_{12}^2+V_2\Ga_{22}^2+1)
-\pa_wV_2(\pa_vV_1+V_1\Ga_{12}^1+V_2\Ga_{22}^1)]}{2K^2(\pa_wV_1V_2-\pa_wV_2V_1)}]\pa_vK
-V^T\pa_vV-(2\frac{[\pa_wV\times dV]^T\wedge(\pa_wV\times
dx_0)}{K(\pa_wV\times V)^T(dx_0\times\wedge dx_0)} +V^T\pa_wV)
\frac{N_0^T(V\times\mathcal{V})}{N_0^T(\pa_wV\times V)}+
(2\frac{[\pa_wV\times d(V+x_0)]^T\wedge(V\times
dx_0)}{K(\pa_wV\times V)^T(dx_0\times\wedge dx_0)}
+|V|^2)\frac{N_0^T(\pa_wV\times\pa_vV)}{N_0^T(\pa_wV\times V)},$

$-\pa_u^2V_1\frac{\pa_wV_1\pa_wV_2}{K(\pa_wV_1V_2-\pa_wV_2V_1)}
+\pa_w^2V_1[[2\frac{N_0^T(\pa_ux_0\times dV)\wedge
N_0^T(\pa_wV\times dx_0)}{KN_0^T(\pa_wV\times
V)N_0^T(dx_0\times\wedge dx_0)}+2\frac{(\pa_wV\times
dV)^TN_0\wedge N_0^T(\pa_ux_0\times dx_0)}{KN_0^T(\pa_wV\times
V)N_0^T(dx_0\times\wedge dx_0)}-2\frac{N_0^T(\pa_ux_0\times
V)}{KN_0^T(\pa_wV\times V)^2}\frac{(\pa_wV\times dV)^TN_0\wedge
N_0^T(\pa_wV\times dx_0)}{N_0^T(dx_0\times\wedge
dx_0)}+V^T\pa_ux_0]\frac{N_0^T(V\times\mathcal{U})}{N_0^T(\pa_wV\times
V)}-[2\frac{(\pa_wV\times dV)^TN_0\wedge N_0^T(\pa_wV\times
dx_0)}{KN_0^T(\pa_wV\times V)N_0^T(dx_0\times\wedge
dx_0)}+V^T\pa_wV]\frac{N_0^T(\pa_ux_0\times
V)N_0^T(V\times\mathcal{U})}{N_0^T(\pa_wV\times V)^2}
-[2\frac{[\pa_wV\times d(V+x_0)]^TN_0\wedge N_0^T(V\times
dx_0)}{KN_0^T(\pa_wV\times V)N_0^T(dx_0\times\wedge
dx_0)}+|V|^2][\frac{N_0^T(\pa_ux_0\times\pa_uV)}{N_0^T(\pa_wV\times
V)}\\-\frac{N_0^T(\pa_ux_0\times
V)N_0^T(\pa_wV\times\pa_uV)}{N_0^T(\pa_wV\times V)^2}]]
-\pa_{vw}^2V_1\frac{\pa_wV_2^2}{K(\pa_wV_1V_2-\pa_wV_2V_1)}
\frac{N_0^T(V\times\mathcal{U})}{N_0^T(\pa_wV\times V)}\\
+\pa_{uw}^2V_1[\frac{2\pa_wV_1(\pa_uV_2+V_1\Ga_{11}^2+V_2\Ga_{12}^2)
-\pa_wV_2(\pa_uV_1+V_1\Ga_{11}^1+V_2\Ga_{12}^1)}{K(\pa_wV_1V_2-\pa_wV_2V_1)}
+\frac{\pa_wV_2(\pa_vV_2+V_1\Ga_{12}^2+V_2\Ga_{22}^2)}{K(\pa_wV_1V_2-\pa_wV_2V_1)}\\
-V_2\pa_wV_1\frac{\pa_wV_1(\pa_uV_2+V_1\Ga_{11}^2+V_2\Ga_{12}^2)
-\pa_wV_2(\pa_uV_1+V_1\Ga_{11}^1+V_2\Ga_{12}^1)}{K(\pa_wV_1V_2-\pa_wV_2V_1)^2}\\
-V_2\pa_wV_2\frac{\pa_wV_1(\pa_vV_2+V_1\Ga_{12}^2+V_2\Ga_{22}^2)
-\pa_wV_2(\pa_vV_1+V_1\Ga_{12}^1+V_2\Ga_{22}^1)}{K(\pa_wV_1V_2-\pa_wV_2V_1)^2}
+V^T\pa_ux_0-\frac{\pa_wV_1\pa_wV_2}{K(\pa_wV_1V_2-\pa_wV_2V_1)}
\frac{N_0^T(V\times\mathcal{U})}{N_0^T(\pa_wV\times V)}\\
+(2\frac{(\pa_wV\times dV)^TN_0\wedge N_0^T(\pa_wV\times
dx_0)}{KN_0^T(\pa_wV\times V)N_0^T(dx_0\times\wedge
dx_0)}+V^T\pa_wV)\frac{N_0^T(V\times\pa_ux_0)}{N_0^T(\pa_wV\times
V)}-(2\frac{[\pa_wV\times d(V+x_0)]^TN_0\wedge N_0^T(V\times
dx_0)}{KN_0^T(\pa_wV\times V)N_0^T(dx_0\times\wedge
dx_0)}\\+|V|^2)\frac{N_0^T(\pa_wV\times\pa_ux_0)}{N_0^T(\pa_wV\times
V)}]
-\pa_{uv}^2V_1\frac{\pa_wV_2^2}{K(\pa_wV_1V_2-\pa_wV_2V_1)}\\
+\pa_u^2V_2\frac{\pa_wV_1^2}{K(\pa_wV_1V_2-\pa_wV_2V_1)}
+\pa_w^2V_2[[2\frac{N_0^T(\pa_vx_0\times dV)\wedge
N_0^T(\pa_wV\times dx_0)}{KN_0^T(\pa_wV\times
V)N_0^T(dx_0\times\wedge dx_0)}+2\frac{(\pa_wV\times
dV)^TN_0\wedge N_0^T(\pa_vx_0\times dx_0)}{KN_0^T(\pa_wV\times
V)N_0^T(dx_0\times\wedge dx_0)}\\-2\frac{N_0^T(\pa_vx_0\times
V)}{KN_0^T(\pa_wV\times V)^2}\frac{(\pa_wV\times dV)^TN_0\wedge
N_0^T(\pa_wV\times dx_0)}{N_0^T(dx_0\times\wedge
dx_0)}+V^T\pa_vx_0]\frac{N_0^T(V\times\mathcal{U})}{N_0^T(\pa_wV\times
V)}\\-[2\frac{(\pa_wV\times dV)^TN_0\wedge N_0^T(\pa_wV\times
dx_0)}{KN_0^T(\pa_wV\times V)N_0^T(dx_0\times\wedge
dx_0)}+V^T\pa_wV]\frac{N_0^T(\pa_vx_0\times
V)N_0^T(V\times\mathcal{U})}{N_0^T(\pa_wV\times V)^2}
-[2\frac{[\pa_wV\times d(V+x_0)]^TN_0\wedge N_0^T(V\times
dx_0)}{KN_0^T(\pa_wV\times V)N_0^T(dx_0\times\wedge
dx_0)}+|V|^2][\frac{N_0^T(\pa_vx_0\times\pa_uV)}{N_0^T(\pa_wV\times
V)}-\frac{N_0^T(\pa_vx_0\times
V)N_0^T(\pa_wV\times\pa_uV)}{N_0^T(\pa_wV\times V)^2}]]
+\pa_{vw}^2V_2\frac{\pa_wV_1\pa_wV_2}{K(\pa_wV_1V_2-\pa_wV_2V_1)}
\frac{N_0^T(V\times\mathcal{U})}{N_0^T(\pa_wV\times V)}\\
+\pa_{uw}^2V_2[-\frac{\pa_wV_1(\pa_uV_1+V_1\Ga_{11}^1
+V_2\Ga_{12}^1)}{K(\pa_wV_1V_2-\pa_wV_2V_1)}
+\frac{\pa_wV_1(\pa_vV_2+V_1\Ga_{12}^2+V_2\Ga_{22}^2)
-2\pa_wV_2(\pa_vV_1+V_1\Ga_{12}^1+V_2\Ga_{22}^1)}{K(\pa_wV_1V_2-\pa_wV_2V_1)}\\
+V_1\pa_wV_1\frac{\pa_wV_1(\pa_uV_2+V_1\Ga_{11}^2+V_2\Ga_{12}^2)
-\pa_wV_2(\pa_uV_1+V_1\Ga_{11}^1+V_2\Ga_{12}^1)}{K(\pa_wV_1V_2-\pa_wV_2V_1)^2}\\
+V_1\pa_wV_2\frac{\pa_wV_1(\pa_vV_2+V_1\Ga_{12}^2+V_2\Ga_{22}^2)
-\pa_wV_2(\pa_vV_1+V_1\Ga_{12}^1+V_2\Ga_{22}^1)}{K(\pa_wV_1V_2-\pa_wV_2V_1)^2}
+V^T\pa_vx_0+\frac{\pa_wV_1^2}{K(\pa_wV_1V_2-\pa_wV_2V_1)}
\frac{N_0^T(V\times\mathcal{U})}{N_0^T(\pa_wV\times V)}\\
+(2\frac{(\pa_wV\times dV)^TN_0\wedge N_0^T(\pa_wV\times
dx_0)}{KN_0^T(\pa_wV\times V)N_0^T(dx_0\times\wedge
dx_0)}+V^T\pa_wV)\frac{N_0^T(V\times\pa_vx_0)}{N_0^T(\pa_wV\times
V)}-(2\frac{[\pa_wV\times d(V+x_0)]^TN_0\wedge N_0^T(V\times
dx_0)}{KN_0^T(\pa_wV\times V)N_0^T(dx_0\times\wedge
dx_0)}\\+|V|^2)\frac{N_0^T(\pa_wV\times\pa_vx_0)}{N_0^T(\pa_wV\times
V)}]
+\pa_{uv}^2V_2\frac{\pa_wV_1\pa_wV_2}{K(\pa_wV_1V_2-\pa_wV_2V_1)}=\\
-\pa_wV_1\frac{\pa_wV_1\pa_u(V_1\Ga_{11}^2+V_2\Ga_{12}^2)
-\pa_wV_2\pa_u(V_1\Ga_{11}^1+V_2\Ga_{12}^1)}{K(\pa_wV_1V_2-\pa_wV_2V_1)}
-\pa_wV_2\frac{\pa_wV_1\pa_u(V_1\Ga_{12}^2+V_2\Ga_{22}^2)
-\pa_wV_2\pa_u(V_1\Ga_{12}^1+V_2\Ga_{22}^1)}{K(\pa_wV_1V_2-\pa_wV_2V_1)}\\
+[\pa_wV_1\frac{\pa_wV_1(\pa_uV_2+V_1\Ga_{11}^2+V_2\Ga_{12}^2)
-\pa_wV_2(\pa_uV_1+V_1\Ga_{11}^1+V_2\Ga_{12}^1)}{K(\pa_wV_1V_2-\pa_wV_2V_1)^2}\\
+\pa_wV_2\frac{\pa_wV_1(\pa_vV_2+V_1\Ga_{12}^2+V_2\Ga_{22}^2)
-\pa_wV_2(\pa_vV_1+V_1\Ga_{12}^1+V_2\Ga_{22}^1)}{K(\pa_wV_1V_2-\pa_wV_2V_1)^2}]
(\pa_wV_1\pa_uV_2-\pa_wV_2\pa_uV_1)\\
+[\pa_wV_1\frac{\pa_wV_1(\pa_uV_2+V_1\Ga_{11}^2+V_2\Ga_{12}^2)
-\pa_wV_2(\pa_uV_1+V_1\Ga_{11}^1+V_2\Ga_{12}^1)}{K^2(\pa_wV_1V_2-\pa_wV_2V_1)}\\
+\pa_wV_2\frac{\pa_wV_1(\pa_vV_2+V_1\Ga_{12}^2+V_2\Ga_{22}^2)
-\pa_wV_2(\pa_vV_1+V_1\Ga_{12}^1+V_2\Ga_{22}^1)}{K^2(\pa_wV_1V_2-\pa_wV_2V_1)}]\pa_uK
-\pa_uV^T\pa_wV-\pa_wV_1V^T(\Ga_{11}^1\pa_ux_0\\+\Ga_{11}^2\pa_vx_0)
-\pa_wV_2V^T(\Ga_{12}^1\pa_ux_0+\Ga_{12}^2\pa_vx_0)
-\pa_wV^T\pa_wV\frac{N_0^T(V\times\mathcal{U})}{N_0^T(\pa_wV\times
V)}-(2\frac{(\pa_wV\times dV)^TN_0\wedge N_0^T(\pa_wV\times
dx_0)}{KN_0^T(\pa_wV\times V)N_0^T(dx_0\times\wedge
dx_0)}+V^T\pa_wV)\frac{N_0^T(\pa_wV\times\mathcal{U})}{N_0^T(\pa_wV\times
V)}+(4\frac{N_0^T(\pa_wV\times dV)\wedge N_0^T(\pa_wV\times
dx_0)}{KN_0^T(\pa_wV\times V)N_0^T(dx_0\times\wedge
dx_0)}+2V^T\pa_wV)\frac{N_0^T(\pa_wV\times\pa_uV)}{N_0^T(\pa_wV\times
V)}\\-\pa_wV_1\frac{\pa_wV_1(\pa_wV_1\Ga_{11}^2+\pa_wV_2\Ga_{12}^2)
-\pa_wV_2(\pa_wV_1\Ga_{11}^1+\pa_wV_2\Ga_{12}^1)}{K(\pa_wV_1V_2-\pa_wV_2V_1)}
\frac{N_0^T(V\times\mathcal{U})}{N_0^T(\pa_wV\times V)}\\
-\pa_wV_2\frac{\pa_wV_1(\pa_wV_1\Ga_{12}^2+\pa_wV_2\Ga_{22}^2)
-\pa_wV_2(\pa_wV_1\Ga_{12}^1+\pa_wV_2\Ga_{22}^1)}{K(\pa_wV_1V_2-\pa_wV_2V_1)}
\frac{N_0^T(V\times\mathcal{U})}{N_0^T(\pa_wV\times V)}
-(2\frac{(\pa_wV\times dV)^TN_0\wedge N_0^T(\pa_wV\times
dx_0)}{KN_0^T(\pa_wV\times V)N_0^T(dx_0\times\wedge
dx_0)}\\+V^T\pa_wV)\frac{N_0^T(V\times[\pa_wV_1(\Ga_{11}^1\pa_ux_0+\Ga_{11}^2\pa_vx_0)
+\pa_wV_2(\Ga_{12}^1\pa_ux_0+\Ga_{12}^2\pa_vx_0)])}{N_0^T(\pa_wV\times
V)}+(2\frac{[\pa_wV\times d(V+x_0)]^TN_0\wedge N_0^T(V\times
dx_0)}{KN_0^T(\pa_wV\times V)N_0^T(dx_0\times\wedge
dx_0)}+|V|^2)\frac{N_0^T(\pa_wV\times[\pa_wV_1(\Ga_{11}^1\pa_ux_0+\Ga_{11}^2\pa_vx_0)
+\pa_wV_2(\Ga_{12}^1\pa_ux_0+\Ga_{12}^2\pa_vx_0)])}{N_0^T(\pa_wV\times
V)},$

$-\pa_v^2V_1\frac{\pa_wV_2^2}{K(\pa_wV_1V_2-\pa_wV_2V_1)}
+\pa_w^2V_1[[2\frac{N_0^T(\pa_ux_0\times dV)\wedge
N_0^T(\pa_wV\times dx_0)}{KN_0^T(\pa_wV\times
V)N_0^T(dx_0\times\wedge dx_0)}+2\frac{(\pa_wV\times
dV)^TN_0\wedge N_0^T(\pa_ux_0\times dx_0)}{KN_0^T(\pa_wV\times
V)N_0^T(dx_0\times\wedge dx_0)}-2\frac{N_0^T(\pa_ux_0\times
V)}{KN_0^T(\pa_wV\times V)^2}\frac{(\pa_wV\times dV)^TN_0\wedge
N_0^T(\pa_wV\times dx_0)}{N_0^T(dx_0\times\wedge
dx_0)}+V^T\pa_ux_0]\frac{N_0^T(V\times\mathcal{V})}{N_0^T(\pa_wV\times
V)}-[2\frac{(\pa_wV\times dV)^TN_0\wedge N_0^T(\pa_wV\times
dx_0)}{KN_0^T(\pa_wV\times V)N_0^T(dx_0\times\wedge
dx_0)}+V^T\pa_wV]\frac{N_0^T(\pa_ux_0\times
V)N_0^T(V\times\mathcal{V})}{N_0^T(\pa_wV\times V)^2}
-[2\frac{[\pa_wV\times d(V+x_0)]^TN_0\wedge N_0^T(V\times
dx_0)}{KN_0^T(\pa_wV\times V)N_0^T(dx_0\times\wedge
dx_0)}+|V|^2][\frac{N_0^T(\pa_ux_0\times\pa_vV)}{N_0^T(\pa_wV\times
V)}\\-\frac{N_0^T(\pa_ux_0\times
V)N_0^T(\pa_wV\times\pa_vV)}{N_0^T(\pa_wV\times V)^2}]]+
\pa_{vw}^2V_1[\frac{2\pa_wV_1(\pa_uV_2+V_1\Ga_{11}^2+V_2\Ga_{12}^2)
-\pa_wV_2(\pa_uV_1+V_1\Ga_{11}^1+V_2\Ga_{12}^1)}{K(\pa_wV_1V_2-\pa_wV_2V_1)}\\
+\frac{\pa_wV_2(\pa_vV_2+V_1\Ga_{12}^2+V_2\Ga_{22}^2)}{K(\pa_wV_1V_2-\pa_wV_2V_1)}
-V_2\pa_wV_1\frac{\pa_wV_1(\pa_uV_2+V_1\Ga_{11}^2+V_2\Ga_{12}^2)
-\pa_wV_2(\pa_uV_1+V_1\Ga_{11}^1+V_2\Ga_{12}^1)}{K(\pa_wV_1V_2-\pa_wV_2V_1)^2}\\
-V_2\pa_wV_2\frac{\pa_wV_1(\pa_vV_2+V_1\Ga_{12}^2+V_2\Ga_{22}^2)
-\pa_wV_2(\pa_vV_1+V_1\Ga_{12}^1+V_2\Ga_{22}^1)}{K(\pa_wV_1V_2-\pa_wV_2V_1)^2}
+V^T\pa_ux_0-\frac{\pa_wV_2^2}{K(\pa_wV_1V_2-\pa_wV_2V_1)}
\frac{N_0^T(V\times\mathcal{V})}{N_0^T(\pa_wV\times
V)}\\-[\frac{2(\pa_wV\times dV)^TN_0\wedge N_0^T(\pa_wV\times
dx_0)}{KN_0^T(\pa_wV\times V)N_0^T(dx_0\times\wedge
dx_0)}+V^T\pa_wV]\frac{V_2}{\pa_wV_1V_2-\pa_wV_2V_1}
+[\frac{2[\pa_wV\times d(V+x_0)]^TN_0\wedge N_0^T(V\times
dx_0)}{KN_0^T(\pa_wV\times V)N_0^T(dx_0\times\wedge
dx_0)}\\+|V|^2]\frac{\pa_wV_2}{\pa_wV_1V_2-\pa_wV_2V_1}]-\pa_{uw}^2V_1\frac{\pa_wV_1\pa_wV_2}{K(\pa_wV_1V_2-\pa_wV_2V_1)}
\frac{N_0^T(V\times\mathcal{V})}{N_0^T(\pa_wV\times
V)}-\pa_{uv}^2V_1\frac{\pa_wV_1\pa_wV_2}{K(\pa_wV_1V_2-\pa_wV_2V_1)}\\
+\pa_v^2V_2\frac{\pa_wV_1\pa_wV_2}{K(\pa_wV_1V_2-\pa_wV_2V_1)}
+\pa_w^2V_2[[2\frac{N_0^T(\pa_vx_0\times dV)\wedge
N_0^T(\pa_wV\times dx_0)}{KN_0^T(\pa_wV\times
V)N_0^T(dx_0\times\wedge dx_0)}+2\frac{(\pa_wV\times
dV)^TN_0\wedge N_0^T(\pa_vx_0\times dx_0)}{KN_0^T(\pa_wV\times
V)N_0^T(dx_0\times\wedge dx_0)}\\-2\frac{N_0^T(\pa_vx_0\times
V)}{KN_0^T(\pa_wV\times V)^2}\frac{(\pa_wV\times dV)^TN_0\wedge
N_0^T(\pa_wV\times dx_0)}{N_0^T(dx_0\times\wedge
dx_0)}+V^T\pa_vx_0]\frac{N_0^T(V\times\mathcal{V})}{N_0^T(\pa_wV\times
V)}\\-[2\frac{(\pa_wV\times dV)^TN_0\wedge N_0^T(\pa_wV\times
dx_0)}{KN_0^T(\pa_wV\times V)N_0^T(dx_0\times\wedge
dx_0)}+V^T\pa_wV]\frac{N_0^T(\pa_vx_0\times
V)N_0^T(V\times\mathcal{V})}{N_0^T(\pa_wV\times V)^2}
-[2\frac{[\pa_wV\times d(V+x_0)]^TN_0\wedge N_0^T(V\times
dx_0)}{KN_0^T(\pa_wV\times V)N_0^T(dx_0\times\wedge
dx_0)}+|V|^2][\frac{N_0^T(\pa_vx_0\times\pa_vV)}{N_0^T(\pa_wV\times
V)}-\frac{N_0^T(\pa_vx_0\times
V)N_0^T(\pa_wV\times\pa_vV)}{N_0^T(\pa_wV\times V)^2}]]
+\pa_{vw}^2V_2[-\frac{\pa_wV_1(\pa_uV_1+V_1\Ga_{11}^1+V_2\Ga_{12}^1)}
{K(\pa_wV_1V_2-\pa_wV_2V_1)}\\+\frac{\pa_wV_1(\pa_vV_2+V_1\Ga_{12}^2+V_2\Ga_{22}^2)
-2\pa_wV_2(\pa_vV_1+V_1\Ga_{12}^1+V_2\Ga_{22}^1)}{K(\pa_wV_1V_2-\pa_wV_2V_1)}\\
+V_1\pa_wV_1\frac{\pa_wV_1(\pa_uV_2+V_1\Ga_{11}^2+V_2\Ga_{12}^2)
-\pa_wV_2(\pa_uV_1+V_1\Ga_{11}^1+V_2\Ga_{12}^1)}{K(\pa_wV_1V_2-\pa_wV_2V_1)^2}\\
+V_1\pa_wV_2\frac{\pa_wV_1(\pa_vV_2+V_1\Ga_{12}^2+V_2\Ga_{22}^2)
-\pa_wV_2(\pa_vV_1+V_1\Ga_{12}^1+V_2\Ga_{22}^1)}{K(\pa_wV_1V_2-\pa_wV_2V_1)^2}
+V^T\pa_vx_0+\frac{\pa_wV_1\pa_wV_2}{K(\pa_wV_1V_2-\pa_wV_2V_1)}
\frac{N_0^T(V\times\mathcal{V})}{N_0^T(\pa_wV\times V)}
+[\frac{2(\pa_wV\times dV)^TN_0\wedge N_0^T(\pa_wV\times
dx_0)}{KN_0^T(\pa_wV\times V)N_0^T(dx_0\times\wedge
dx_0)}+V^T\pa_wV]\frac{V_1}{\pa_wV_1V_2-\pa_wV_2V_1}
-[\frac{2[\pa_wV\times d(V+x_0)]^TN_0\wedge N_0^T(V\times
dx_0)}{KN_0^T(\pa_wV\times V)N_0^T(dx_0\times\wedge
dx_0)}\\+|V|^2]\frac{\pa_wV_1}{\pa_wV_1V_2-\pa_wV_2V_1}
]+\pa_{uw}^2V_2\frac{\pa_wV_1^2}{K(\pa_wV_1V_2-\pa_wV_2V_1)}
\frac{N_0^T(V\times\mathcal{V})}{N_0^T(\pa_wV\times
V)}+\pa_{uv}V_2\frac{\pa_wV_1^2}{K(\pa_wV_1V_2-\pa_wV_2V_1)}\\
=[\pa_wV_1\frac{\pa_wV_1(\pa_uV_2+V_1\Ga_{11}^2+V_2\Ga_{12}^2)
-\pa_wV_2(\pa_uV_1+V_1\Ga_{11}^1+V_2\Ga_{12}^1)}{K(\pa_wV_1V_2-\pa_wV_2V_1)^2}\\
+\pa_wV_2\frac{\pa_wV_1(\pa_vV_2+V_1\Ga_{12}^2+V_2\Ga_{22}^2)
-\pa_wV_2(\pa_vV_1+V_1\Ga_{12}^1+V_2\Ga_{22}^1)}
{K(\pa_wV_1V_2-\pa_wV_2V_1)^2}](\pa_wV_1\pa_vV_2-\pa_wV_2\pa_vV_1)\\
+[\pa_wV_1\frac{\pa_wV_1(\pa_uV_2+V_1\Ga_{11}^2+V_2\Ga_{12}^2)
-\pa_wV_2(\pa_uV_1+V_1\Ga_{11}^1+V_2\Ga_{12}^1)}{K^2(\pa_wV_1V_2-\pa_wV_2V_1)}\\
+\pa_wV_2\frac{\pa_wV_1(\pa_vV_2+V_1\Ga_{12}^2+V_2\Ga_{22}^2)
-\pa_wV_2(\pa_vV_1+V_1\Ga_{12}^1+V_2\Ga_{22}^1)}
{K^2(\pa_wV_1V_2-\pa_wV_2V_1)}]\pa_vK
-\pa_wV_1\frac{\pa_wV_1\pa_v(V_1\Ga_{11}^2+V_2\Ga_{12}^2)
-\pa_wV_2\pa_v(V_1\Ga_{11}^1+V_2\Ga_{12}^1)}{K(\pa_wV_1V_2-\pa_wV_2V_1)}
-\pa_wV_2\frac{\pa_wV_1\pa_v(V_1\Ga_{12}^2+V_2\Ga_{22}^2)
-\pa_wV_2\pa_v(V_1\Ga_{12}^1+V_2\Ga_{22}^1)}
{K(\pa_wV_1V_2-\pa_wV_2V_1)}-\pa_vV^T\pa_wV
-\pa_wV_1V^T(\Ga_{12}^1\pa_ux_0+\Ga_{12}^2\pa_vx_0)
-\pa_wV_2V^T(\Ga_{22}^1\pa_ux_0+\Ga_{22}^2\pa_vx_0)
-\pa_wV_1\frac{\pa_wV_1(\pa_wV_1\Ga_{11}^2+\pa_wV_2\Ga_{12}^2)
-\pa_wV_2(\pa_wV_1\Ga_{11}^1+\pa_wV_2\Ga_{12}^1)}{K(\pa_wV_1V_2-\pa_wV_2V_1)}
\frac{N_0^T(V\times\mathcal{V})}{N_0^T(\pa_wV\times V)}\\
-\pa_wV_2\frac{\pa_wV_1(\pa_wV_1\Ga_{12}^2+\pa_wV_2\Ga_{22}^2)
-\pa_wV_2(\pa_wV_1\Ga_{12}^1+\pa_wV_2\Ga_{22}^1)}{K(\pa_wV_1V_2-\pa_wV_2V_1)}
\frac{N_0^T(V\times\mathcal{V})}{N_0^T(\pa_wV\times V)}
-[\frac{2(\pa_wV\times dV)^TN_0\wedge N_0^T(\pa_wV\times
dx_0)}{KN_0^T(\pa_wV\times V)N_0^T(dx_0\times\wedge
dx_0)}\\+V^T\pa_wV]\frac{V_1(\pa_wV_1\Ga_{12}^2+\pa_wV_2\Ga_{22}^2)
-V_2(\pa_wV_1\Ga_{12}^1+\pa_wV_2\Ga_{22}^1)}{\pa_wV_1V_2-\pa_wV_2V_1}
+[\frac{2[\pa_wV\times d(V+x_0)]^TN_0\wedge N_0^T(V\times
dx_0)}{KN_0^T(\pa_wV\times V)N_0^T(dx_0\times\wedge
dx_0)}\\+|V|^2]\frac{\pa_wV_1(\pa_wV_1\Ga_{12}^2+\pa_wV_2\Ga_{22}^2)
-\pa_wV_2(\pa_wV_1\Ga_{12}^1+\pa_wV_2\Ga_{22}^1)}{\pa_wV_1V_2-\pa_wV_2V_1}
-|\pa_wV|^2\frac{N_0^T(V\times\mathcal{V})}{N_0^T(\pa_wV\times
V)}-[2\frac{(\pa_wV\times dV)^TN_0\wedge N_0^T(\pa_wV\times
dx_0)}{KN_0^T(\pa_wV\times V)N_0^T(dx_0\times\wedge
dx_0)}+V^T\pa_wV]\frac{N_0^T(\pa_wV\times\mathcal{V})}{N_0^T(\pa_wV\times
V)}+2[2\frac{(\pa_wV\times dV)^TN_0\wedge N_0^T(\pa_wV\times
dx_0)}{KN_0^T(\pa_wV\times V)N_0^T(dx_0\times\wedge
dx_0)}+V^T\pa_wV]\frac{N_0^T(\pa_wV\times\pa_vV)}{N_0^T(\pa_wV\times
V)},$

$\mathcal{\ti
V}=\frac{\pa_wV_1(\pa_vV_2+V_1\Ga_{12}^2+V_2\Ga_{22}^2+1)
-\pa_wV_2(\pa_vV_1+V_1\Ga_{12}^1+V_2\Ga_{22}^1)}{(\pa_uV_1+V_1\Ga_{11}^1+V_2\Ga_{12}^1+1)
(\pa_vV_2+V_1\Ga_{12}^2+V_2\Ga_{22}^2+1)-(\pa_uV_2+V_1\Ga_{11}^2+V_2\Ga_{12}^2)
(\pa_vV_1+V_1\Ga_{12}^1+V_2\Ga_{22}^1)},$

$\mathcal{\ti
U}=\frac{\pa_wV_1(\pa_uV_2+V_1\Ga_{11}^2+V_2\Ga_{12}^2)
-\pa_wV_2(\pa_uV_1+V_1\Ga_{11}^1+V_2\Ga_{12}^1+1)}{(\pa_uV_1+V_1\Ga_{11}^1+V_2\Ga_{12}^1+1)
(\pa_vV_2+V_1\Ga_{12}^2+V_2\Ga_{22}^2+1)-(\pa_uV_2+V_1\Ga_{11}^2+V_2\Ga_{12}^2)
(\pa_vV_1+V_1\Ga_{12}^1+V_2\Ga_{22}^1)},$

\section*{Acknowledgements}
The research has been partially supported by the University of Bucharest.


\begin{thebibliography}{99}
\def\topset{0pt}
\def\parsep{0pt plus 5pt minus 1pt}
\def\itemsep{-0.5ex}
\small

\bibitem{B1} L. Bianchi {\it Sur la d\`{e}formation des quadriques,} Comptes rendus de
l'Acad\'{e}mie,
{\href{http://gallica.bnf.fr/ark:/12148/bpt6k3096x/f562.image.langEN}
{{\bf 142}, (1906), 562-564}}; and
{\href{http://gallica.bnf.fr/ark:/12148/bpt6k30977/f633.image.langEN}
{{\bf 143}, (1906), 633-635}}.

\bibitem{B2} L. Bianchi {\it Lezioni Di Geometria Differenziale, Teoria delle Transformazioni
delle Superficie applicabili sulle quadriche,} Vol
{\href{http://gallica.bnf.fr/ark:/12148/bpt6k99687j.capture} {{\bf
3}}}, Enrico Spoerri Libraio-Editore, Pisa (1909).

\bibitem{B3} L. Bianchi {\it Concerning Singular Transformations $B_k$ of surfaces applicable
to quadrics,} Transactions of the American Mathematical Society,
{\href{http://www.ams.org/journals/tran/1917-018-03/S0002-9947-1917-1501075-2/S0002-9947-1917-1501075-2.pdf}
{{\bf 18} (1917), 379-401}}.

\bibitem{D1} G. Darboux  {\href{http://fr.dleex.com/details/?11194}
{\it Le\c{c}ons Sur La Th\'{e}orie G\'{e}n\'{e}rale Des
Surfaces,}} Vol {\bf 1-4}, Gauthier-Villars, Paris (1894-1917).

\bibitem{D2} G. Darboux {\href{http://fr.dleex.com/details/?11198}
{\it Le\c{c}ons Sur Les Syst\`{e}mes Orthogonaux Et Les
Coordon\'{e}es Curvilignes,}} Gauthier-Villars, Paris (1910).

\bibitem{D3} I. Dinc\u{a} {\it On Bianchi's B\"{a}cklund transformation of quadrics,}
{\href{http://arxiv.org/abs/1110.5474} {arxiv:1110.5474v2}}.

\bibitem{E1} L. P. Eisenhart {\it A Treatise on the Differential
Geometry of Curves and Surfaces,} Dover Publications, Inc., New
York, New York (1909, republished 1960).

\bibitem{E2} L. P. Eisenhart {\it Transformation of Surfaces,} Princeton University Press,
Princeton, (1922).
\end{thebibliography}
\end{document}